\documentclass[10pt]{article}

\usepackage[a4paper,margin=2cm]{geometry}
\usepackage[T1]{fontenc}
\usepackage{lmodern}
\usepackage{microtype}
\usepackage[english]{babel}
\usepackage{xcolor}
\usepackage{setspace}
\usepackage{amsmath,amssymb,amsthm,mathtools,amsfonts}
\usepackage{bbm} 
\numberwithin{equation}{section}

\usepackage{graphicx}
\usepackage[export]{adjustbox} 
\usepackage{tikz}
\usepackage{subcaption}

\usepackage{enumitem}
\setlist[enumerate,1]{label=(\roman*),ref=\roman*}

\usepackage[toc,page]{appendix}

\usepackage{thmtools}
\usepackage[
    colorlinks=true,
    linkcolor=blue,
    citecolor=red,
    urlcolor=blue
]{hyperref}
\usepackage[nameinlink,capitalize]{cleveref}

\theoremstyle{plain}
\newtheorem{theorem}{Theorem}[section]
\newtheorem{conjecture}[theorem]{Conjecture}
\newtheorem{corollary}[theorem]{Corollary}
\newtheorem{lemma}[theorem]{Lemma}
\newtheorem{proposition}[theorem]{Proposition}
\newtheorem{remark}[theorem]{Remark}

\theoremstyle{definition}
\newtheorem{definition}[theorem]{Definition}

\crefname{enumi}{}{}
\Crefname{enumi}{}{}
\crefformat{enumi}{(#2#1#3)}
\Crefformat{enumi}{(#2#1#3)}

\title{On Romik's rationality conjectures for critical percolation}
\author{Melda Akyazi-Zavadski, Nathana\"el Berestycki, Sunil Chhita}
\date{November 19, 2026}

\begin{document}
\maketitle
\begin{abstract}
The completely packed $O(1)$ loop model can be viewed as an interface representation of critical bond percolation on the square lattice. We consider submatching events, which  prescribe how these interfaces pair finitely many consecutive boundary points. Romik's rationality conjectures \cite{romikConnectivityPatternsLoop2014} predict that the probabilities of these events on the semi-infinite cylinder are given by specific rational functions, which even become dyadic rationals in the half-plane limit.

We derive explicit Pfaffian formulae for the probability of $k$ nested arches on the semi-infinite cylinder of circumference $2N$, valid for $k<N$. We deduce that these probabilities are rational functions of $N$ and, by taking the limit $N\to\infty$, that the corresponding half-plane probabilities remain rational.
This establishes key aspects of Romik's conjectures. 

Our proof uses identities of Fonseca and Zinn-Justin \cite{fonsecaGround2009}, derived through the qKZ formalism, which express these probabilities as ratios of enumeration numbers for totally symmetric self-complementary plane partitions (TSSCPPs). These plane partitions admit a dimer representation, allowing us to apply Kasteleyn theory. We evaluate the resulting hypergeometric sum using Zeilberger's and Gosper's algorithms and then take the half-plane limit $N\to\infty$. Along the way, we derive a somewhat  simplified expression for the inverse Kasteleyn matrix entries (encoding correlations in TSSCPPs) of Ayyer and Chhita \cite{ayyerCorrelations2021}, which might be of independent interest. 
\end{abstract}

\tableofcontents

\section{Introduction}
\subsection{Background and motivation}
A striking feature of many exactly solvable lattice models is that non-trivial finite-size observables can sometimes be described by explicit closed-form expressions. The completely packed $O(1)$ loop model provides a particularly clear instance of this phenomenon. On the semi-infinite cylinder and in the half-plane, the interfaces induce random noncrossing matchings of boundary points. Of particular interest are submatching events, which prescribe the pairing of a fixed number of consecutive boundary points. Motivated by earlier numerical observations and conjectures, Romik conjectured a specific rational-function form for the probabilities of these events on the cylinder, as well as dyadic rationality of the corresponding half-plane probabilities; see \cite{romikConnectivityPatternsLoop2014} and references therein.

In this paper, we study submatching events consisting of nested arches. For each fixed $k\geq 1$, we derive an explicit Pfaffian representation for the probability of $k$ nested arches on a semi-infinite cylinder of circumference $2N$, for $N>k$, and show that this probability is a rational function of $N$. Taking the limit $N\to\infty$ then establishes rationality of the corresponding half-plane probabilities. Our approach exploits a connection to totally symmetric self-complementary plane partitions (TSSCPPs) and their dimer representation.

The first step in this approach is provided by a result of Fonseca and Zinn-Justin \cite{fonsecaGround2009}. Using integral formulae for a polynomial solution of the quantum Knizhnik--Zamolodchikov (qKZ) equation, they relate certain components and partial sums of the ground-state vector of the $O(1)$ loop model to enumerations of families of non-intersecting lattice paths. These path ensembles admit a natural interpretation in terms of TSSCPPs. In particular, for the nested submatching events considered here, their results express the corresponding cylindrical probabilities as ratios of TSSCPP enumeration numbers, with the numerator counting configurations whose associated path ensemble has a prescribed number of initial horizontal paths.

A further advantage of this reformulation is that TSSCPPs admit a natural dimer representation. A plane partition in a $2N\times2N\times2N$ box can be represented as a lozenge tiling of a regular hexagon, and the symmetries of a TSSCPP allow one to restrict this description to a fundamental one-twelfth sector of the hexagon with a free boundary. Equivalently, the resulting tilings can be regarded as perfect matchings of an associated planar graph. Ayyer and Chhita \cite{ayyerCorrelations2021} established the Pfaffian structure of this dimer model and derived explicit formulae for its inverse Kasteleyn matrix, with different expressions depending on the parities of the vertices. Building on these results, we derive explicit finite-sum expressions for the inverse Kasteleyn entries required to compute the cylindrical nested-arch probabilities through a Pfaffian.

The appearance of TSSCPPs in this problem is part of a broader context of connections between integrable lattice models and enumerative combinatorics. In particular, TSSCPPs have a long-standing relation to alternating sign matrices (ASMs). Mills, Robbins, and Rumsey \cite{millsSelfcomplementaryTotallySymmetric1986} conjectured that TSSCPPs in a $2N\times2N\times2N$ box are counted by the same product formula as $N\times N$ ASMs. The TSSCPP enumeration was subsequently proved by Andrews \cite{andrews1994}, while the ASM enumeration was proved by Zeilberger \cite{zeilbergerProofAlternatingSign1996}; Kuperberg later gave an alternative proof of the latter using the six-vertex model \cite{kuperbergAnotherProofAlternativesign1996}.

By contrast, the connection between ASMs and the $O(1)$ loop model is more direct. ASMs are in bijection with six-vertex configurations with domain-wall boundary conditions and, equivalently, with fully packed loop (FPL) configurations with alternating boundary conditions. The Razumov--Stroganov conjecture \cite{razumovCombinatorialNatureGroundState2004},
proved by Cantini and Sportiello \cite{cantiniProofRazumovStroganov2011}, identifies the appropriately normalised ground-state components of the periodic $O(1)$ loop model with the numbers of FPL configurations having prescribed boundary link patterns. Together with the identities of Fonseca and Zinn-Justin, this theorem yields exact enumerative relations between constrained TSSCPP configurations and FPL configurations subject to prescribed boundary connectivities.

Beyond their combinatorial interest, these connectivity probabilities also raise the question of how exact lattice results relate to conformal structures in the continuum. The connectivity states of the completely packed $O(1)$ loop model are naturally encoded by noncrossing link patterns, with the finite-size Hamiltonian acting on this space through Temperley--Lieb generators. 

In the continuum, noncrossing link patterns also index the pure partition functions of multiple Schramm--Löwner evolutions (SLEs) \cite{kytolaPurePartitionFunctions2016}. It is therefore natural to ask whether the exact boundary connectivity probabilities studied here admit an interpretation within this continuum framework. Such an interpretation is not immediate: our infinite-circumference limit is taken at fixed lattice spacing, whereas a continuum scaling limit involves sending the lattice spacing to zero. In the latter limit, the neighbouring boundary points involved in a fixed submatching would coalesce, so any connection with multiple SLEs would have to account for this coalescence.

Further motivation for studying exact connectivity probabilities and their asymptotics comes from recent results on critical site percolation on the triangular lattice. Sun, Xu, and Zhuang \cite{sun2024annulus} derived an exact formula for the probability of two disjoint open paths crossing an annulus in the continuum scaling limit. Nolin, Qian, Sun, and Zhuang \cite{nolinBackboneExponentAnnulus2025} discuss the asymptotic expansion of this probability, whose leading behaviour is governed by the backbone exponent. The authors interpret the subleading terms as suggesting a connection with an underlying
conformal field theory. Although these results concern a different observable, they illustrate how exact connectivity probabilities and their asymptotics may provide insight into continuum conformal structures.

\subsection{Main results}
We consider the completely packed $O(1)$ loop model on the semi-infinite cylinder $\mathcal{C}_N$ of circumference $2N$ and on the upper half-plane $\mathbb{H}$. In both geometries, the loop configuration consists almost surely of closed loops and paths connecting pairs of boundary points. These paths induce a random noncrossing matching of the boundary points; see \cref{sec:preliminaries}. 

For $k\geq 1$, let $F_k$ denote the event that the $2k$ consecutive boundary points labelled $1,\ldots,2k$ are paired in a fully nested manner, with point $i$ connected to point $2k+1-i$ for $i=1,\ldots,k$. We write $\mathbb{P}_{\mathcal{C}_N}(F_k)$ and $\mathbb{P}_{\mathbb{H}}(F_k)$ for the corresponding probabilities, with $N\geq k$ in the cylindrical case.

\begin{theorem}\label{thm:half-plane}
For every $k\geq 1$, the probability $\mathbb P_{\mathbb H}(F_k)$ is rational.
\end{theorem}

The key step in proving \cref{thm:half-plane} is the following stronger result for the cylindrical model.

\begin{theorem}
\label{thm:cylinder-rationality}
For every fixed $k\geq 1$, there exists a rational function $R_k\in\mathbb{Q}(N)$ such that
\begin{equation}
    \mathbb{P}_{\mathcal{C}_N}(F_k)=R_k(N),
    \qquad N>k.
\end{equation}
\end{theorem}

Our proof applies the inverse Kasteleyn formulae of Ayyer and Chhita \cite{ayyerCorrelations2021} to the TSSCPP events arising from the Fonseca--Zinn-Justin correspondence \cite{fonsecaGround2009}. This yields an explicit Pfaffian representation for the cylindrical probabilities, whose matrix size depends only on $k$. We use Gosper's and Zeilberger's algorithms to obtain a closed-form evaluation of the hypergeometric term appearing in the matrix entries, thereby establishing that these entries are rational functions of $N$. The restriction $k<N$ reflects the range of the Fonseca--Zinn-Justin correspondence.

By the convergence of finite boundary connectivity probabilities as $N\to\infty$ \cite[Lemma 1.17]{romikConnectivityPatternsLoop2014}, the half-plane probabilities are obtained as limits of the corresponding cylindrical probabilities. For fixed $k$, each entry of the matrix in our Pfaffian representation converges to a rational number as $N\to\infty$. Since the Pfaffian is a polynomial in its matrix entries with integer coefficients, the half-plane probabilities are Pfaffians of matrices with rational entries. In particular, they are rational, proving \cref{thm:half-plane}.

\cref{thm:cylinder-rationality} establishes rational-function dependence for the family of nested-arch events, as predicted by Romik's conjecture \cite[Conjecture 1.12]{romikConnectivityPatternsLoop2014}. However, the prescribed denominator and the additional properties of the numerator remain unproved. Likewise, our half-plane result establishes rationality of the corresponding probabilities, but not the dyadic rationality conjectured by Romik \cite[Conjecture 1.3]{romikConnectivityPatternsLoop2014}. The precise statements of these conjectures are recalled in \cref{subsec:rationality-conj}.

The boundary case $N=k$ is not covered by \cref{thm:cylinder-rationality}. In this case, $F_k$ specifies the fully nested matching, whose probability is $1/A_k$, where $A_k$ denotes the number of alternating-sign matrices of order $k$; see \cite{difrancescoRazumovStroganov2005} and \cite[Theorem 1.9]{romikConnectivityPatternsLoop2014}. For $k\geq 2$, we have $R_{k-1}(k)=1/A_k$; see \eqref{eq:rotational-inv}. Exact symbolic computations in Mathematica verify that $R_k(k)=R_{k-1}(k)$ for $2\leq k\leq 8$, suggesting the following conjecture.

\begin{conjecture}
\label{conj:cylinder-boundary} 
For every $k\geq 2$, the rational function $R_k$ from \cref{thm:cylinder-rationality} satisfies 
\begin{equation}     
    R_k(k)=R_{k-1}(k). 
\end{equation} 
\end{conjecture}  

This would extend the rational-function identity in \cref{thm:cylinder-rationality} to $N=k$ for $k\geq 2$, as predicted by Romik's cylindrical rationality conjecture. The case $k=1$ follows from the explicit formula in \cite[Eq. (7)]{romikConnectivityPatternsLoop2014}.

The explicit Pfaffian representation also allows us to evaluate the half-plane probabilities for small values of $k$. The first eight probabilities, obtained by exact symbolic computations in Mathematica, are listed in \cref{tab:half-plane-probabilities}. The corresponding inverse Kasteleyn matrix entries are given in \cref{appendix-A}.

 \begin{table}[htbp]
    \centering
    \renewcommand{\arraystretch}{2}
    \begin{tabular}{c@{\qquad\qquad}c}
        $k$ & $\mathbb{P}_{\mathbb{H}}(F_k)$ \\ \hline
        $1$ & $\dfrac{3}{2^3}$ \\
        $2$ & $\dfrac{59}{2^{10}}$ \\
        $3$ & $\dfrac{7737}{2^{21}}$ \\
        $4$ & $\dfrac{6803145}{2^{36}}$ \\
        $5$ & $\dfrac{40210542795}{2^{55}}$ \\
        $6$ & $\dfrac{1599678636998003}{2^{78}}$ \\
        $7$ & $\dfrac{428672479245743736177}{2^{105}}$ \\
        $8$ & $\dfrac{774167330660186496842976913}{2^{136}}$
    \end{tabular}
    \caption{Probabilities of the nested-arch events in the half-plane.}
    \label{tab:half-plane-probabilities}
\end{table}

The values for $k=1,2,3$ agree with the limits of the normalised enumeration formulae in \cite[Conjecture 8(ii)]{zuberCounting2004}, where Zuber credits the cases $k=1,2$ to Wilson. For $k=1,2$, the same limiting values are predicted by the conjectured formulae of Mitra, Nienhuis, de Gier and Batchelor \cite[Eqs. (30) and (34)]{mitraExact2004}. All three values are listed in \cite[Table 1]{romikConnectivityPatternsLoop2014}; the case $k=1$ was proved in \cite[Theorem 1.4]{romikConnectivityPatternsLoop2014}.

The denominators of these probabilities suggest additional structure beyond rationality, leading to the following conjecture.

\begin{conjecture}\label{conj:half-plane}
For every $k\geq 1$, there exists an odd integer $a_k$ such that
\begin{equation}
    \mathbb{P}_{\mathbb{H}}(F_k)= \frac{a_k}{2^{2k^2+k}}.
\end{equation}
\end{conjecture}

The conjecture has also been verified by exact symbolic computation for $k=17$. The requirement that $a_k$ be odd implies that $2^{2k^2+k}$ is the exact denominator of $\mathbb{P}_{\mathbb{H}}(F_k)$ in lowest terms. For nested-arch events, this conjecture therefore strengthens the dyadic-rationality assertion of Romik's half-plane conjecture by prescribing the exact denominator.

The remainder of the paper is organised as follows. In \cref{sec:preliminaries}, we explain the connection between critical bond percolation and the completely packed $O(1)$ loop model, introduce the relevant boundary connectivity events, and recall the two rationality conjectures of Romik relevant to our setting. We then introduce TSSCPPs and their representation as lozenge tilings and dimer configurations, before recalling the Fonseca--Zinn-Justin correspondence relating nested-arch probabilities to constrained TSSCPP enumerations.

In \cref{sec:simplified-K-inv}, we simplify the inverse Kasteleyn matrix formulae of Ayyer and Chhita, obtaining results that may be of independent interest. In \cref{sec:dimer-events}, we translate the TSSCPP path constraints into prescribed-dimer events, express their probabilities in terms of Pfaffians, and identify the hypergeometric term requiring further analysis. In \cref{sec:hypergeomtric}, we rewrite this term in a suitable form, recall the necessary background on Gosper's and Zeilberger's algorithms, and use these algorithms to obtain a closed-form expression. This allows us to prove \cref{thm:cylinder-rationality} and deduce \cref{thm:half-plane}. Finally, in \cref{sec:observations}, we formulate a conjectured denominator bound that would imply dyadic rationality of the limiting inverse Kasteleyn entries and hence, via the Pfaffian representation, of the nested-arch probabilities. We verify this bound in special cases and discuss the difficulties in extending the argument to general parameters.

\begin{remark} As we were finishing this draft we learnt from Dan Romik that he was able to obtain a full solution to the rationality conjectures with the help of AI. We will report on this progress and further results in forthcoming work \cite{ABCR}. 
\end{remark}

\subsection*{Acknowledgements}
This research was supported by the Austrian Science Fund (FWF), grant DOI 10.55776/F1002, within the Special Research Area ``Discrete random structures: enumeration and scaling limits''. NB also acknowledges the support of grant 10.55776/PAT1878824 on ``Random Conformal Fields'' from the FWF. SC acknowledges the support of ``Engineering and Physical Sciences Research Council'' (EPSRC) EP/T004290/1.

For symbolic computations, the authors used the \texttt{fastZeil} package (version 3.61) from the RISCErgoSum collection, developed by Peter Paule, Markus Schorn, and Axel Riese.

Part of this paper was written while all three authors were visiting the Erwin Schrödinger Institute on the occasion of the programme ``Discrete random structures'' (BFT26). The authors gratefully acknowledge hospitality, and, in the case of SC, financial support from the institution. We are also grateful to Dan Romik for sharing with us at the end of this project his preliminary work solving his conjecture with the help of AI, see \cite{ABCR}.

\medskip \noindent \textbf{AI statement and methodology}. The ideas for the results and the proofs themselves were entirely conceived without any use of AI tools. Different versions of ChatGPT (ranging from 5.2--6) were consulted after writing for assistance with language editing and improving the clarity and organisation of the exposition, as well as for occasional clarification of mathematical background and guidance on the use of Mathematica and the \texttt{fastZeil} package. The authors take full responsibility for the manuscript.

\subsection*{Notations and conventions} 

Throughout the paper, we write $[j]_2$ for $j \mod 2$ and denote the $n$-th Catalan number by $C_n=\frac{1}{n+1}\binom{2n}{n}$. Since conventions for binomial coefficients vary across the literature, we state ours explicitly to avoid ambiguity. For integers $a,b$, the binomial coefficient $\binom{a}{b}$ is nonzero if and only if $a \ge b >0$ or $b = 0$. In the latter case, we define $\binom{a}{0}=1$ for every $a \in \mathbb{Z}$, including negative $a$. 
Furthermore, we adopt the convention $\operatorname{sgn}(0)=0$ and take empty sums and products to be $0$ and $1$, respectively.

\section{Preliminaries}
\label{sec:preliminaries}

\subsection{Critical bond percolation and the $O(1)$ loop model}
The main object of study is the completely packed $O(1)$ loop model. It can be constructed from critical bond percolation on the square lattice $\mathbb{Z}^2$ by drawing the interfaces between open clusters and dual open clusters. We first recall bond percolation and then describe this construction.

In Bernoulli bond percolation with parameter $p$, each edge of $\mathbb{Z}^2$ is independently declared open with probability $p$ and closed with probability $1-p$. The connected components of the subgraph formed by the open edges are called \emph{open clusters}. \emph{Critical} bond percolation corresponds to $p=p_c=\frac{1}{2}$ \cite{harrisLower1960,kestenCritical1980}. This critical value separates the regime $p\leq p_c$, in which almost surely no infinite open cluster exists, from the regime $p>p_c$, in which an infinite open cluster exists almost surely.

\begin{figure}[ht]
    \centering
    \begin{subfigure}{0.31\textwidth}
        \centering
         \raisebox{0.45cm}{%
        \includegraphics[width=0.8\linewidth]{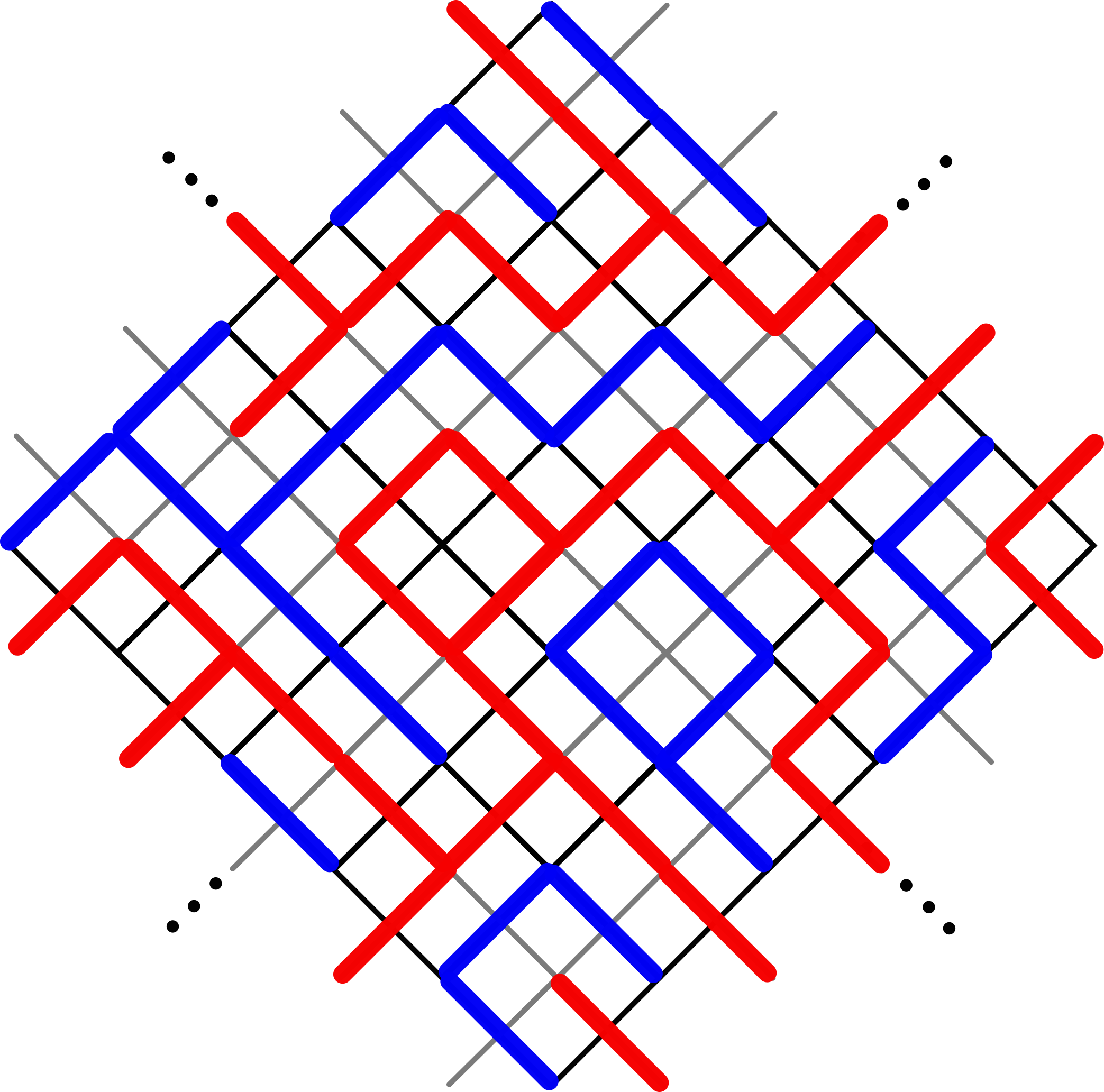}}
        \caption{Open primal (blue) and dual (red) bonds.}
        \label{fig:primal-dual}
    \end{subfigure}
    \hfill
    \begin{subfigure}{0.31\textwidth}
        \centering
        \raisebox{0.5cm}{%
        \includegraphics[width=0.8\linewidth]{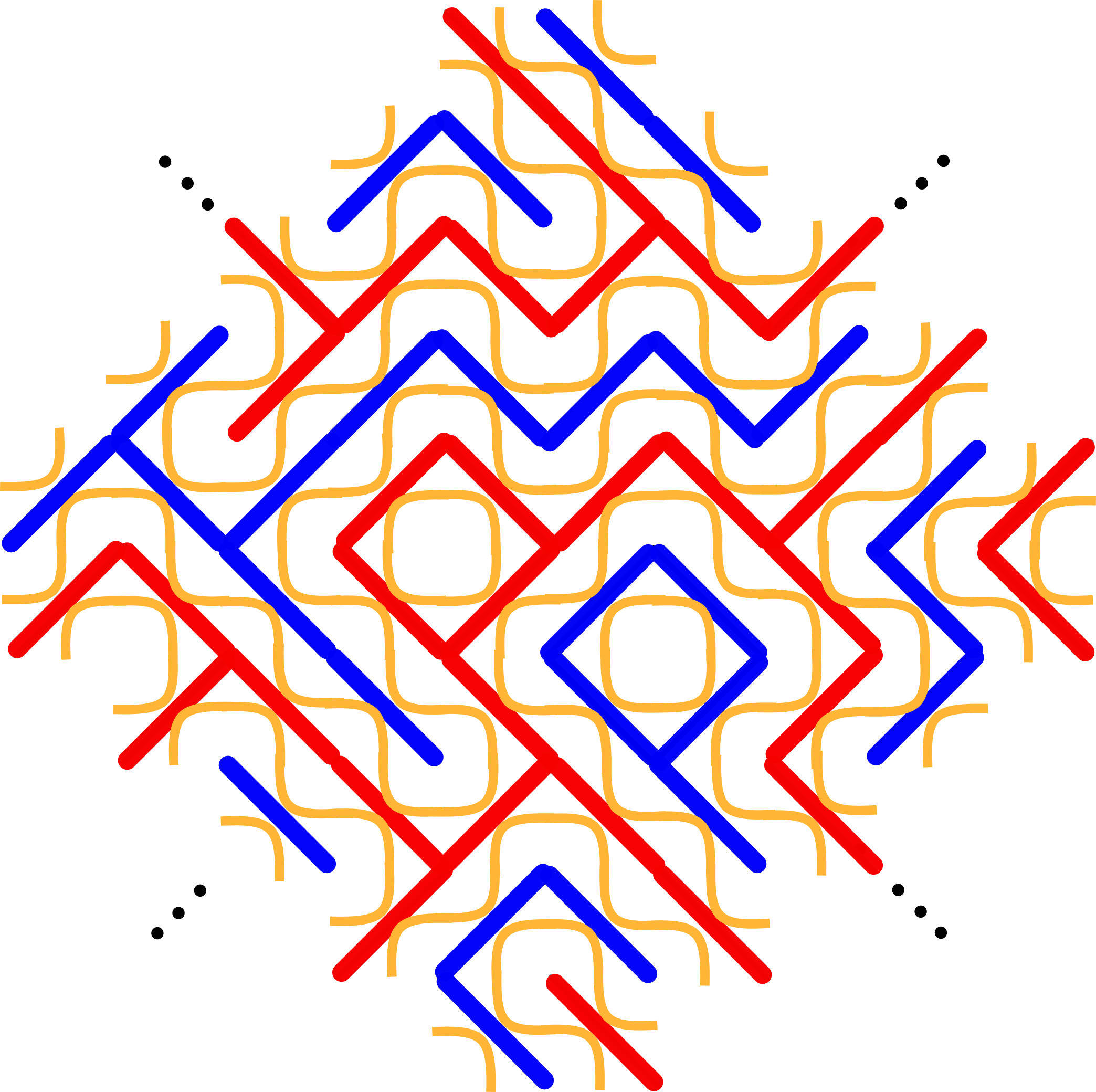}}
        \caption{Interfaces (yellow) separating primal and dual open clusters.}
        \label{fig:interfaces-with-bond}
    \end{subfigure}
    \hfill
    \begin{subfigure}{0.31\textwidth}
        \centering
        \raisebox{0.5cm}{%
        \includegraphics[width=0.8\linewidth]{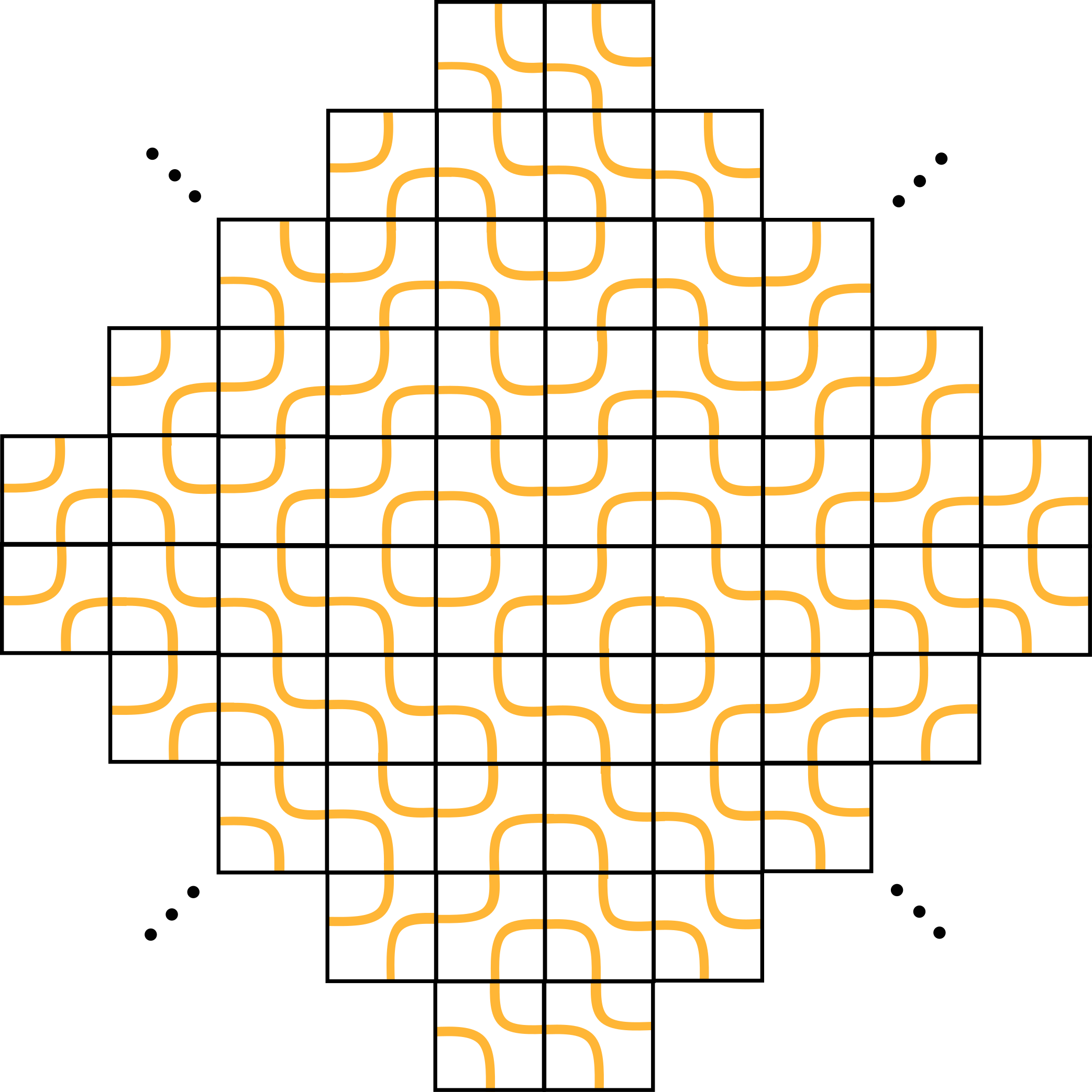}}
        \caption{The same interfaces on the plaquette lattice $M^*$ (black).}
        \label{fig:loops}
    \end{subfigure}
    \caption{Correspondence between bond percolation and the loop representation.}
    \label{fig:from-bonds-to-loops}
\end{figure}

To describe the corresponding loop representation, we introduce the primal and dual lattices, together with the medial lattice and its dual. We refer to the original square lattice as the \emph{primal lattice}; in \cref{fig:primal-dual}, it is drawn in black. The open edges in the bond-percolation configuration are shown in blue. The \emph{dual lattice} is obtained by placing one vertex in each face of the primal lattice and joining vertices in adjacent faces, so that each dual edge crosses exactly one primal edge; in \cref{fig:primal-dual}, it is drawn in grey. The corresponding dual configuration is obtained by declaring a dual edge open precisely when the primal edge it crosses is closed; we draw these open dual edges in red. Consequently, no open primal edge crosses an open dual edge, and the primal and dual configurations determine each other uniquely. Throughout, we use the rotated representation of the primal and dual lattices shown in \cref{fig:primal-dual}.

To construct the plaquette representation of the loop model, let $M$ denote the medial lattice of the primal lattice. Its vertices are the midpoints of the primal edges, with two vertices joined by an edge whenever the corresponding primal edges are adjacent sides of a common face of the primal lattice. We call the dual lattice $M^*$ of $M$ the \emph{plaquette lattice}, and call its faces \emph{plaquettes}. Both $M$ and $M^*$ are square lattices rotated by $45^\circ$ relative to the primal and dual lattices. The lattice $M^*$ is shown in black in \cref{fig:loops}. Each medial vertex lies at the intersection of a primal edge and its dual edge, exactly one of which is open. We pair the four incident medial half-edges by joining the two lying on each side of the open edge. These pairings partition the edges of $M$ into edge-disjoint trails whose consecutive edges meet at right angles. Smoothing the pairings at the medial vertices yields noncrossing interfaces separating primal open clusters from dual open clusters, which are drawn in yellow in \cref{fig:interfaces-with-bond,fig:loops}. \cref{fig:interfaces-with-bond} shows the interfaces together with the open primal and dual edges, while \cref{fig:loops} shows the same interfaces in the plaquette representation on $M^*$.

\begin{figure}[ht]
	\centering
	\includegraphics[width=0.2\textwidth]{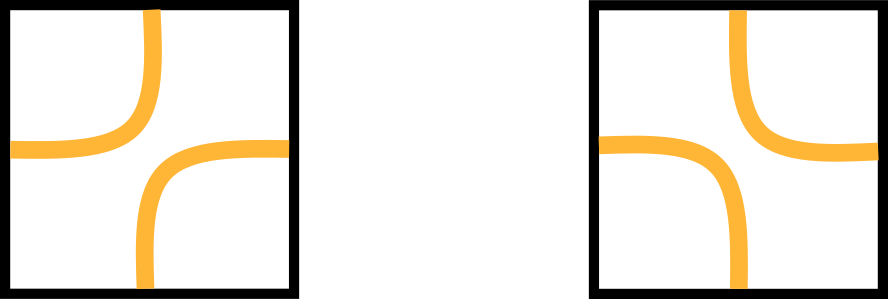}
	\caption{The two plaquette types of the completely packed $O(1)$ loop model.}
	\label{fig:plaquettes}
\end{figure}

Each face of $M^*$ contains exactly one medial vertex, so the two possible pairings at that vertex can be represented by the two plaquette types shown in \cref{fig:plaquettes}. From now on, we use this plaquette representation on $M^*$. The two plaquette types correspond to the standard local connectivities of the planar Temperley--Lieb algebra; see, for example, \cite{morin-duchesneConformal2017}. Under periodic boundary conditions, as in the cylindrical geometry considered below, they are naturally described within the affine Temperley--Lieb algebra.

At criticality, $p=\frac{1}{2}$, and each face of $M^*$ independently receives either of the two plaquette types with equal probability. In the general $O(\tau)$ loop model, each closed loop carries an additional weight $\tau$. In our case, $\tau=1$, so no additional weighting of loops occurs. Since every face of $M^*$ is occupied by one of the two plaquette types, the resulting loop configuration is completely packed. Following \cite{zinn-justinIntegrability2009}, we refer to this model as the completely packed $O(1)$ loop model. 

We denote the loop configuration by $\omega$.

\subsection{Connectivity events and Romik's rationality conjectures}
\label{subsec:rationality-conj}
To formulate Romik's rationality conjectures, we consider the loop model on the semi-infinite cylinder $\mathcal{C}_N$ or the upper half-plane $\mathbb{H}$. In the plaquette representation, we identify the underlying lattices with the graphs having vertex sets
\begin{equation*}
    V(\mathcal{C}_N) = (\mathbb{Z}/(2N\mathbb{Z}))\times\mathbb Z_{\geq 0}, \qquad
    V(\mathbb{H}) = \mathbb{Z}\times\mathbb{Z}_{\geq 0},
\end{equation*}
respectively. Both lattices have nearest-neighbour edges, with horizontal coordinates understood modulo $2N$ on $\mathcal{C}_N$. For $N=1$, the periodic construction retains the two distinct horizontal edges joining the same pair of vertices.

In both geometries, the interfaces may form closed loops, paths connecting two boundary points, or infinite paths. However, infinite paths occur with probability zero in either geometry. In the half-plane, this follows by restricting the full-plane configuration to the half-plane. Indeed, at criticality, the full-plane interfaces are almost surely finite, since neither the primal nor the dual percolation configuration contains an infinite open cluster. On $\mathcal{C}_N$, an alternating row of plaquettes prevents interfaces from crossing the row. Such rows occur independently with probability $2^{1-2N}>0$, and hence infinitely often almost surely by the second Borel--Cantelli lemma; see \cite[Lemma 1.6]{romikConnectivityPatternsLoop2014} for details. Consequently, almost surely every interface is either a closed loop or a finite path connecting two boundary points.

We regard the half-plane model as the infinite-circumference limit of the cylindrical model. This convergence relies on the almost sure finiteness of the half-plane interface paths established above. More precisely, for any fixed finite collection of prescribed connections between boundary points, the corresponding probabilities on $\mathcal{C}_N$ converge to those on $\mathbb{H}$ as $N\to\infty$; see \cite[Lemma 1.17]{romikConnectivityPatternsLoop2014} for a detailed argument.

\begin{figure}[ht]
     \centering
     \begin{subfigure}{0.46\textwidth}
         \centering
         \includegraphics[width=0.95\linewidth]{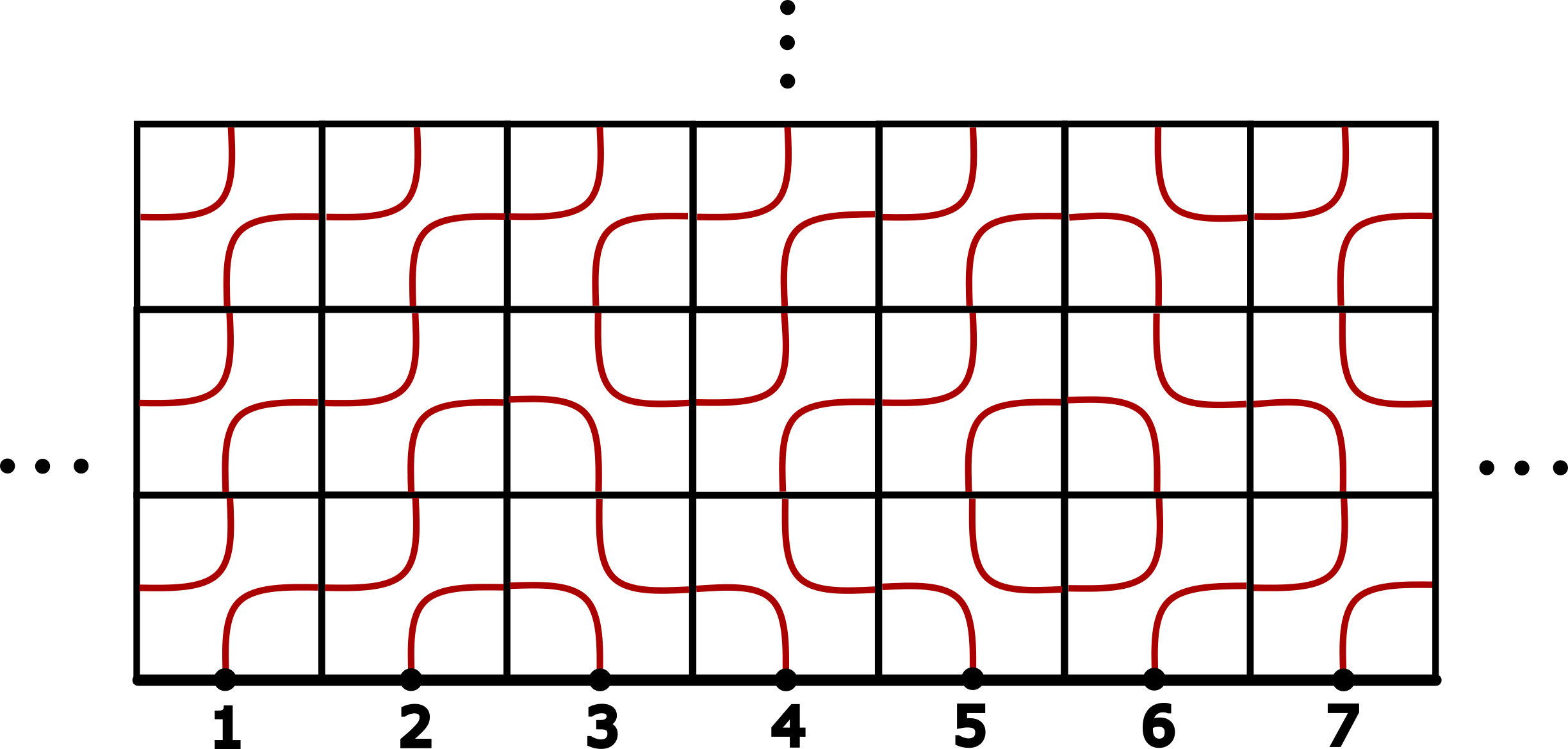}
         \caption{Portion of a loop configuration.}
         \label{fig:example-link-pattern-loop}
     \end{subfigure}
     \hfill
     \begin{subfigure}[b]{0.46\textwidth}
         \centering
         \includegraphics[width=0.8\linewidth]{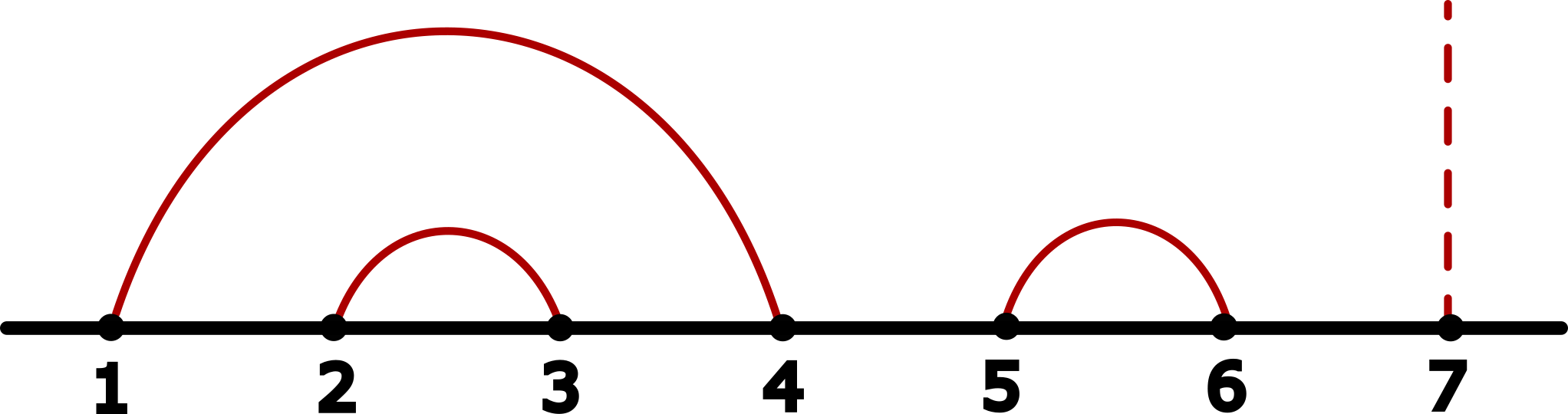}
         \caption{Corresponding boundary connectivities.}
         \label{fig:example-link-pattern}
     \end{subfigure}
        \caption{Example of a loop configuration in the half-plane and its corresponding link pattern.}
        \label{fig:example-link-pattern-full}
\end{figure}

We now focus on the connectivity of interface paths with endpoints on the boundary. In both geometries, the boundary vertices are precisely those with second coordinate zero. In the plaquette representation, the interface endpoints lie at the midpoints of the horizontal edges along this boundary; see \cref{fig:example-link-pattern-loop} for an example.

Fix $x\in\mathbb{Z}$ and set
\begin{equation*}
    b_j=\left(x+j-\frac12,0\right).
\end{equation*}
On $\mathcal{C}_N$, we take $j=1,\ldots,2N$, with horizontal coordinates understood modulo $2N$; on $\mathbb{H}$, we take $j\in\mathbb{Z}$. In either case, we label the boundary point $b_j$ by $j$.

\begin{definition}
    A \emph{noncrossing matching}, or \emph{link pattern}, of order $k$ is a perfect matching $\Pi$ of $\{1,\dots,2k\}$, which we represent by a fixed-point-free involution 
    \begin{equation*}
        \pi_{\Pi} \colon \{1,\dots, 2k\} \to \{1,\dots, 2k\},
    \end{equation*}
    such that there do not exist $a<b<c<d$ with 
    \begin{equation*}
        \pi_{\Pi}(a)=c \qquad\text{and}\qquad \pi_{\Pi}(b)=d.
    \end{equation*}
    If $\pi_{\Pi}(i)=j$, we say that $i$ is matched to $j$ in $\Pi$.
    
An \emph{infinite noncrossing matching} $\Pi$ is a perfect matching of $\mathbb{Z}$, represented by a fixed-point-free involution
\begin{equation*}
    \pi_{\Pi} \colon \mathbb{Z} \to \mathbb{Z},
\end{equation*}
satisfying the same noncrossing condition.
\end{definition}

Both finite and infinite link patterns can be represented by arranging their vertices in increasing order on a horizontal line and joining matched pairs with noncrossing arcs in the upper half-plane, as illustrated in \cref{fig:example-link-pattern}. In the finite case, an equivalent representation places the vertices in cyclic order on a circle, as is natural for the cylindrical geometry.

Each boundary point is incident to exactly one interface path. Since infinite paths are almost surely absent, these paths pair all boundary points. By construction, the pairing is noncrossing and therefore defines a link pattern of order $N$ on $\mathcal{C}_N$ and an infinite link pattern on $\mathbb{H}$. We denote the resulting connectivity pattern of a loop configuration $\omega$ by $\Pi_\omega$, and write $\pi_\omega$ for its associated involution.

On $\mathcal{C}_N$, the probability distribution of the boundary link pattern $\Pi_\omega$ is the stationary distribution of the Markov chain induced by adding independent rows of plaquettes; see \cite[Section 1.4]{romikConnectivityPatternsLoop2014}.

\begin{definition}
    Let $\Pi'$ be a link pattern of order $\ell$, and let $\Pi$ be either a link pattern of order $N\geq\ell$ or an infinite link pattern. We say that $\Pi'$ is a \emph{submatching} of $\Pi$ if
    \begin{equation*}
        \pi_{\Pi}(i)=\pi_{\Pi'}(i), \qquad i=1,\dots,2\ell.
    \end{equation*}
\end{definition}

For instance, the link pattern in \cref{fig:example-link-pattern} does not have a submatching of order $1$, since boundary points $1$ and $2$ are not matched to each other. However, it has submatchings of orders $2$ and $3$, corresponding to the first four and six boundary points, respectively. The dashed line at point $7$ indicates that its matching partner lies outside the displayed portion of the link pattern.

\begin{definition}
    For $k\geq 1$ (with $k\leq N$ on $\mathcal{C}_N$), the \emph{$k$-nested-arch event}, denoted by $F_k$, is
    \begin{equation*}
       F_k:= \left\{\omega : \pi_\omega(i)=2k+1-i \text{ for } i=1,\dots,k \right\}.
    \end{equation*}
    Equivalently, $F_k$ is the event that the link pattern of order $k$ consisting of the nested pairs $\{i,2k+1-i\}$, $i=1,\ldots,k$, is a submatching of $\Pi_\omega$.
\end{definition}

The event $F_2$ is illustrated by the link pattern in \cref{fig:example-link-pattern-full}. We first recall Romik's conjecture on the rational-function dependence of general submatching probabilities on the cylinder circumference.

\begin{conjecture}[Romik {\cite[Conjecture 1.12]{romikConnectivityPatternsLoop2014}}]
\label{conj:Romik-finite-n}
Let $\Pi'$ be a link pattern of order $k$. Then there exists an even polynomial $Q_{\Pi'}$ of degree $k(k+1)$ with dyadic rational coefficients such that, for all $N\geq k$,
\begin{equation}
    \mathbb{P}_{\mathcal{C}_N} \left( \Pi' \text{ is a submatching of } \Pi_\omega \right) = \frac{Q_{\Pi'}(N)}{\prod_{j=1}^{k}(4N^2-j^2)^{k+1-j}}.
\end{equation}
\end{conjecture}

By the convergence of finite boundary connectivity probabilities as $N\to\infty$, \cref{conj:Romik-finite-n} implies that the corresponding half-plane submatching probabilities are dyadic rationals. This is precisely the rationality assertion of Romik's half-plane conjecture.

\begin{conjecture}[Romik {\cite[Conjecture 1.3]{romikConnectivityPatternsLoop2014}}, rationality assertion]
\label{conj:dyadic-rat}
Let $\Pi'$ be a finite link pattern. Then
\begin{equation}
     \mathbb{P}_{\mathbb H} \left(\Pi' \text{ is a submatching of } \Pi_\omega\right) \in \mathbb{Z}[1/2].
\end{equation}
\end{conjecture}

We next introduce TSSCPPs and their dimer representation, which form the basis for the subsequent analysis of nested-arch probabilities.

\subsection{Totally symmetric self-complementary plane partitions}
\subsubsection{From boxed plane partitions to lozenge tilings}
\begin{figure}[ht]
     \centering
     \begin{subfigure}{0.46\textwidth}
         \centering
         \includegraphics[width=0.6\linewidth]{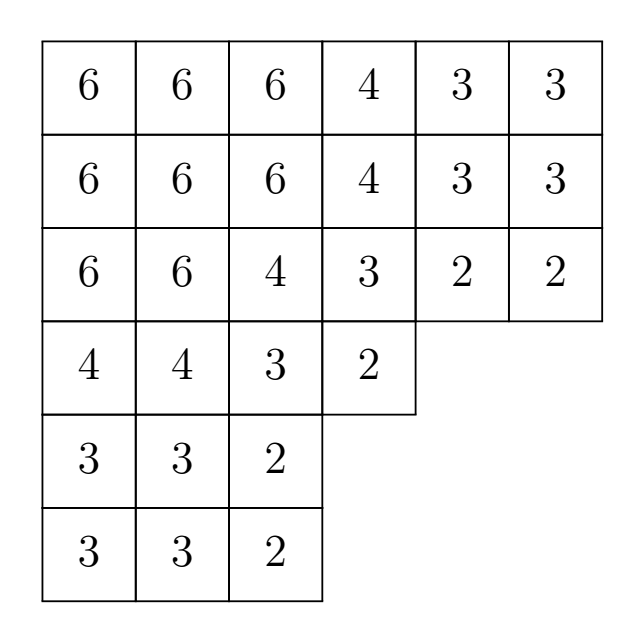}
         \caption{Two-dimensional array representation.}
         \label{fig:array}
     \end{subfigure}
     \hfill
     \begin{subfigure}[b]{0.46\textwidth}
         \centering
         \includegraphics[width=0.55\linewidth]{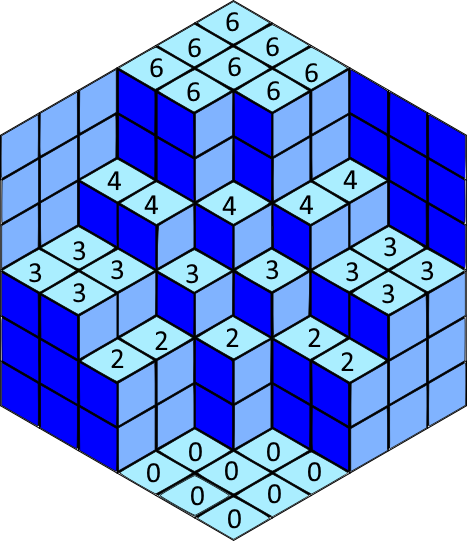}
         \caption{Unit-cube representation with indicated heights.}
         \label{fig:cubes}
     \end{subfigure}
        \caption{A plane partition shown in two equivalent representations. This example is a TSSCPP of order $3$.}
        \label{fig:example-TSSCPP}
\end{figure}

A \emph{plane partition} is a two-dimensional array of non-negative integers with finitely many non-zero entries that is weakly decreasing along each row and each column. Equivalently, a plane partition can be represented as a stack of unit cubes, with each entry specifying the height of the corresponding column of cubes, as illustrated in \cref{fig:example-TSSCPP}. A plane partition is \emph{contained in an $a\times b\times c$ box} if it has at most $a$ non-zero rows, at most $b$ non-zero columns, and all entries are at most $c$. A plane partition contained in a prescribed finite box is called a
\emph{boxed plane partition}. In the unit-cube representation, we indicate the ambient box by drawing
its walls, as shown in \cref{fig:cubes}.

A boxed plane partition is \emph{totally symmetric} if its associated configuration of unit cubes is invariant under every permutation of the three coordinate axes. The \emph{complement} of a boxed plane partition, with respect to its ambient box, is obtained by interchanging occupied and unoccupied cubes and reflecting the resulting configuration through the centre of the box. A boxed plane partition is \emph{self-complementary} if it coincides with its complement. A \emph{totally symmetric self-complementary plane partition} (TSSCPP) of order $N$ is a plane partition contained in a $2N \times 2N \times 2N$ box that is both totally symmetric and self-complementary.

As \cref{fig:cubes} suggests, boxed plane partitions can equivalently be viewed as \emph{lozenge tilings} (or \emph{rhombus tilings}) of hexagonal regions on the triangular lattice. Since such tilings are in bijection with perfect matchings of the dual honeycomb graph, they also admit a dimer-model representation, which we describe below.

\subsubsection{TSSCPPs as dimer configurations and Kasteleyn theory}
We recall the standard correspondence between lozenge tilings of a finite region $\mathcal{R}$ of the triangular lattice and perfect matchings of its dual graph $\mathcal{G}$. A lozenge tiling of $\mathcal{R}$ is a covering of $\mathcal{R}$ by $60^{\circ}$-rhombi, called \textit{lozenges}, with no gaps or overlaps. The dual graph $\mathcal{G}$ has one vertex for each unit triangle of $\mathcal{R}$, and two vertices are joined by an edge whenever the corresponding unit triangles share a side. Thus, $\mathcal{G}$ is a finite subgraph of the honeycomb lattice.

Each lozenge consists of two adjacent unit triangles and therefore determines a unique edge of $\mathcal{G}$, namely the edge joining the vertices corresponding to its two constituent triangles. Selecting this edge for every lozenge produces a perfect matching of $\mathcal{G}$, since every vertex of $\mathcal{G}$ is incident to exactly one selected edge. Conversely, every perfect matching of $\mathcal{G}$ determines a unique lozenge tiling of $\mathcal{R}$. The selected edges are called \textit{dimers}, and a perfect matching is therefore also referred to as a \textit{dimer configuration}. We denote the set of dimer configurations by $\mathcal{M}(\mathcal{G})$.

\begin{figure}[ht]
    \centering
    \begin{subfigure}{0.31\textwidth}
        \centering
        \includegraphics[width=0.8\linewidth]{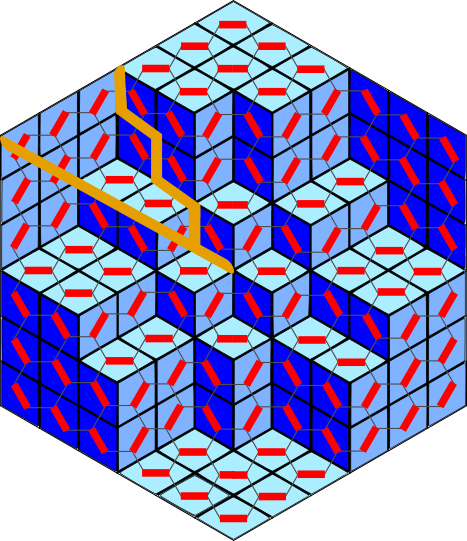}
        \caption{TSSCPP of order $3$ with highlighted fundamental domain in yellow.}
        \label{fig:TSSCPP-lozenges}
    \end{subfigure}
    \hfill
    \begin{subfigure}{0.31\textwidth}
        \centering
        \raisebox{0.8cm}{%
        \includegraphics[width=0.5\linewidth]{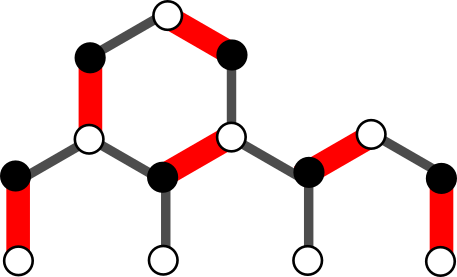}}
        \caption{Maximum matching on the fundamental domain.}
        \label{fig:free-boundary-dimers}
    \end{subfigure}
    \hfill
    \begin{subfigure}{0.31\textwidth}
        \centering
        \raisebox{0.8cm}{%
        \includegraphics[width=0.7\linewidth]{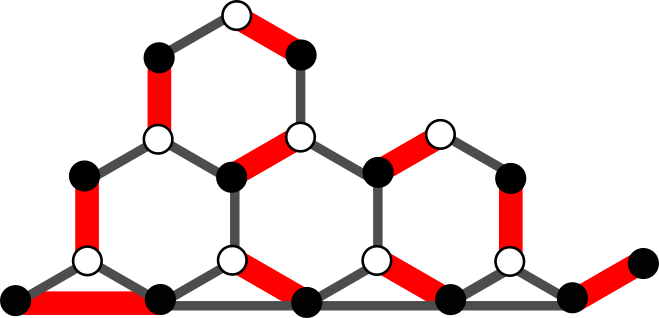}}
        \caption{Corresponding perfect matching on the augmented graph.}
        \label{fig:perfect-matching-dimers}
    \end{subfigure}
    \caption{Extracting the relevant dimer model from TSSCPPs. Dimers are shown in red.}
    \label{fig:TSSCPP-dimers}
\end{figure}

Applying the preceding correspondence to TSSCPPs, we may represent them by matching configurations on an appropriate fundamental domain. A TSSCPP of order $N$ corresponds to a lozenge tiling of a regular hexagon of side length $2N$. By total symmetry and self-complementarity, the full tiling is uniquely determined by its restriction to a one-twelfth fundamental domain of the hexagon. Moreover, the outermost lozenges and those along the three symmetry axes are forced. \cref{fig:TSSCPP-lozenges} shows an example, with the corresponding dimers superimposed on the lozenge tiling. 

Along the straight radial side of the one-twelfth fundamental domain, lozenges may cross the boundary, giving rise to a free boundary condition. In the restricted dual graph, these lozenges correspond to edges incident to degree-one vertices along this boundary. Although every vertex is matched in the dual graph of the full hexagon, these boundary vertices may remain unmatched after restricting to the fundamental domain. The admissible configurations on the fundamental domain thus correspond to \emph{maximum matchings} of the restricted dual graph, that is, matchings containing the largest possible number of edges. The maximum matching corresponding to the TSSCPP shown in \cref{fig:TSSCPP-lozenges} is illustrated in \cref{fig:free-boundary-dimers}.

Following \cite{ayyerCorrelations2021}, we augment the restricted dual graph by adding auxiliary vertices and edges along its free boundary. To follow the indexing convention of Ayyer and Chhita, we set $n:=N-1$ for a TSSCPP of order $N$ and denote the augmented graph by $\mathcal{G}_n$. Its vertex and edge sets are specified below. We add $2n+2$ auxiliary vertices when $n$ is even and $2n+1$ when $n$ is odd, connecting them in a triangular arrangement as illustrated in \cref{fig:perfect-matching-dimers}. These added triangles make the graph $\mathcal{G}_n$ non-bipartite. Their presence, which is related to the reflection symmetry imposed on TSSCPPs, arises for reasons similar to those in the model of free boundary dimers analysed in \cite{BLQ}. 
When $n$ is even (equivalently, when $N$ is odd), the rightmost added vertex has degree one. Consequently, its incident edge belongs to every perfect matching, as illustrated in \cref{fig:perfect-matching-dimers}. As proved in \cite[Proposition 2.2]{ayyerCorrelations2021}, the perfect matchings of $\mathcal{G}_n$ are in bijection with the maximum matchings of the restricted dual graph, and hence with TSSCPPs of order $N$.

\begin{figure}[ht]
     \centering
     \begin{subfigure}{0.2\textwidth}
         \centering
         \includegraphics[scale=0.4]{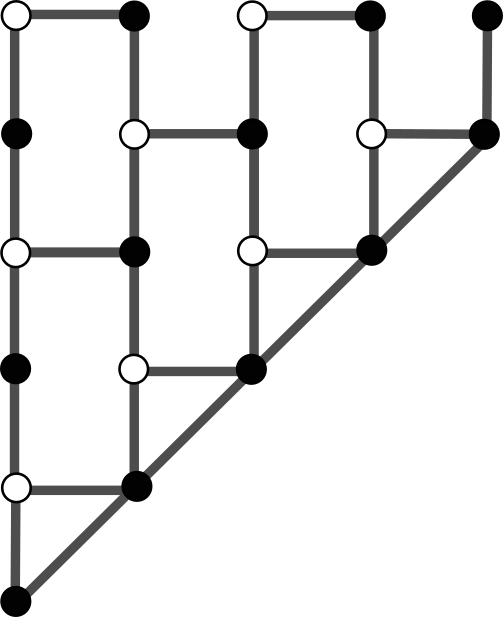}
         \caption{Graph $\mathcal{G}_2$.}
         \label{fig:G2}
     \end{subfigure}
     \hfill
     \begin{subfigure}[b]{0.24\textwidth}
         \centering
         \includegraphics[scale=0.4]{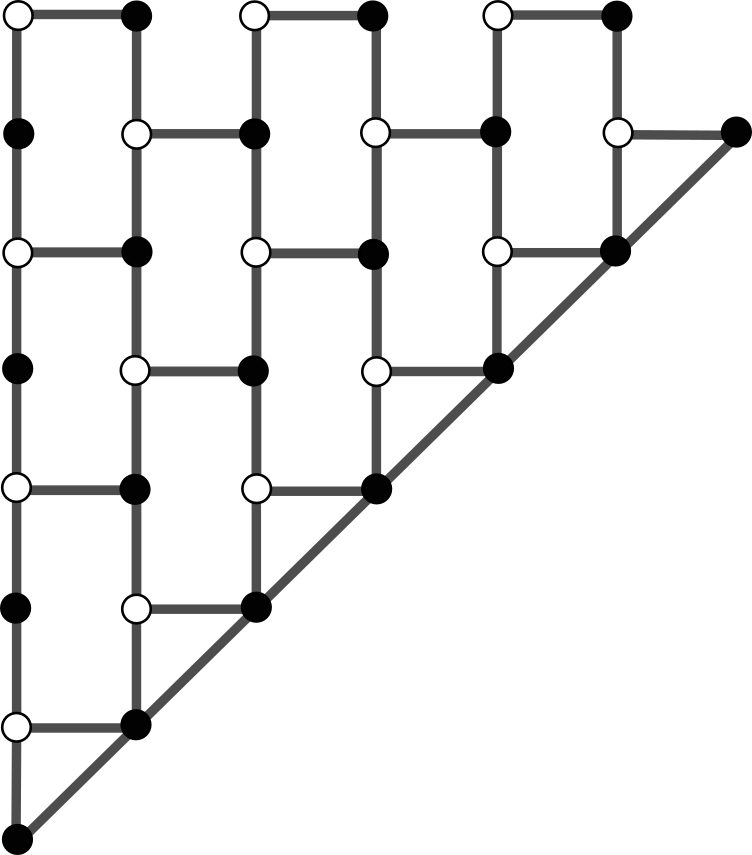}
         \caption{Graph $\mathcal{G}_3$.}
         \label{fig:G3}
     \end{subfigure}
          \hfill
     \begin{subfigure}[b]{0.24\textwidth}
         \centering
         \includegraphics[scale=0.4]{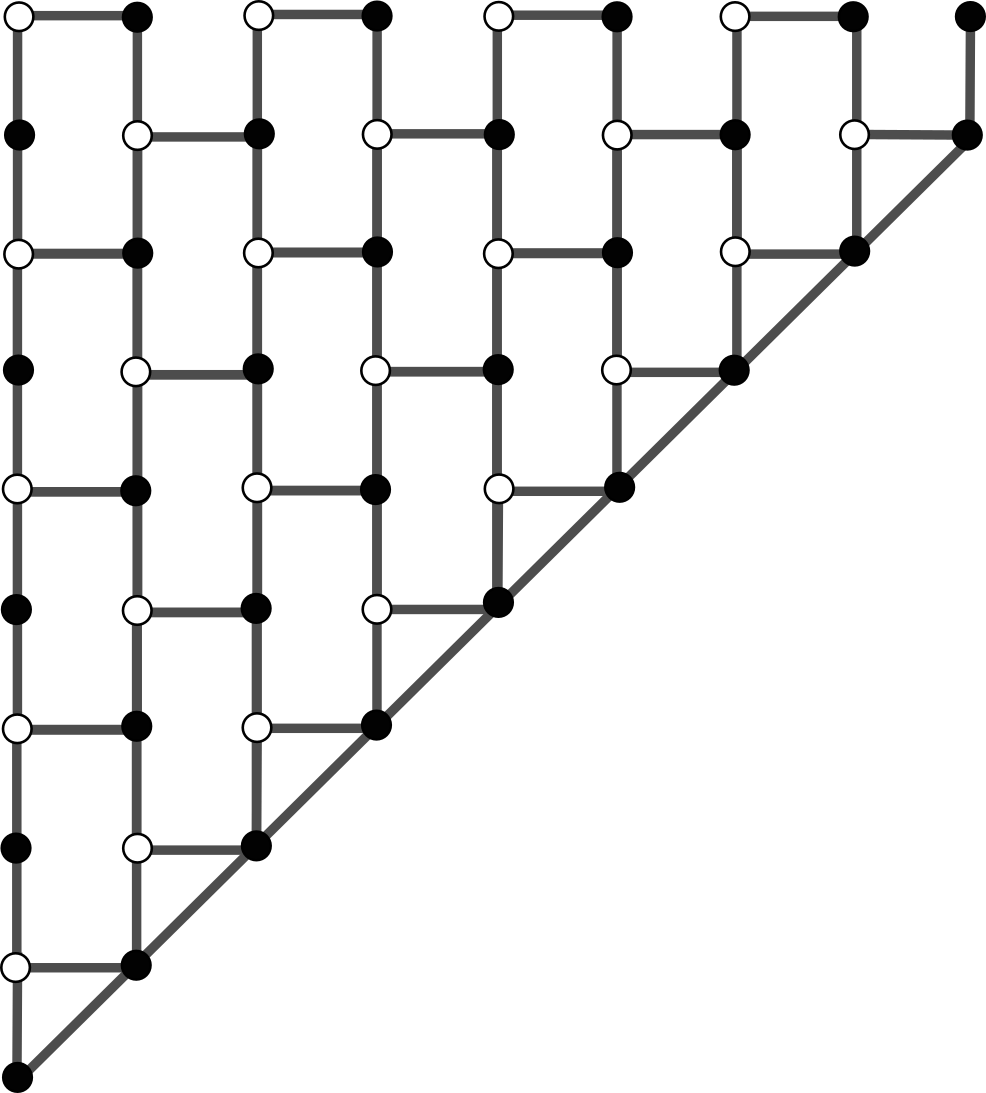}
         \caption{Graph $\mathcal{G}_4$.}
         \label{fig:G4}
     \end{subfigure}
          \hfill
     \begin{subfigure}[b]{0.24\textwidth}
         \centering
         \includegraphics[scale=0.4]{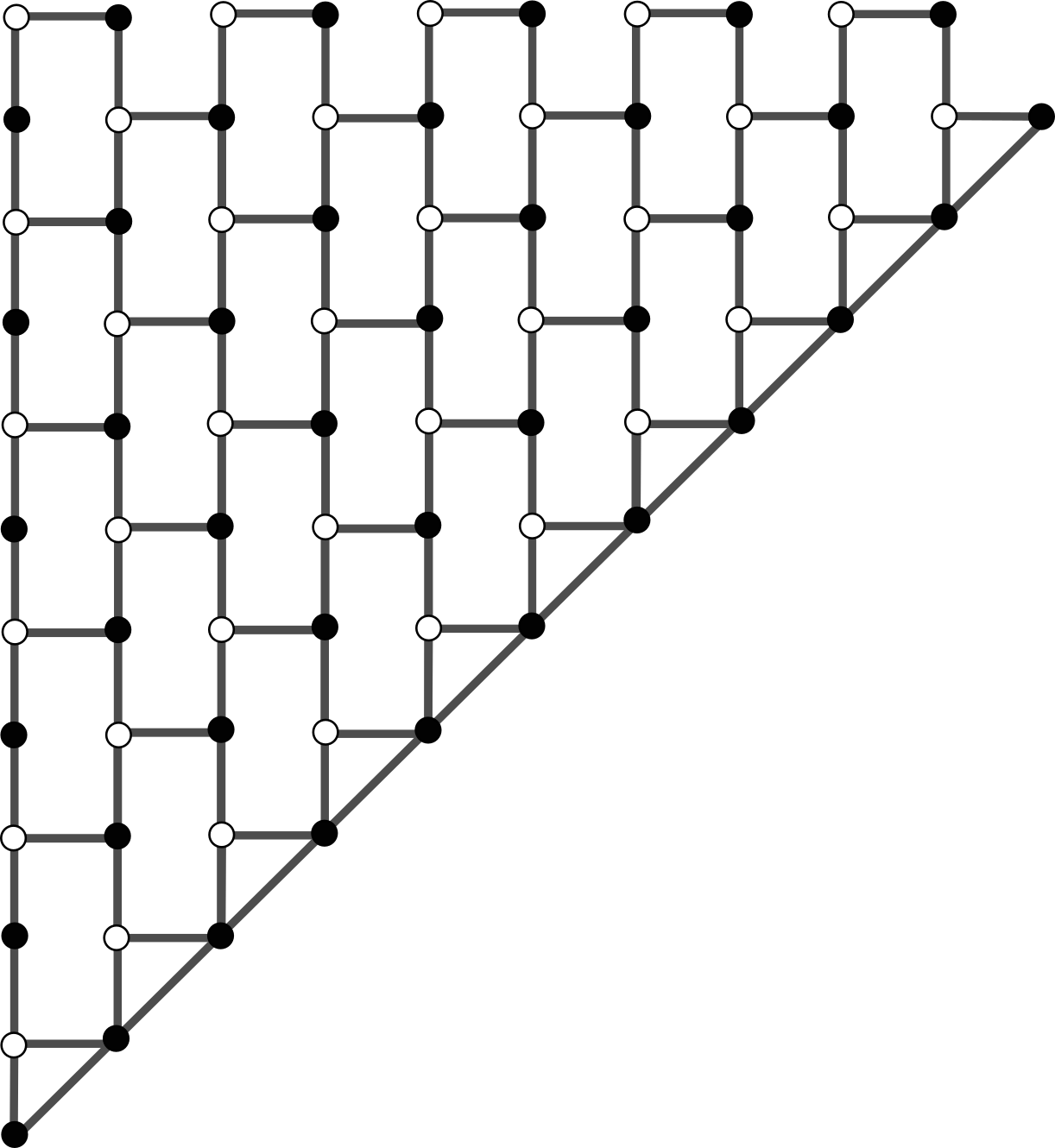}
         \caption{Graph $\mathcal{G}_5$.}
         \label{fig:G5}
     \end{subfigure}
        \caption{The graph $\mathcal{G}_n$ for $n=2,3,4,5$.}
        \label{fig:G}
\end{figure}

After rectifying the honeycomb lattice as illustrated in \cref{fig:G}, we represent the resulting graph by $\mathcal{G}_n=(V_n,E_n)$, with vertex set
\begin{align}
    V_n = \{(x_1,x_2)\mid 0\le x_1\le 2n,\ x_1\le x_2\le 2n+1\} \setminus \left\{ (2n,2n+1)\;\mid\; n\text{ is odd} \right\}.
\end{align}
The edge set is
\begin{align}
    E_n =   & \{((x_1,x_2),(x_1+1,x_2))\mid 0\le x_1\le 2n-1,\ x_1\le x_2\le 2n+1,\ [x_1+x_2]_2=1\} \notag \\
            &\cup \{((x_1,x_2),(x_1,x_2+1))\mid 0\le x_1\le 2n-\mathbbm{1}_{n\in 2\mathbb{Z}+1},\ x_1\le x_2\le 2n\} \notag \\
            &\cup \{((x_1,x_1),(x_1+1,x_1+1))\mid 0\le x_1\le 2n-1\}.
\end{align}
We use the standard chessboard colouring, taking the origin $(0,0)$ to be black. Thus, a vertex $(x_1,x_2)$ is black if $[x_1+x_2]_2=0$ and white otherwise. In particular, all diagonal vertices are black. The diagonal edges create odd cycles (triangles); in particular, the graph $\mathcal{G}_n$ is not bipartite, as already mentioned.

The graphs $\mathcal{G}_n$ for $n=2,3,4,5$ are illustrated in \cref{fig:G}. \cref{fig:G2} shows the rectification of the underlying graph in \cref{fig:perfect-matching-dimers}.

We now equip $\mathcal{G}_n$ with the Kasteleyn orientation used in \cite{ayyerCorrelations2021}. The associated Kasteleyn matrix $K_n = \left(  K_n(x,y) \right)_{x,y\in V_n}$ is the skew-symmetric signed adjacency matrix defined by $K_n(x,y)=1$ if the edge between $x$ and $y$ is oriented from $x$ to $y$, $K_n(x,y)=-1$ if it is oriented from $y$ to $x$, and $K_n(x,y)=0$ if $x$ and $y$ are not adjacent. In the convention adopted here, an orientation is Kasteleyn if every bounded face has an odd number of boundary edges oriented counter-clockwise. In the unweighted setting considered here, this condition ensures that the number of dimer configurations on $\mathcal{G}_n$ is given by $|\operatorname{Pf}(K_n)|$.

\begin{figure}[ht]
	\centering
	\includegraphics[width=0.3\textwidth]{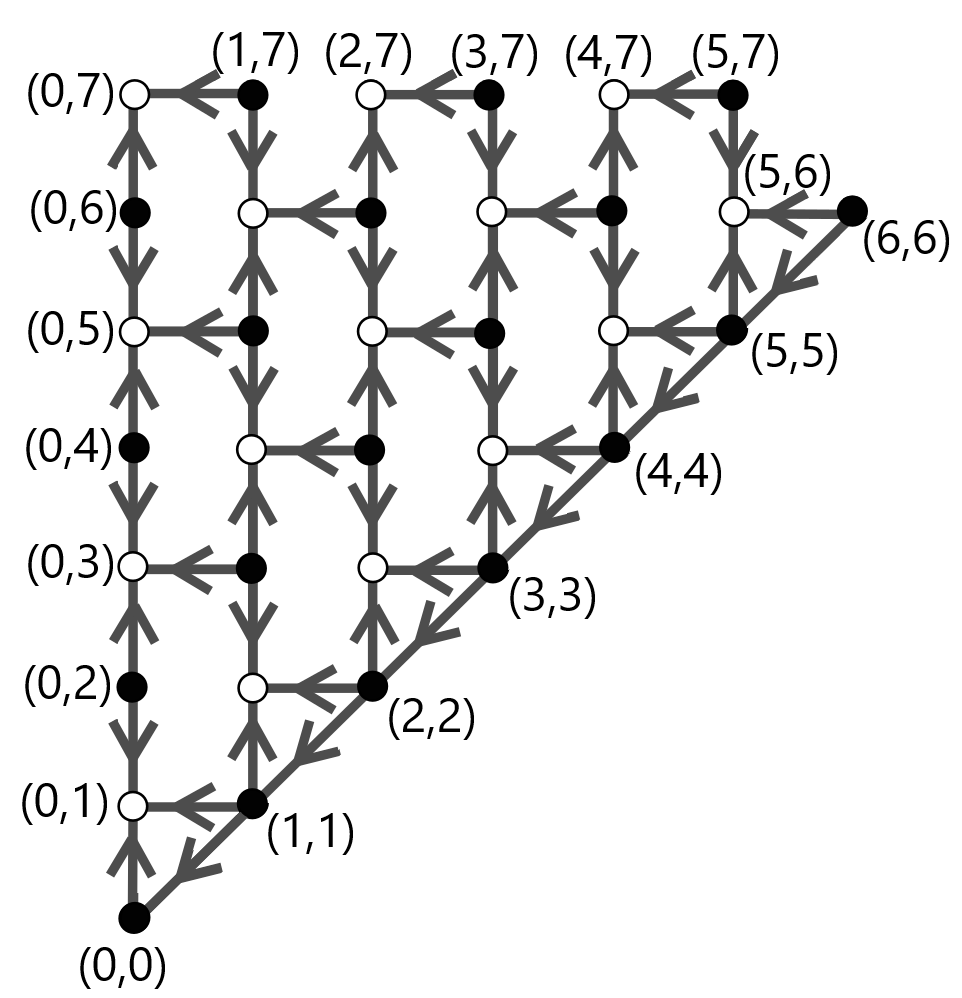}
	\caption{Kasteleyn orientation of $\mathcal{G}_3$, with selected boundary vertex coordinates indicated.}
	\label{fig:Oriented-TSSCPP}
\end{figure}

To specify the matrix explicitly, we write
\begin{align*}
    K_n(x,y) = k_n(x,y) - k_n(y,x)
\end{align*}
where, for $x=(x_1,x_2)$ and $y=(y_1,y_2)$,
\begin{align}
    k_n(x,y) = 
    \begin{cases} 
        1 & \text{if } [x_1+x_2]_2=0,\ x_2=y_2,\ x_1-y_1=1,\\
        1 & \text{if } [x_1+x_2]_2=0,\ |y_2-x_2|=1,\ x_1=y_1,\\
        1 & \text{if } x_1=x_2,\ y_1=y_2=x_1-1,\\
        0 & \text{otherwise.}
        \label{eq:Kasteleyn-orientation}
\end{cases}
\end{align}
By construction, $K_n$ is skew-symmetric and has entries in $\{-1,0,1\}$. The first case in \eqref{eq:Kasteleyn-orientation} orients every horizontal edge from a black vertex to a white vertex. The second case orients every vertical edge from black to white. Finally, the third case orients each diagonal edge from $(k,k)$ to $(k-1,k-1)$, for $k \in \{1,\dots,2n\}$. The resulting orientation of $\mathcal{G}_3$ is illustrated in \cref{fig:Oriented-TSSCPP}. 

Ayyer and Chhita \cite{ayyerCorrelations2021} derived explicit formulae for the inverse Kasteleyn matrix. Together with the standard Pfaffian correlation formula \cite[Theorem 3.3]{ayyerCorrelations2021}, these expressions allow us to compute probabilities of prescribed dimer events. Let $D$ be a uniformly random dimer configuration on $\mathcal{G}_n$, and let $U=\{e_1,\ldots,e_m\}$ be a set of edges of $\mathcal{G}_n$, where $e_r=\{v_{2r-1},v_{2r}\}$ for $r=1,\ldots,m$, then
\begin{equation}
    \mathbb{P}(U\subseteq D) = \prod_{r=1}^{m}K_n(v_{2r-1},v_{2r})\operatorname{Pf}\left( \left(K_n^{-1}(v_i,v_j)\right)_{1\leq i,j\leq 2m}^{\mathsf{T}} \right).
    \label{eq:dimer-correlations}
\end{equation}

We next recall the Fonseca--Zinn-Justin correspondence relating nested-arch probabilities in the $O(1)$ loop model to constrained TSSCPP enumerations.

\subsection{The Fonseca--Zinn-Justin correspondence for nested submatching probabilities}
\label{sec:fonseca-zinn-justin-correspondence}
\begin{figure}[ht]
     \centering
     \begin{subfigure}{0.46\textwidth}
        \centering
        \raisebox{1cm}{%
        \includegraphics[width=0.19\linewidth]{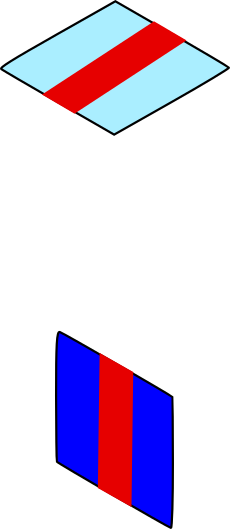}}
        \caption{The two local path segments.}
        \label{fig:lozenges-paths}
     \end{subfigure}
     \hfill
     \begin{subfigure}[b]{0.46\textwidth}
         \centering
         \includegraphics[width=0.64\linewidth]{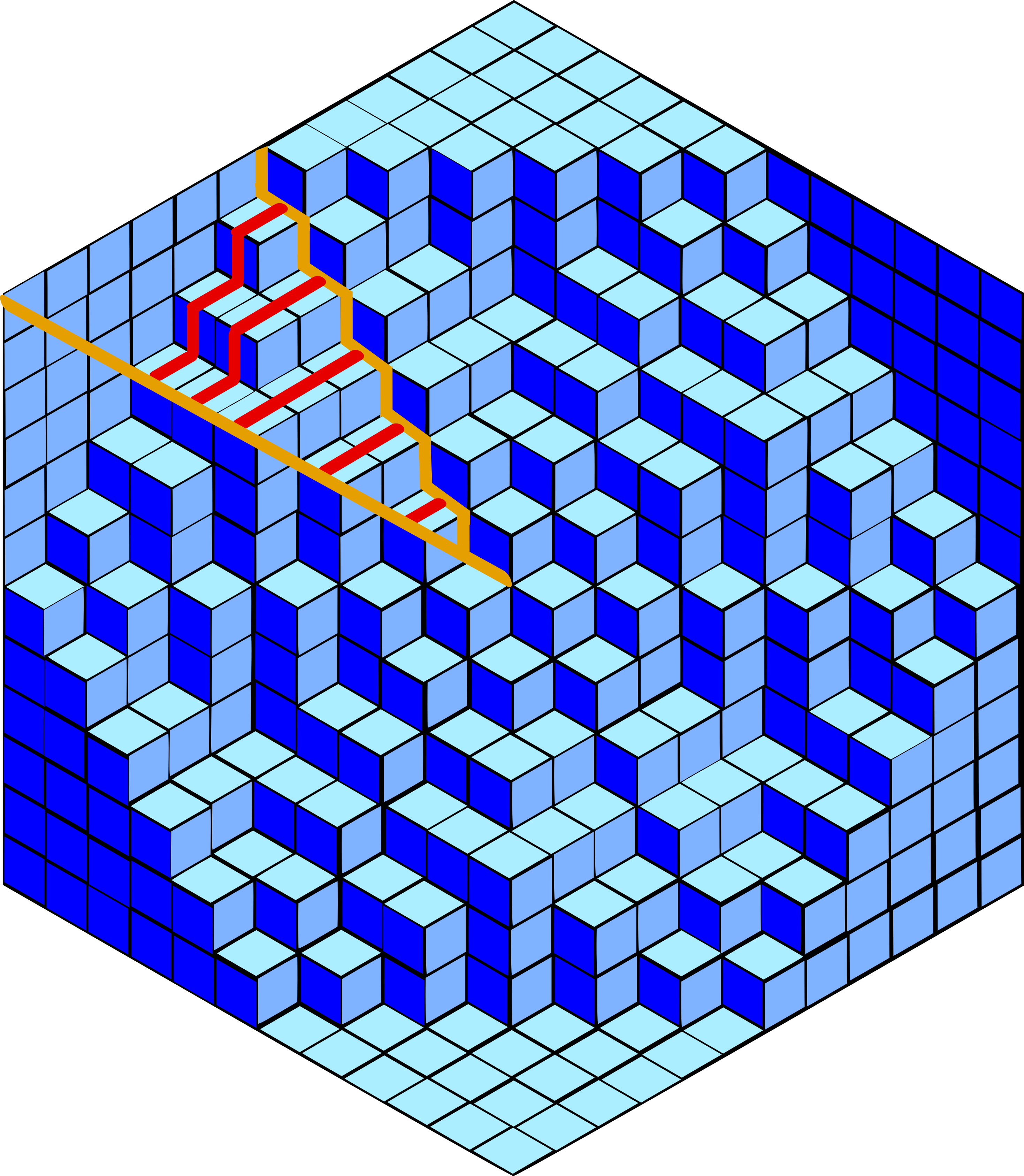}
         \caption{NILPs associated with a TSSCPP.}
        \label{fig:example-paths-TSSCPP}
     \end{subfigure}
        \caption{Construction of NILPs from a TSSCPP.}
    \label{fig:NILP}
\end{figure}

We adopt the non-intersecting lattice-path representation of TSSCPPs used by Fonseca and Zinn-Justin \cite{fonsecaGround2009}. Within the fundamental domain introduced above, we draw path segments through lozenges of two of the three possible orientations, as shown in \cref{fig:lozenges-paths}. Lozenges of the third orientation contribute no path segments. Joining these segments produces a family of $N-1$ non-intersecting lattice paths (NILPs), which we number from the centre outwards; see \cref{fig:example-paths-TSSCPP}.

Of particular interest are configurations whose first $k$ paths are horizontal, that is, contain no downward steps. Fonseca and Zinn-Justin include a formal null path of length zero, indexed by $0$, in their enumeration \cite[Section 2.3]{fonsecaGround2009}. We omit this degenerate path. Since it imposes no constraint, requiring their first $k+1$ paths to be horizontal is equivalent to requiring our first $k$ paths to be horizontal, for $1\leq k<N$. In the corresponding TSSCPPs, this constraint fixes the tiling of a central hexagonal region. An example with the first three paths horizontal is shown in \cref{fig:example-paths-TSSCPP}.

Let $\mathcal{T}_N$ denote the set of TSSCPPs of order $N$. For every $N>k\geq 1$, let $\mathcal{T}_{N,k}\subseteq\mathcal{T}_N$ consist of those whose first $k$ paths are horizontal. By the enumeration theorems of Andrews \cite{andrews1994} and Zeilberger \cite{zeilbergerProofAlternatingSign1996}, the number of TSSCPPs of order $N$ equals the number $A_N$ of alternating-sign matrices of order $N$: \begin{equation*}     
    |\mathcal{T}_N|= A_N =\prod_{j=0}^{N-1}\frac{(3j+1)!}{(N+j)!}
\end{equation*}

Fonseca and Zinn-Justin \cite{fonsecaGround2009} work with a ground-state vector of the $O(1)$ loop model, whose components are indexed by link patterns and, when normalised to sum to one, give the stationary distribution described before. They choose an overall normalisation in which certain sums of these components count constrained NILP configurations. Their proof compares contour-integral expressions for the loop-model sums, derived from a polynomial solution of the qKZ equation, with those for the NILP generating functions, obtained using the Lindström--Gessel--Viennot lemma. At $\tau=1$, every NILP configuration has weight one, so these generating functions give unweighted counts. The resulting correspondence is established without constructing a direct bijection between loop configurations and TSSCPPs.

The following formulation of the Fonseca--Zinn-Justin correspondence \cite{fonsecaGround2009} is adapted to our conventions.
\begin{theorem}[Fonseca--Zinn-Justin {\cite{fonsecaGround2009}}]
\label{thm:fonseca-zinn-justin}
For every $N>k\geq 1$,
    \begin{equation}
        \mathbb{P}_{\mathcal{C}_N}(F_k) = \frac{|\mathcal{T}_{N,k}|}{|\mathcal{T}_N|}.
    \end{equation}
\end{theorem}

Normalising the ground-state components by their total weight identifies the proportion of TSSCPPs satisfying this constraint with the probability of $k$ nested arches enclosing an arbitrary link pattern in the cylindrical $O(1)$ loop model. A cyclic rotation of the boundary labels maps these enclosing arches to the nested pairs defining $F_k$. By rotational invariance of the cylindrical model, the corresponding probabilities agree.

The endpoint $k=N$ can be treated separately. For $N\geq 2$, the event $F_{N-1}$ fixes $N-1$ pairs, forcing the two remaining boundary points, $2N-1$ and $2N$, to be paired. Thus, $F_{N-1}$ specifies a unique link pattern, whose cyclic rotation by one site gives the link pattern specified by $F_N$. By rotational invariance and \cref{thm:fonseca-zinn-justin},
\begin{equation}
    \mathbb{P}_{\mathcal{C}_N}(F_N) = \mathbb{P}_{\mathcal{C}_N}(F_{N-1}) = \frac{|\mathcal{T}_{N,N-1}|}{|\mathcal{T}_N|} = \frac{1}{|\mathcal{T}_N|}= \frac{1}{A_N},
    \label{eq:rotational-inv}
\end{equation}
since requiring all $N-1$ nontrivial paths to be horizontal determines a unique TSSCPP. For $N=1$, we have $\mathbb{P}_{\mathcal{C}_1}(F_1)=1$.

In the subsequent analysis, we identify the constraint defining $\mathcal{T}_{N,k}$ with the presence of prescribed dimers, allowing the probability in \cref{thm:fonseca-zinn-justin} to be evaluated using the Pfaffian correlation formula.

\section{Simplification of the inverse Kasteleyn matrix}
\label{sec:simplified-K-inv}
We simplify the inverse Kasteleyn matrix formulae of Ayyer and Chhita \cite{ayyerCorrelations2021}, in preparation for their application to the vertices associated with our prescribed-dimer events. 

We adopt their coordinate convention, writing vertices as
\begin{align}     
    & (x_1,x_2) = (i_1,i_1+2i_2+\epsilon_i), \nonumber \\
    & \text{ where $0\leq i_1\leq 2n-1$, $\epsilon_i\in\{0,1\}$, and $0\leq i_2\leq n-\lfloor(i_1+\epsilon_i)/2\rfloor$.} \label{eq:rangei_1-i_2}
\end{align} 
 Here, $\epsilon_i$ records the parity of the vertex under the chessboard colouring introduced above. Analogously, we write a second vertex as $(y_1,y_2)=(j_1,j_1+2j_2+\epsilon_j)$, with the corresponding restrictions on $j_1,j_2,\epsilon_j$.

The formulae of Ayyer and Chhita involve the auxiliary function 
\begin{align}     
    p(n,k,\ell) = \frac{(n+k-2\ell+1)!(2n-k-\ell+1)!}{(k-\ell)!(3n-k+2-2\ell)!}(-1)^k(3n-3k+2)C_{n-k}.
    \label{eq:p-def} 
\end{align} 
Recall that $C_n=\frac{1}{n+1}\binom{2n}{n}$ denotes the $n$-th Catalan number. We consider $p(n,k,\ell)$ for integers $0\leq \ell \leq k \leq n$ and extend it by zero outside this range. 

The identity
\begin{align}     
    p(n,k,\ell)=(-1)^\ell p(n-\ell,k-\ell,0),\qquad 0\leq \ell\leq k\leq n    
    \label{eq:p-ell} 
\end{align} 
follows directly from \eqref{eq:p-def} and shows that $p$ is determined by its values at $\ell=0$.

The following quantity will appear repeatedly in our formulae for the entries of the inverse Kasteleyn matrix: 
\begin{equation}\label{eq:S_n}
S_n(k_1, k_2) = \sum_{ \ell_1 =0}^{k_1}\sum_{  \ell_2 =0}^{ k_2} p(n,k_1,\ell_1) \, p(n, k_2, \ell_2) \left( \mathbbm{1}_{\ell_2 = \ell_1+1}\! -  \!\mathbbm{1}_{\ell_2 = \ell_1-1}  \right).
\end{equation}
In particular, this quantity will give rise to a hypergeometric sum that we evaluate using Zeilberger’s and Gosper’s algorithms.

The next proposition gives the entries in the three independent parity cases; the fourth follows from skew-symmetry. 

\begin{proposition}
\label{prop:simple-inverse-Kasteleyn}
Let $x=(x_1,x_2)$ and $y=(y_1,y_2)$ be vertices represented in the coordinates above. 
 \begin{itemize}[leftmargin=*]
        \item For $\epsilon_i=1$ and $\epsilon_j = 1$,                \begin{align}
                 K_n^{-1} (x,y)\label{eq:K11Lemma}
                &=  (-1)^{i_2 + j_2} \sum_{k_1= 0}^n \sum_{k_2= 0}^n   S_n(k_1, k_2) \binom{n-k_1-i_2}{i_1-2k_1} \binom{n-k_2-j_2}{j_1-2k_2}.
            \end{align}            
        \item For $\epsilon_i=0$ and $\epsilon_j = 0$, $K_n^{-1} (x,y) $ is given by 
            \begin{align}
            \label{eq:K00Lemma}
                 & (-1)^{i_2 + j_2} \Biggl( \sum_{k_1 = 0}^n \sum_{k_2 = 0}^n S_n(k_1,k_2) \sum_{m_1=0}^{\infty} \binom{n-k_1+i_2}{i_1+ i_2 - 2k_1 - m_1 - 1} \sum_{m_2=0}^{\infty} \binom{n-k_2+j_2}{j_1+ j_2 - 2k_2 - m_2 - 1} \notag \\ 
                  & \quad +  ([j_1 + 1]_2 \, \mathbbm{1}_{j_2 = 0} + 2^{j_2 - 1} \, \mathbbm{1}_{j_2 \ge 1}) \sum_{k=0}^n p(n,k,0) \sum_{m = 0}^{\infty} \binom{n - k + i_2}{i_1 + i_2 - 2k - m - 1} \notag \\
                  &\quad -  ([i_1 + 1]_2 \, \mathbbm{1}_{i_2 = 0} + 2^{i_2 - 1} \, \mathbbm{1}_{i_2 \ge 1}) \sum_{k=0}^n p(n,k,0) \sum_{m = 0}^{\infty} \binom{n - k + j_2}{j_1 + j_2 - 2k - m - 1}  \notag \\
                  & \quad - \sum_{\ell_1'=0}^{i_2-1}\sum_{\ell_2'=0}^{j_2-1} \binom{i_2-1}{\ell_1'}\binom{j_2-1}{\ell_2'} \operatorname{sgn}(i_1\!+\!\ell_1'\!-\!j_1\!-\!\ell_2') \!+\!\mathbbm{1}_{j_2=0} \Xi^{j}(x,y) \!+\!\mathbbm{1}_{i_2=0} \Xi^{i}(x,y) \!+\!\mathbbm{1}_{i_2=0}\mathbbm{1}_{j_2=0} \Xi^{i,j}(x,y) \Biggr),
            \end{align}
            where the ``error terms'' $\Xi^i, \Xi^j$ and $\Xi^{i,j}$ are explicitly given by 
            \begin{align}
                \Xi^{i}(x,y)
                    &:=\sum_{\ell_2'=0}^{j_2-1} \binom{j_2-1}{\ell_2'} \Bigl( \mathbbm{1}_{i_1\le j_1+\ell_2'}[i_1+1]_2 - \mathbbm{1}_{i_1>j_1+\ell_2'}[i_1]_2 \Bigr), \\
                \Xi^{j}(x,y)
                    &:= \sum_{\ell_1'=0}^{i_2-1} \binom{i_2-1}{\ell_1'} \Bigl( \mathbbm{1}_{i_1+\ell_1'<j_1}[j_1]_2 - \mathbbm{1}_{i_1+\ell_1'\ge j_1}[j_1+1]_2 \Bigr),\\
                \Xi^{i,j}(x,y)
                    &:= \mathbbm{1}_{i_1<j_1}[i_1+1]_2[j_1]_2 - \mathbbm{1}_{i_1>j_1}[i_1]_2[j_1+1]_2.
            \end{align}
     \item For $\epsilon_i=1$ and $\epsilon_j = 0$,
            \begin{align}
                \label{eq:K10Lemma}
                K_n^{-1} (x,y) 
                &= (-1)^{i_2 + j_2} \Biggl( \sum_{k_1 = 0}^n \sum_{k_2 = 0}^n S_n(k_1,k_2) \binom{n-k_1 - i_2}{i_1 - 2k_1} \sum_{m=0}^{\infty} \binom{n-k_2 + j_2}{j_1 + j_2 - 2k_2 - m -1} \notag \\
                &\quad + ([j_1 + 1]_2 \, \mathbbm{1}_{j_2 = 0} + 2^{j_2 - 1} \, \mathbbm{1}_{j_2 \ge 1}) \sum_{k=0}^n p(n,k,0) \binom{n-k-i_2}{i_1 - 2k} - \mathbbm{1}_{i_2<j_2}\binom{j_2 - i_2 - 1}{i_1 - j_1} \Biggr).
            \end{align}
        \item For $\epsilon_i=0$ and $\epsilon_j = 1$,
            \begin{align}
                K_n^{-1}(x,y) = - K_n^{-1}(y,x).
            \end{align}
    \end{itemize}
\end{proposition}

In the above statement, we view the terms $\Xi^i, \Xi^j$ and $\Xi^{i,j}$ as error terms because they occur only when evaluating the inverse Kasteleyn matrix at vertices belonging to the triangular faces adjoining the diagonal boundary of $\mathcal{G}_n$. These will not be needed in our analysis of nested-arch events. 

Note that the preceding proposition is independent of the convention adopted for $\binom{a}{0}$ when $a<0$, owing to the indicator $\mathbbm{1}_{i_2<j_2}$ in \eqref{eq:K10Lemma}. We verify this separately for each parity case in the proof.

\begin{proof}[Proof of \cref{prop:simple-inverse-Kasteleyn}]

We begin with an auxiliary lemma needed for the case $\epsilon_i,\epsilon_j=0$. The resulting simplification is particularly useful when the ranges of $i_2$ and $j_2$ are further restricted, as is the case for nested-arch events.

\begin{lemma} 
\label{lem:indicators-sgn}
Let $i_1,j_1$ be non-negative integers. The following identities hold. 

\begin{enumerate}
    \item \label{itm:boundary-case1} For every integer $j_2\ge 1$,
        \begin{align}
            &\sum_{\ell_2' = 0}^{j_2} \binom{j_2}{\ell_2'} \left( \mathbbm{1}_{i_1 < j_1 + \ell_2'} [i_1 + 1]_2 \, [j_1 + \ell_2']_2 - \mathbbm{1}_{i_1 > j_1 + \ell_2'} [i_1]_2 \, [j_1 + \ell_2' + 1]_2 \right)\notag \\
            &= \sum_{\ell_2'=0}^{j_2 - 1} \binom{j_2 - 1}{\ell_2'} \left( \mathbbm{1}_{i_1 \le j_1 + \ell_2'} [i_1 + 1]_2 - \mathbbm{1}_{i_1 > j_1 + \ell_2'} [i_1]_2 \right).
        \end{align}
    \item \label{itm:boundary-case2} For every integer $i_2 \ge 1$,
        \begin{align}
            &\sum_{\ell_1' = 0}^{i_2} \binom{i_2}{\ell_1'} \left( \mathbbm{1}_{i_1 + \ell_1' < j_1} [i_1 + \ell_1' + 1]_2 \, [j_1]_2 - \mathbbm{1}_{i_1 + \ell_1' > j_1 } [i_1 + \ell_1']_2 \, [j_1 + 1]_2 \right)\notag \\
            &= \sum_{\ell_1'=0}^{i_2 - 1} \binom{i_2 - 1}{\ell_1'} \left( \mathbbm{1}_{i_1 + \ell_1' < j_1} [j_1]_2 - \mathbbm{1}_{i_1 + \ell_1' \ge j_1 } [j_1 + 1]_2 \right).
        \end{align}
    \item \label{itm:lemma-indicators-sgn-1} For every integer $j_2 \ge 1$,
        \begin{align}
            \sum_{\ell_2' = 0}^{j_2} \binom{j_2}{\ell_2'} \left( \mathbbm{1}_{i_1 < j_1 + \ell_2'} [j_1 + \ell_2']_2 - \mathbbm{1}_{i_1 \geq j_1 + \ell_2'} [j_1 + \ell_2' + 1]_2 \right) = - \sum_{\ell_2' = 0}^{j_2 - 1} \binom{j_2 - 1}{\ell_2'} \, \operatorname{sgn}(i_1 - j_1 - \ell_2'). 
            \label{eq: indicators-sign-identity with one variable}
            \end{align}        
    \item \label{itm:lemma-indicators-sgn-2}     For all integers $i_2,j_2\geq1$,
         \begin{align}
            &\sum_{\ell'_1=0}^{i_2} \sum_{\ell'_2=0}^{j_2}  \binom{i_2}{\ell'_1} \binom{j_2}{\ell'_2} \biggl( \mathbbm{1}_{i_1 + \ell_1'<j_1 + \ell_2'} [i_1 + \ell'_1 + 1]_2 \, [j_1 + \ell'_2]_2 -  \mathbbm{1}_{i_1 + \ell'_1> j_1 + \ell_2'} [i_1+\ell'_1]_2 \, [j_1 + \ell_2' + 1]_2  \biggr) \notag \\
            &= - \sum_{\ell'_1 = 0}^{i_2 - 1}  \sum_{\ell'_2 = 0}^{j_2 - 1}  \binom{i_2 - 1}{\ell'_1} \binom{j_2 - 1}{\ell'_2} \, \operatorname{sgn} (i_1 + \ell_1' - j_1 - \ell_2'). \label{eq:indicators-sign-identity}
        \end{align}
     \end{enumerate}
\end{lemma}

Note that in this lemma, the convention about the value of the binomial coefficient $\binom{a}{0}$ for $a<0$ is irrelevant since $i_2, j_2 \ge 1$. 

\begin{proof}
\begin{itemize}[leftmargin=*]
\item We first prove \cref{itm:boundary-case1}. Using Pascal's rule and then reindexing the second
resulting sum, we obtain
\begin{align}
    &\sum_{\ell_2' = 0}^{j_2} \binom{j_2}{\ell_2'} \left( \mathbbm{1}_{i_1 < j_1 + \ell_2'} [i_1 + 1]_2 \, [j_1 + \ell_2']_2 - \mathbbm{1}_{i_1 > j_1 + \ell_2'} [i_1]_2 \, [j_1 + \ell_2' + 1]_2 \right) \notag \\
    &= \sum_{\ell_2' = 0}^{j_2} \binom{j_2 - 1}{\ell_2'} \left( \mathbbm{1}_{i_1 < j_1 + \ell_2'} [i_1 + 1]_2 \, [j_1 + \ell_2']_2 - \mathbbm{1}_{i_1 > j_1 + \ell_2'} [i_1]_2 \, [j_1 + \ell_2' + 1]_2 \right) \notag \\
    &\quad + \sum_{\ell_2' = 0}^{j_2} \binom{j_2 - 1}{\ell_2' - 1} \left( \mathbbm{1}_{i_1 < j_1 + \ell_2'} [i_1 + 1]_2 \, [j_1 + \ell_2']_2 - \mathbbm{1}_{i_1 > j_1 + \ell_2'} [i_1]_2 \, [j_1 + \ell_2' + 1]_2 \right) \notag \\
    &= \sum_{\ell_2' = 0}^{j_2 - 1} \binom{j_2 - 1}{ \ell_2'} \Biggl\{ \left(  \mathbbm{1}_{i_1 < j_1 + \ell_2'} [i_1 + 1]_2 \, [j_1 + \ell_2']_2 - \mathbbm{1}_{i_1 > j_1 + \ell_2'} [i_1]_2 \, [j_1 + \ell_2' + 1]_2  \right) \notag \\
    &\quad + \left(  \mathbbm{1}_{i_1 < j_1 + \ell_2' + 1} [i_1 + 1]_2 \, [j_1 + \ell_2' + 1]_2 - \mathbbm{1}_{i_1 > j_1 + \ell_2' + 1} [i_1]_2 \, [j_1 + \ell_2' + 2]_2  \right) \Biggr\}. \label{line:parity-j1}
\end{align}
Here, $\binom{j_2-1}{j_2}=0$ since $j_2\ge 1$, while $\binom{j_2-1}{-1}=0$ by our convention.

To simplify this further, we use
\begin{align}
    \mathbbm{1}_{i_1 < j_1 + \ell_2'} [i_1 + 1]_2 \, [j_1 + \ell_2']_2 &= \mathbbm{1}_{i_1 \le j_1 + \ell_2'} [i_1 + 1]_2 \, [j_1 + \ell_2']_2,  \label{eq:314}\\
    \mathbbm{1}_{i_1 > j_1 + \ell_2' + 1} [i_1]_2 \, [j_1 + \ell_2' + 2]_2 &=  \mathbbm{1}_{i_1 > j_1 + \ell_2'} [i_1]_2 \, [j_1 + \ell_2']_2, \label{eq:315}
\end{align}
since the extra cases that were added on the right side do not contribute due to the vanishing product of parity factors.

We use \eqref{eq:314} to rewrite the first of the four terms in the curly brackets of \eqref{line:parity-j1}, and we use \eqref{eq:315} to rewrite the fourth term. 
After noting that $[j_1 + \ell'_2]_2 + [j_1 + \ell'_2 +1]_2 = 1$,
the four terms reduce to
\begin{align*}
        \mathbbm{1}_{i_1 \le j_1 + \ell_2'} [i_1 + 1]_2 - \mathbbm{1}_{i_1 > j_1 + \ell_2'} [i_1]_2.
\end{align*}
Substituting this expression into the preceding identity yields \cref{itm:boundary-case1}.

\item The proof of \cref{itm:boundary-case2} is analogous. Applying Pascal's rule and reindexing the second sum,
we use the identities
\begin{align}
    \mathbbm{1}_{i_1 + \ell_1' +1 < j_1 } [i_1 + \ell_1' + 2]_2 \, [j_1]_2 &= \mathbbm{1}_{i_1 + \ell_1' < j_1 } [i_1 + \ell_1' ]_2 \, [j_1]_2, \\
    \mathbbm{1}_{i_1 + \ell_1' + 1 > j_1 } [i_1 + \ell_1' + 1]_2 \, [j_1 + 1]_2 &=  \mathbbm{1}_{i_1 + \ell_1' \ge j_1 } [i_1 + \ell_1' + 1]_2 \, [j_1+1]_2.
\end{align}
to obtain
\begin{align}
     &\sum_{\ell_1' = 0}^{i_2} \binom{i_2}{\ell_1'} \left( \mathbbm{1}_{i_1 + \ell_1' < j_1} [i_1 + \ell_1' + 1]_2 \, [j_1]_2 - \mathbbm{1}_{i_1 + \ell_1' \ge j_1 } [i_1 + \ell_1']_2 \, [j_1 + 1]_2 \right)\notag \\
            &= \sum_{\ell_1'=0}^{i_2 - 1} \binom{i_2 - 1}{\ell_1'} \left( \mathbbm{1}_{i_1 + \ell_1' < j_1} [j_1]_2 - \mathbbm{1}_{i_1 + \ell_1' \ge j_1 } [j_1 + 1]_2 \right),
\end{align}
which proves \cref{itm:boundary-case2}.

\item We prove \cref{itm:lemma-indicators-sgn-1} by induction on $j_2$. For $j_2=1$, and for all $i_1, j_1 \ge 0$ we have
    \begin{align}
        &\sum_{\ell_2' = 0}^{1} \binom{1}{\ell_2'} \left( \mathbbm{1}_{i_1 < j_1 + \ell_2'} [j_1 + \ell_2']_2 - \mathbbm{1}_{i_1 \geq j_1 + \ell_2'} [j_1 + \ell_2' + 1]_2 \right) \notag \\
        &= \binom{1}{0} \left( \mathbbm{1}_{i_1 < j_1 } [j_1 ]_2 - \mathbbm{1}_{i_1 \geq j_1 } [j_1 + 1]_2 \right) + \binom{1}{1} \bigl( \mathbbm{1}_{i_1 < j_1 + 1 } [j_1 + 1]_2 - \mathbbm{1}_{i_1 \geq j_1 + 1} [j_1 + 2]_2 \bigr) \notag \\
        &= \mathbbm{1}_{i_1 < j_1 } [j_1 ]_2 - \mathbbm{1}_{i_1 > j_1} [j_1 + 1]_2 - \mathbbm{1}_{i_1 = j_1} [j_1 + 1]_2 + \mathbbm{1}_{i_1 < j_1} [j_1 + 1]_2 + \mathbbm{1}_{i_1 = j_1} [j_1 + 1]_2 - \mathbbm{1}_{i_1 > j_1} [j_1]_2 \notag \\
        &= \mathbbm{1}_{i_1 < j_1 } - \mathbbm{1}_{i_1 > j_1 }\notag  \\
        &= - \operatorname{sgn}(i_1 - j_1).\notag 
    \end{align}
For the induction step, assume that the identity holds for a given $j_2\geq 1$ and all $i_1,j_1\geq 0$. Set $\bar{j}=j_1+1$ and $\bar{\ell}=\ell_2'-1$. Pascal's rule gives
    \begin{align*}
        &\sum_{\ell_2' = 0}^{j_2 + 1} \binom{j_2 + 1}{\ell'_2}   \left( \mathbbm{1}_{i_1 < j_1 + \ell_2'} [j_1 + \ell_2']_2 - \mathbbm{1}_{i_1 \geq j_1 + \ell_2'} [j_1 + \ell_2' + 1]_2 \right) \notag \\
        &= \sum_{\ell'_2=0}^{j_2} \binom{j_2}{\ell'_2}  \left( \mathbbm{1}_{i_1 < j_1 + \ell_2'} [j_1 + \ell_2']_2 - \mathbbm{1}_{i_1 \geq j_1 + \ell_2'} [j_1 + \ell_2' + 1]_2 \right) \notag \\ 
        &\quad + \sum_{\ell'_2=0}^{j_2 + 1} \binom{j_2}{\ell'_2 - 1}  \left( \mathbbm{1}_{i_1 < j_1 + \ell_2'} [j_1 + \ell_2']_2 - \mathbbm{1}_{i_1 \geq j_1 + \ell_2'} [j_1 + \ell_2' + 1]_2 \right) \notag \\
        &= - \sum_{\ell_2' = 0}^{j_2} \binom{j_2 - 1}{\ell'_2} \, \operatorname{sgn}(i_1 - j_1 - \ell_2') + \sum_{\bar{\ell} = 0}^{j_2} \binom{j_2}{\bar{\ell}}  \left( \mathbbm{1}_{i_1 < \bar{j} + \bar{\ell}} [\bar{j} + \bar{\ell}]_2 - \mathbbm{1}_{i_1 \geq \bar{j} + \bar{\ell}} [\bar{j} + \bar{\ell} + 1]_2 \right) \notag \\
        &= - \sum_{\ell_2' = 0}^{j_2}\binom{j_2 - 1}{\ell'_2} \, \operatorname{sgn}(i_1 - j_1 - \ell_2') - \sum_{\bar{\ell} = 0}^{j_2 - 1}\binom{j_2 - 1}{\bar{\ell}} \, \operatorname{sgn}(i_1 - \bar{j} - \bar{\ell}) \notag \\
        &= - \sum_{\ell_2' = 0}^{j_2} \binom{j_2 - 1}{\ell'_2} \, \operatorname{sgn}(i_1 - j_1 - \ell_2') - \sum_{\ell'_2 = 0}^{j_2}\binom{j_2 - 1}{\ell'_2 - 1} \, \operatorname{sgn}(i_1 - j_1 - \ell'_2) \notag  \\
        &= - \sum_{\ell_2' = 0}^{j_2} \binom{j_2}{\ell'_2} \, \operatorname{sgn}(i_1 - j_1 - \ell_2').
    \end{align*}
The second and third equalities follow from the induction hypothesis applied with $j_1$ and $\bar{j}$, respectively. The fourth equality is obtained by reindexing, and the final equality follows from another application of Pascal's rule. All adjustments to the summation bounds involve only vanishing boundary terms as before.

\item Finally, we prove \cref{itm:lemma-indicators-sgn-2} by induction on $i_2$, with \cref{itm:lemma-indicators-sgn-1} as the base case. Indeed, after realising that
\begin{align}
    \mathbbm{1}_{i_1 + 1 < j_1 + \ell_2'} [i_1 + 2]_2 \, [j_1 + \ell_2']_2 = \mathbbm{1}_{i_1 < j_1 + \ell_2'} [i_1]_2 \, [j_1 + \ell_2']_2 ,\\
     \mathbbm{1}_{i_1 > j_1 + \ell_2'} [i_1]_2 \, [j_1 + \ell_2' + 1]_2 =  \mathbbm{1}_{i_1 \ge j_1 + \ell_2'} [i_1]_2 \, [j_1 + \ell_2' + 1]_2, 
\end{align}
by the same reasoning as before, \cref{itm:lemma-indicators-sgn-2} coincides with \cref{itm:lemma-indicators-sgn-1} for $i_2=1$:
    \begin{align*}
        &\sum_{\ell'_1=0}^{1} \sum_{\ell'_2=0}^{j_2}  \binom{1}{\ell'_1} \binom{j_2}{\ell'_2} \biggl( \mathbbm{1}_{i_1 + \ell_1'<j_1 + \ell_2'} [i_1 + \ell'_1 + 1]_2 \, [j_1 + \ell'_2]_2 -  \mathbbm{1}_{i_1 + \ell'_1 > j_1 + \ell_2'} [i_1+\ell'_1]_2 \, [j_1 + \ell_2' + 1]_2  \biggr) \\
        &=  \sum_{\ell'_2=0}^{j_2}  \binom{1}{0} \binom{j_2}{\ell'_2} \biggl( \mathbbm{1}_{i_1 <j_1 + \ell_2'} [i_1 + 1]_2 \, [j_1 + \ell'_2]_2 -  \mathbbm{1}_{i_1 > j_1 + \ell_2'} [i_1]_2 \, [j_1 + \ell_2' + 1]_2  \biggr) \\
        &\quad + \sum_{\ell'_2=0}^{j_2}  \binom{1}{1} \binom{j_2}{\ell'_2} \biggl( \mathbbm{1}_{i_1 + 1 <j_1 + \ell_2'} [i_1 + 2]_2 \, [j_1 + \ell'_2]_2 - \mathbbm{1}_{i_1 + 1 > j_1 + \ell_2'}  [i_1 + 1]_2 \, [j_1 + \ell_2' + 1]_2\biggr) \\
        &= \sum_{\ell'_2=0}^{j_2} \binom{j_2}{\ell'_2} \biggl( \mathbbm{1}_{i_1 <j_1 + \ell_2'} [i_1 + 1]_2 \, [j_1 + \ell'_2]_2 -  \mathbbm{1}_{i_1 \ge j_1 + \ell_2'} [i_1]_2 \, [j_1 + \ell_2' + 1]_2  \biggr) \\
        &\quad + \sum_{\ell'_2=0}^{j_2} \binom{j_2}{\ell'_2} \biggl( \mathbbm{1}_{i_1 < j_1 + \ell_2'} [i_1]_2 \, [j_1 + \ell'_2]_2 - \mathbbm{1}_{i_1 \ge j_1 + \ell_2'}  [i_1 + 1]_2 \, [j_1 + \ell_2' + 1]_2 \biggr) \\
        &= \sum_{\ell'_2=0}^{j_2} \binom{j_2}{\ell'_2} \biggl( \mathbbm{1}_{i_1 <j_1 + \ell_2'} [j_1 + \ell'_2]_2 -  \mathbbm{1}_{i_1 \geq j_1 + \ell_2'} [j_1 + \ell_2' + 1]_2  \biggr).
    \end{align*}
Now we fix $j_2\geq 1$ and continue with the induction on $i_2$. By \cref{itm:lemma-indicators-sgn-1}, the assertion holds for $i_2 = 1$. 

Suppose now that it holds for some $i_2 \ge 1$. To establish the case $i_2+1$, set $\bar{i} = i_1 + 1$ and $\bar{\ell} = \ell'_1 - 1$. Pascal's rule gives
    \begin{align*}
        &\sum_{\ell'_1=0}^{i_2 + 1} \sum_{\ell'_2=0}^{j_2}  \binom{i_2 + 1}{\ell'_1} \binom{j_2}{\ell'_2} \biggl( \mathbbm{1}_{i_1 + \ell_1'<j_1 + \ell_2'} [i_1 + \ell'_1 + 1]_2 \, [j_1 + \ell'_2]_2 -  \mathbbm{1}_{i_1 + \ell'_1> j_1 + \ell_2'} [i_1+\ell'_1]_2 \, [j_1 + \ell_2' + 1]_2  \biggr) \notag \\
        &= \sum_{\ell'_1=0}^{i_2} \sum_{\ell'_2=0}^{j_2}  \binom{i_2}{\ell'_1} \binom{j_2}{\ell'_2} \biggl( \mathbbm{1}_{i_1 + \ell_1'<j_1 + \ell_2'} [i_1 + \ell'_1 + 1]_2 \, [j_1 + \ell'_2]_2 -  \mathbbm{1}_{i_1 + \ell'_1> j_1 + \ell_2'} [i_1+\ell'_1]_2 \, [j_1 + \ell_2' + 1]_2  \biggr) \notag \\
        &\quad + \sum_{\ell'_1=0}^{i_2 + 1} \sum_{\ell'_2=0}^{j_2}  \binom{i_2}{\ell'_1 - 1} \binom{j_2}{\ell'_2} \biggl( \mathbbm{1}_{i_1 + \ell_1'<j_1 + \ell_2'} [i_1 + \ell'_1 + 1]_2 \, [j_1 + \ell'_2]_2 -  \mathbbm{1}_{i_1 + \ell'_1> j_1 + \ell_2'} [i_1+\ell'_1]_2 \, [j_1 + \ell_2' + 1]_2  \biggr)  \notag \\
        &= - \sum_{\ell'_1 = 0}^{i_2 - 1}  \sum_{\ell'_2 = 0}^{j_2 - 1}  \binom{i_2 - 1}{\ell'_1} \binom{j_2 - 1}{\ell'_2} \, \operatorname{sgn} (i_1 + \ell_1' - j_1 - \ell_2') \notag \\
        &\quad +  \sum_{\bar{\ell}=0}^{i_2} \sum_{\ell'_2=0}^{j_2}  \binom{i_2}{\bar{\ell}} \binom{j_2}{\ell'_2} \biggl( \mathbbm{1}_{\bar{i} + \bar{\ell}<j_1 + \ell_2'} [\bar{i} + \bar{\ell} + 1]_2 \, [j_1 + \ell'_2]_2 -  \mathbbm{1}_{\bar{i} + \bar{\ell} > j_1 + \ell_2'} [\bar{i} +\bar{\ell}]_2 \, [j_1 + \ell_2' + 1]_2  \biggr) \notag \\
        &= - \sum_{\ell'_1 = 0}^{i_2 - 1}  \sum_{\ell'_2 = 0}^{j_2 - 1}  \binom{i_2 - 1}{\ell'_1} \binom{j_2 - 1}{\ell'_2} \,  \operatorname{sgn} (i_1 + \ell_1' - j_1 - \ell_2') - \sum_{\bar{\ell} = 0}^{i_2 - 1}  \sum_{\ell'_2 = 0}^{j_2 - 1}  \binom{i_2 - 1}{\bar{\ell}} \binom{j_2 - 1}{\ell'_2} \, \operatorname{sgn} (\bar{i} + \bar{\ell} - j_1 - \ell_2')  \notag \\
        &= - \sum_{\ell'_1 = 0}^{i_2 }  \sum_{\ell'_2 = 0}^{j_2}  \binom{i_2 - 1}{\ell'_1} \binom{j_2 - 1}{\ell'_2} \, \operatorname{sgn} (i_1 + \ell_1' - j_1 - \ell_2') - \sum_{\ell'_1 = 0}^{i_2}  \sum_{\ell'_2 = 0}^{j_2 - 1}  \binom{i_2 - 1}{\ell'_1 - 1} \binom{j_2 - 1}{\ell'_2} \, \operatorname{sgn} (i_1 + \ell'_1 - j_1 - \ell_2') \notag \\
        &= - \sum_{\ell'_1 = 0}^{i_2}  \sum_{\ell'_2 = 0}^{j_2 - 1}  \binom{i_2}{\ell'_1} \binom{j_2 - 1}{\ell'_2} \, \operatorname{sgn} (i_1 + \ell_1' - j_1 - \ell_2')
    \end{align*}
The second and third equalities follow from the induction hypothesis, applied with $i_1$ and $\bar i$, respectively. The fourth equality is obtained by reindexing, and the final equality follows from another application of Pascal's rule. As before, adjustments to the summation bounds involve only vanishing boundary terms. This completes the proof.
\end{itemize}
\end{proof}

Let us now return to the proof of \cref{prop:simple-inverse-Kasteleyn}. The main difference of our simplifications compared to the original form in \cite{ayyerCorrelations2021} and to what was already done in \cite{ayyerGOE2023} is that we additionally use the negative binomial theorem for $x_2\ge 0$ and $|r|<1$ in the form
\begin{align}
    \label{eq:negative-binomial}
    \sum_{x_1=0}^{\infty} (-1)^{x_1} \binom{x_1 + x_2 -1}{x_1} r^{x_1} = (1+r)^{-x_2}.
\end{align}
The case $x_1 = x_2=0$ uses our convention $\binom{-1}{0}=1$. 

Let $\Gamma_0$ denote a positively oriented contour enclosing $0$ and no other integer. We will repeatedly use the standard contour integral representation of binomial coefficients,
\begin{align}
    \binom{s}{x} = \frac{1}{2\pi \mathrm{i}} \int\limits_{\Gamma_0} \mathrm{d}r \ \frac{(1+r)^{s}}{r^{x+1}}, \label{eq:contour-integral-binomial}
\end{align}
valid for integers $s,x$ satisfying $s\geq 0$ or $x\leq 0$. In particular, the case $s<0$ and $x=0$ relies on our convention $\binom{s}{0}=1$.

\begin{itemize}[leftmargin=*]
\item 
Suppose first that $\epsilon_i=\epsilon_j=1$. We start from \cite[Theorem 3.2]{ayyerCorrelations2021}. Substituting the formula for $t_n^{1,1}$ from \cite[(3.13)]{ayyerCorrelations2021} into \cite[(3.14)]{ayyerCorrelations2021}, and then using \cite[Lemma 5.4]{ayyerGOE2023} gives
\begin{align*}
    K_n^{-1} (x,y)
    &= \sum_{\ell'_1=0}^{i_1} \sum_{\ell'_2=0}^{j_1} (-1)^{i_2 + j_2} (-1)^{\ell'_1 + \ell'_2} \binom{i_2-1+\ell'_1}{\ell'_1} \binom{j_2-1+\ell'_2}{\ell'_2} \notag \\
    &\qquad \times \sum_{k_1 = 0}^n \sum_{k_2=0}^n \sum_{\ell_1=0}^{k_1} \sum_{\ell_2=0}^{k_2} p(n,k_1,\ell_1) \, p(n, k_2, \ell_2) \left( \mathbbm{1}_{\ell_2 = \ell_1+1} -  \mathbbm{1}_{\ell_2 = \ell_1-1}  \right) \notag \\
    & \qquad \times \frac{1}{(2\pi \mathrm{i})^2} \int\limits_{\Gamma_0} \mathrm{d}r_1 \int\limits_{\Gamma_0} \mathrm{d}r_2 \, \frac{(1+r_1)^{n-k_1} (1+r_2)^{n-k_2}}{r_1^{i_1 - \ell_1' - 2k_1 + 1} r_2^{j_1 - \ell_2' - 2k_2 + 1}}.
\end{align*}
Note that in the above expression, the binomial coefficients from \cite{ayyerCorrelations2021} and \cite{ayyerGOE2023} both use the convention that $\binom{a}{0} = 1$ for $a\in \mathbb{Z}$.

We now bring the $\ell_1'$- and $\ell_2'$-dependent factors into the contour integrals. The sums over $\ell_1'$ and $\ell_2'$ can then be extended to infinity, since the extra terms have zero residue at the origin. Indeed, for fixed $k_1$, the $r_1$-integral vanishes whenever the integrand is analytic at the origin, which happens whenever $\ell_1' \ge i_1-2k_1+1$. In particular, extending the $\ell_1'$-sum beyond $\ell_1'=i_1$ introduces only vanishing terms, and the same holds for the $\ell_2'$-sum. Since the convention $\binom{a}{0}=1$ for $a\in \mathbb{Z}$ ensures the applicability of \eqref{eq:negative-binomial}, we may evaluate both sums to obtain
\begin{align*}
    K_n^{-1} (x,y) & = \sum_{k_1 = 0}^n \sum_{k_2=0}^n \sum_{\ell_1=0}^{k_1} \sum_{\ell_2=0}^{k_2} p(n,k_1,\ell_1) \, p(n, k_2, \ell_2) \left( \mathbbm{1}_{\ell_2 = \ell_1+1} -  \mathbbm{1}_{\ell_2 = \ell_1-1}  \right) \notag \\
    & \qquad \times \frac{1}{(2\pi \mathrm{i})^2} \int\limits_{\Gamma_0} \mathrm{d}r_1 \int\limits_{\Gamma_0} \mathrm{d}r_2  \frac{(1+ r_1)^{n-k_1 -i_2}}{ r_1^{i_1 -2k_1 +1}} \frac{(1+r_2)^{n-k_2 + j_2}}{r_2^{j_1 -2k_2 +1}}.
\end{align*}
Rewriting the contour integrals as binomial coefficients using \eqref{eq:contour-integral-binomial}, we get
\begin{align}
    K_n^{-1} (x,y) &= (-1)^{i_2 + j_2} \sum_{k_1 = 0}^n \sum_{k_2=0}^n \sum_{\ell_1=0}^{k_1} \sum_{\ell_2=0}^{k_2}  p(n,k_1,\ell_1) \, p(n, k_2, \ell_2) \left( \mathbbm{1}_{\ell_2 = \ell_1+1} -  \mathbbm{1}_{\ell_2 = \ell_1-1}  \right) \notag \\
    & \qquad \times  \binom{n - k_1 - i_2}{i_1 - 2k_1} \binom{n - k_2 - j_2}{j_1 - 2k_2}.
    \label{eq:no_geometric_sum}
\end{align}
This proves the case $\epsilon_i = \epsilon_j = 1$ of \cref{prop:simple-inverse-Kasteleyn}, since the inner double sum equals $S_n(k_1, k_2)$ by \eqref{eq:S_n}. We note that the result is independent of the convention for ${a \choose 0}$ when $a < 0$. Indeed, the only way this convention comes into play for the first binomial coefficient is if $i_1 = 2k_1$ (hence $i_1$ is even) and $n - k_1 - i_2 = a< 0$. Due to the assumed range of $i_1$ and $i_2$ (see \eqref{eq:rangei_1-i_2}, in particular $i_2 \le n- \lfloor (i_1 + 1)/ 2\rfloor = n-i_1/2$, as $\epsilon_i = 1$ and $i_1$ is even) this is however impossible, i.e. no coefficient of the form ${a\choose 0}$ with $a<0$ ever occurs in this sum. The same argument applies to the second binomial coefficient.

\item 
Suppose now that $\epsilon_i=\epsilon_j=0$. Starting from \cite[Theorem 3.2]{ayyerCorrelations2021}, substituting the formula for $t_n^{0,0}$ from \cite[(3.10)]{ayyerCorrelations2021} as well as the formula for $h_n^{0,b}$ from \cite[(3.6)]{ayyerCorrelations2021} into \cite[(3.15)]{ayyerCorrelations2021}, and then using \cite[Lemma 5.4]{ayyerGOE2023} gives
\begin{align}
\label{eq:Kinv00}
    K_n^{-1} (x,y)
    &= \sum_{\ell'_1=0}^{i_2} \sum_{\ell'_2=0}^{j_2} (-1)^{i_2 + j_2} \binom{i_2}{\ell_1'} \binom{j_2}{\ell_2'}  \Biggl( \Bigg\{ \sum_{k_1 = 0}^n \sum_{k_2=0}^n \sum_{\ell_1=0}^{k_1} \sum_{\ell_2=0}^{k_2} p(n,k_1,\ell_1) \, p(n, k_2, \ell_2) \left( \mathbbm{1}_{\ell_2 = \ell_1+1} -  \mathbbm{1}_{\ell_2 = \ell_1-1}  \right) \notag \\
    &\qquad \times \frac{1}{(2\pi \mathrm{i})^2}  \int\limits_{\Gamma_0} \mathrm{d}r_1 \int\limits_{\Gamma_0} \mathrm{d}r_2 \, \frac{(1+r_1)^{n-k_1}(1+r_2)^{n-k_2}}{(1-r_1) \, r_1^{i_1 + \ell_1' - 2k_1} (1-r_2) \, r_2^{j_1 + \ell_2' - 2k_2}}  \Bigg\}\notag \\
    &\quad + [j_1 + \ell_2' + 1]_2 \sum_{k=0}^n p(n,k,0) \frac{1}{2\pi \mathrm{i}} \int \limits_{\Gamma_0} \mathrm{d}r \frac{(1+r)^{n-k}}{(1-r)\, r^{i_1 + \ell_1' - 2k}} \notag \\
    &\quad - [i_1 + \ell_1' + 1]_2 \sum_{k=0}^n p(n,k,0) \frac{1}{2\pi \mathrm{i}} \int \limits_{\Gamma_0} \mathrm{d}r \frac{(1+r)^{n-k}}{(1-r) \, r^{j_1 + \ell_2' - 2k}} \notag \\
     &\quad + \mathbbm{1}_{i_1 + \ell_1' < j_1 + \ell_2'} [i_1 + \ell_1' + 1]_2 \, [j_1 + \ell_2']_2 - \mathbbm{1}_{i_1 + \ell_1' > j_1 + \ell_2'} [i_1 + \ell_1' ]_2 \, [j_1 + \ell_2' +1]_2 \Biggr).
\end{align}
We now simplify the four contributions in \eqref{eq:Kinv00} separately, counting the two indicator terms in the last line as a single contribution. Throughout this discussion, each contribution is understood to include the common prefactor
\begin{align*}
    \sum_{\ell'_1=0}^{i_2} \sum_{\ell'_2=0}^{j_2} (-1)^{i_2 + j_2} \binom{i_2}{\ell_1'} \binom{j_2}{\ell_2'}
\end{align*}
multiplying the corresponding bracketed term in \eqref{eq:Kinv00}. 

The first contribution can be rewritten using the binomial theorem in the form
\begin{align}
    \sum_{\ell_1' = 0}^{i_2} \binom{i_2}{\ell_1'} \frac{1}{r^{\ell'_1}} 
    = \frac{(1+r)^{i_2}}{r^{i_2}}.
    \label{eq:binomial-thm}
\end{align}

Expanding the factors $(1-r_1)^{-1}$ and $(1-r_2)^{-1}$ as geometric series, and rewriting the resulting coefficients as binomial coefficients using \eqref{eq:contour-integral-binomial}, the first contribution becomes
\begin{align}
    &(-1)^{i_2 + j_2} \sum_{k_1 = 0}^n \sum_{k_2=0}^n \sum_{\ell_1=0}^{k_1} \sum_{\ell_2=0}^{k_2} p(n,k_1,\ell_1) \, p(n, k_2, \ell_2) \left( \mathbbm{1}_{\ell_2 = \ell_1+1} -  \mathbbm{1}_{\ell_2 = \ell_1-1}  \right) \notag \\
    &\quad \times \sum_{m_1=0}^{\infty} \binom{n-k_1+i_2}{i_1+ i_2 - 2k_1 - m_1 - 1} \sum_{m_2=0}^{\infty} \binom{n-k_2+j_2}{j_1+ j_2 - 2k_2 - m_2 - 1}.
    \label{eq:geometric-sum}
\end{align}

For the second contribution, \eqref{eq:binomial-thm} gives
\begin{align*}
    &\sum_{\ell_1' =0}^{i_2} \sum_{\ell_2' = 0}^{j_2} (-1)^{i_2 + j_2} \binom{i_2}{\ell_1'} \binom{j_2}{\ell_2'}  [j_1 + \ell_2'+ 1]_2 \sum_{k=0}^n p(n,k,0) \frac{1}{2\pi \mathrm{i}} \int\limits_{\Gamma_0} \mathrm{d}r \, \frac{(1 + r)^{n - k}}{(1 - r) \, r^{i_1 + \ell'_1 - 2k}} \\
    &= \sum_{\ell'_2 = 0}^{j_2} (-1)^{i_2 + j_2}\binom{j_2}{\ell'_2} [j_1 + \ell_2'+ 1]_2 \sum_{k=0}^n p(n,k,0) \frac{1}{2\pi \mathrm{i}} \int\limits_{\Gamma_0} \mathrm{d}r \, \frac{(1+r)^{n-k+i_2}}{(1-r) \, r^{i_1 + i_2 - 2k}}
\end{align*}
The factor $[j_1 + \ell_2'+ 1]_2$ selects precisely one parity class of $\ell_2'$, hence for $j_2 \ge 1$,
\begin{equation*}
    \sum_{\ell_2'=0}^{j_2}\binom{j_2}{\ell_2'}[j_1 + \ell_2' + 1]_2 = 2^{j_2-1}.
\end{equation*}
Expanding $(1-r)^{-1}$ as a geometric series, and rewriting the resulting coefficients as binomial coefficients using \eqref{eq:contour-integral-binomial} gives that the second contribution equals
\begin{align}\label{eq:second_00}
    (-1)^{i_2 + j_2} \, ([j_1 + 1]_2 \, \mathbbm{1}_{j_2 = 0} + 2^{j_2 - 1} \, \mathbbm{1}_{j_2 \ge 1}) \sum_{k=0}^n p(n,k,0) \sum_{m=0}^{\infty} \binom{n-k+i_2}{i_1 + i_2 - 2k - m - 1}.
\end{align}

The third contribution is the negative of the following
expression:
\begin{align*}
    &\sum_{\ell_1' =0}^{i_2} \sum_{\ell_2' = 0}^{j_2} (-1)^{i_2 + j_2} \binom{i_2}{\ell_1'} \binom{j_2}{\ell_2'} [i_1 + \ell_1' + 1]_2 \sum_{k=0}^n p(n,k,0) \frac{1}{2\pi \mathrm{i}} \int \limits_{\Gamma_0} \mathrm{d}r \frac{(1+r)^{n-k}}{(1-r) \, r^{j_1 + \ell_2' - 2k}} \\
    &= (-1)^{i_2 + j_2} \, ([i_1 + 1]_2 \, \mathbbm{1}_{i_2 = 0} + 2^{i_2 - 1} \, \mathbbm{1}_{i_2 \ge 1}) \sum_{k=0}^n p(n,k,0) \frac{1}{2\pi \mathrm{i}} \int\limits_{\Gamma_0} \mathrm{d}r \, \sum_{m=0}^{\infty} \frac{(1+r)^{n-k+j_2}}{r^{j_1 + j_2 - 2k - m}} \\
    &= (-1)^{i_2 + j_2} \, ([i_1 + 1]_2 \, \mathbbm{1}_{i_2 = 0} + 2^{i_2 - 1} \, \mathbbm{1}_{i_2 \ge 1}) \sum_{k=0}^n p(n,k,0) \sum_{m=0}^{\infty} \binom{n-k+j_2}{j_1 + j_2 - 2k - m - 1}.
\end{align*}
For the fourth contribution, we first factor out the common sign $(-1)^{i_2+j_2}$ and then apply \cref{lem:indicators-sgn} for the different ranges: \cref{itm:lemma-indicators-sgn-2} for $i_2,j_2 \ge 1$, \cref{itm:boundary-case2} for $i_2\ge1, j_2=0$, and \cref{itm:boundary-case1} for $i_2=0,j_2\ge 1$. Recall that for this lemma, the convention for the value of ${a \choose 0}$ for $a<0$ is irrelevant. For $i_2=j_2=0$, we just evaluate the sum directly. This gives

\begin{align*}
    &\sum_{\ell'_1=0}^{i_2} \sum_{\ell'_2=0}^{j_2} \binom{i_2}{\ell'_1} \binom{j_2}{\ell'_2} \biggl( \mathbbm{1}_{i_1 + \ell_1'<j_1 + \ell_2'} [i_1 + \ell'_1 + 1]_2 \, [j_1 + \ell'_2]_2 -  \mathbbm{1}_{i_1 + \ell'_1 > j_1 + \ell_2'} [i_1+\ell'_1]_2 \, [j_1 + \ell_2' + 1]_2  \biggr) \\
    &= - \sum_{\ell_1'=0}^{i_2-1}\sum_{\ell_2'=0}^{j_2-1} \binom{i_2-1}{\ell_1'}\binom{j_2-1}{\ell_2'} \operatorname{sgn}(i_1+\ell_1'-j_1-\ell_2') \notag \\
    &\quad +\mathbbm{1}_{i_2\ge 1}\mathbbm{1}_{j_2=0} \sum_{\ell_1'=0}^{i_2-1} \binom{i_2-1}{\ell_1'} \Big( \mathbbm{1}_{i_1+\ell_1'<j_1}[j_1]_2 - \mathbbm{1}_{i_1+\ell_1' \ge j_1}[j_1+1]_2 \Big) \notag \\
    &\quad +\mathbbm{1}_{i_2=0}\mathbbm{1}_{j_2\ge 1} \sum_{\ell_2'=0}^{j_2-1} \binom{j_2-1}{\ell_2'} \Big( \mathbbm{1}_{i_1\le j_1+\ell_2'}[i_1+1]_2 - \mathbbm{1}_{i_1>j_1+\ell_2'}[i_1]_2 \Big) \notag \\
    &\quad +\mathbbm{1}_{i_2=0}\mathbbm{1}_{j_2=0}\Bigl( \mathbbm{1}_{i_1<j_1}[i_1+1]_2 \, [j_1]_2 - \mathbbm{1}_{i_1>j_1}[i_1]_2 \, [j_1+1]_2 \Bigr).
    \end{align*}
    Under the standard convention that empty sums vanish, the indicators $\mathbbm{1}_{i_2\ge 1}$ in the second term and $\mathbbm{1}_{j_2\ge 1}$ in the third are redundant and hence omitted.
    
    Combining the four contributions yields \eqref{eq:K00Lemma}.

\item 
Suppose that $\epsilon_i=1, \epsilon_j=0$. After using \cite[Theorem 3.2]{ayyerCorrelations2021}, substituting the formula for $t_n^{1,0}$ from \cite[(3.11)]{ayyerCorrelations2021} as well as the formula for $h_n^{1,b}$ from \cite[(3.7)]{ayyerCorrelations2021} into \cite[(3.16)]{ayyerCorrelations2021}, and then using \cite[Lemma 5.4]{ayyerGOE2023} gives
\begin{align}
\label{eq:Kinv10}
    K_n^{-1} (x,y)
    &= \sum_{\ell'_1=0}^{i_1} \sum_{\ell'_2=0}^{j_2} (-1)^{i_2 + j_2} (-1)^{\ell_1'} \binom{i_2-1+\ell_1'}{\ell_1'} \binom{j_2}{\ell_2'} \notag \\
    &\quad \times \Biggl( \Biggl\{ \sum_{k_1 = 0}^n \sum_{k_2=0}^n \sum_{\ell_1=0}^{k_1} \sum_{\ell_2=0}^{k_2} p(n,k_1,\ell_1) \, p(n, k_2, \ell_2) \left( \mathbbm{1}_{\ell_2 = \ell_1+1} -  \mathbbm{1}_{\ell_2 = \ell_1-1}  \right) \notag \\
    &\qquad \times \frac{1}{(2\pi \mathrm{i})^2}  \int\limits_{\Gamma_0} \mathrm{d}r_1 \int\limits_{\Gamma_0} \mathrm{d}r_2 \, \frac{(1+r_1)^{n-k_1}}{r_1^{i_1 - \ell_1' - 2k_1 + 1}} \cdot \frac{(1+r_2)^{n-k_2}}{  (1-r_2) \, r_2^{j_1 + \ell_2' - 2k_2}} \Biggr\}  \notag \\
    &\qquad + [j_1 + \ell_2' + 1]_2 \sum_{k=0}^n p(n,k,0) \frac{1}{2\pi \mathrm{i}} \int \limits_{\Gamma_0} \mathrm{d}r \frac{(1+r)^{n-k}}{ r^{i_1 - \ell_1' - 2k +1}} \Biggr) \notag \\
     &\quad - (-1)^{i_2 + j_2} \mathbbm{1}_{i_1 \ge j_1 } \mathbbm{1}_{2i_1 + 2i_2 + 1 < 2j_1 + 2j_2 } \binom{j_2 - i_2 - 1}{i_1 - j_1}.
\end{align}
    The first term of this mixed case contains one binomial coefficient of the type appearing in the case $\epsilon_i=\epsilon_j=1$ and one of the type appearing in $\epsilon_i=\epsilon_j=0$. More precisely, in the double contour integral appearing in the first contribution above, 
   note that the $r_2$ fraction has in the denominator a factor $(1-r_2)$ which is handled using a geometric series, similar to \eqref{eq:geometric-sum}, whereas the $r_1$ fraction does not have such a term and is thus handled similarly to \eqref{eq:no_geometric_sum}. We skip the details as the calculations and arguments are identical. Likewise, the second term can be simplified using the arguments we used to obtain \eqref{eq:second_00} (but without a geometric series).
   
   Finally, consider the third and final contribution above. A priori the convention used for this binomial coefficient is ${a \choose 0 } = 1$ for $a \in \mathbb{Z}$. We claim that 
\begin{align*}
   \mathbbm{1}_{i_1 \ge j_1 } \mathbbm{1}_{2i_1 + 2i_2 + 1 < 2j_1 + 2j_2 } \binom{j_2 - i_2 - 1}{i_1 - j_1} =  \mathbbm{1}_{i_2<j_2} \binom{j_2 - i_2 - 1}{i_1 - j_1}.
\end{align*}
   Indeed, if $i_2\ge j_2$, then one immediately checks that indicators on both sides are equal to 0. On the other hand, if $i_2<j_2$, then the binomial coefficient cannot be of the form $\binom{a}{0}$ for $a<0$, hence it is nonzero if and only if  $j_2 -i_2 -1 \ge i_1 - j_1 \ge 0$, and this is equivalent (in this case) to both indicators on the left hand side being nonzero. The indicator $\mathbbm{1}_{i_2<j_2}$ ensures that the convention on $\binom{a}{0}$ for $a<0$ does not matter. This yields \eqref{eq:K10Lemma}.
\end{itemize}
\end{proof}

We next rewrite the entries in a form better suited to the TSSCPP events of interest. While not shorter, the resulting expressions contain factors of $(-1)^n$ which, combined with the factors $(-1)^k$ in $p(n,k,0)$, yield $(-1)^{n-k}$. As we shall see, for the vertices considered below, only terms with $n-k=\mathcal{O}(1)$ will end up contributing, which will be important since we will be interested in the limit as $n\to \infty$.

Now we can state the form of the entries of the inverse Kasteleyn matrix that we will use in this paper.
\begin{proposition}
\label{prop:final-invK}
Suppose $x=(x_1,x_2) = (i_1,i_1 + 2i_2 + \epsilon_i), \ y=(y_1,y_2) = (j_1,j_1 + 2j_2 + \epsilon_j)$ with $0 \le i_1,j_1 \le 2n-1$ and $\epsilon_i,\epsilon_j \in \{ 0,1 \}$, where $0 \le i_2 \le n - \left\lfloor \frac{i_1 + \epsilon_i}{2} \right\rfloor$ and $0 \le j_2 \le n - \left\lfloor \frac{j_1 + \epsilon_j}{2} \right\rfloor$.
\begin{itemize}[leftmargin=*]
    \item For $\epsilon_i=\epsilon_j = 1$, \eqref{eq:K11Lemma} remains unchanged.
    \item For $\epsilon_i=\epsilon_j = 0$,
    \begin{align}
        K_n^{-1} (x,y) 
        &=(-1)^{i_2+j_2} \Biggl( \sum_{k_1=0}^n \sum_{k_2=0}^n 
        S_n(k_1, k_2)  \sum_{m_1=0}^{\infty} \binom{n-k_1+i_2}{i_1+i_2-2k_1+m_1} \sum_{m_2=0}^{\infty} \binom{n-k_2+j_2}{j_1+j_2-2k_2+m_2} \notag \\
        &\quad +\mathbbm{1}_{i_2\ge 1} \mathbbm{1}_{j_2\ge 1} \Biggl[ -(-1)^n \, 2^{j_2-1} \sum_{k=0}^n p(n,k,0) \sum_{m=0}^{\infty} \binom{n-k+i_2}{i_1+i_2-2k+m} \notag \\
        &\hspace{2.57cm} + (-1)^n \, 2^{i_2-1} \sum_{k=0}^n p(n,k,0) \sum_{m=0}^{\infty} \binom{n-k+j_2}{j_1+j_2-2k+m} \notag \\
        &\hspace{2.57cm} - \sum_{\ell_1'=0}^{i_2-1}\sum_{\ell_2'=0}^{j_2-1} \binom{i_2-1}{\ell_1'}\binom{j_2-1}{\ell_2'} \operatorname{sgn}(i_1+\ell_1'-j_1-\ell_2') \Biggr] \notag \\ 
        &\quad +\mathbbm{1}_{i_2=0} \ \tilde \Xi^i_n(x,y)  +  \mathbbm{1}_{j_2=0} \tilde \Xi^j_n(x,y) + \mathbbm{1}_{i_2 = j_2=0} \tilde \Xi^{i,j}_n(x,y) \Biggr),
    \end{align}
    where $\tilde \Xi^i, \tilde \Xi^j$ and $\tilde \Xi^{i,j}$ are ``error'' terms explicitly given by 
    \begin{align}
    \tilde \Xi^i_n(x,y) 
        &: = (-1)^{n+1} \mathbbm{1}_{j_2\ge1} \sum_{k=0}^n 2^{j_2-1}p(n,k,0) \sum_{m=0}^{\infty} \binom{n-k}{i_1-2k+m} \notag \\
        &\quad + \left( \frac{(-1)^n\! -\!1}{2} \!+\! [i_1\!+\!1]_2 \right) \sum_{k=0}^n p(n,k,0) \sum_{m=0}^{\infty} \binom{n-k+j_2}{j_1+j_2-2k+m} \notag \\
        &\quad + \left(1+(-1)^n \right) 2^{j_2-2} \left( \mathbbm{1}_{j_2\ge1}- 2 \, [i_1+1]_2 \right) + \Xi^i_n(x,y), \\
    \tilde \Xi^j_n(x,y) 
        &:= (-1)^n \mathbbm{1}_{i_2\ge1} \sum_{k=0}^n 2^{i_2-1}p(n,k,0) \sum_{m=0}^{\infty} \binom{n-k}{j_1-2k+m} \notag \\
        &\quad + \left( \frac{1-(-1)^n}{2} -[j_1+1]_2 \right) \sum_{k=0}^n p(n,k,0) \sum_{m=0}^{\infty} \binom{n-k+i_2}{i_1+i_2-2k+m} \notag \\
        &\quad -\left( 1+(-1)^n \right) 2^{i_2-2} \left( \mathbbm{1}_{i_2\ge 1} - 2\, [j_1+1]_2 \right) + \Xi^j_n(x,y), \\
    \tilde \Xi^{i,j}_n(x,y)
        &:= \Xi^{i,j}_n(x,y). 
\end{align}
    
    \item For $\epsilon_i=1$ and $\epsilon_j = 0$,
    \begin{align}
        K_n^{-1} (x,y) 
                &= (-1)^{i_2 + j_2} \Biggl( -\sum_{k_1 = 0}^n \sum_{k_2 = 0}^n S_n(k_1, k_2) \binom{n-k_1 - i_2}{i_1 - 2k_1} \sum_{m=0}^{\infty} \binom{n-k_2 + j_2}{j_1 + j_2 - 2k_2 +m} \notag \\
                &\quad + \left( \mathbbm{1}_{j_2=0} \left( [j_1 + 1]_2 + \frac{(-1)^n -1}{2} \right) + \mathbbm{1}_{j_2 \ge 1} (-1)^n \, 2^{j_2 - 1} \right) \notag \\
                &\qquad \times \sum_{k=0}^n p(n,k,0) \binom{n-k-i_2}{i_1-2k} - \mathbbm{1}_{i_2<j_2} \binom{j_2 - i_2 - 1}{i_1 - j_1}
                \Biggr).
    \end{align}
    \item For $\epsilon_i=0$ and $\epsilon_j = 1$,
            \begin{align}
                K_n^{-1}(x,y) = - K_n^{-1}(y,x).
            \end{align}
\end{itemize}
\end{proposition}

\begin{proof}[Proof of \cref{prop:final-invK}]
We will need the following lemma.

\begin{lemma}
\label{lem:p-simplifications}
Let $S_n(k_1,k_2)$ be defined as in \eqref{eq:S_n}. The following identities hold:
    \begin{enumerate}
        \item \label{itm:simplified-2p}
        For every fixed $0 \le k_1 \le n$,
            \begin{align}
                &\sum_{k_2 = 0}^n 2^{n - k_2 + 1} S_n(k_1,k_2) = ((-1)^n - 1) p(n,k_1,0).
            \end{align}
            Moreover, by interchanging the roles of $k_1$ and $k_2$, we obtain for every fixed $0 \le k_2 \le n$,
            \begin{align}
                &\sum_{k_1 = 0}^n 2^{n - k_1 + 1} S_n(k_1,k_2) = (1-(-1)^n) p(n,k_2,0).
            \end{align}
        \item \label{itm:double-simplified-2p}
        The fully summed identity is
            \begin{align}
                \sum_{k_1=0}^n\sum_{k_2 = 0}^n 2^{n - k_1 + 1} 2^{n - k_2 + 1} S_n(k_1,k_2) =0.
            \end{align}
    \end{enumerate}
\end{lemma}
\begin{proof}
In the following, we repeatedly use \cite[Theorem 4.1]{ayyerCorrelations2021} in the form
    \begin{align}
        \sum_{k=0}^n 2^{n-k} p(n,k,0) = \frac{1+(-1)^n}{2}.
        \label{eq:Ayyer-Chhita-Thm4-1}
    \end{align}
\begin{enumerate}
    \item Fix $0\le k_1\le n$. By \eqref{eq:S_n} and \eqref{eq:p-ell},
        \begin{align*}
            &\sum_{k_2=0}^n \sum_{\ell_1=0}^{k_1} \sum_{\ell_2=0}^{k_2}  2^{n-k_2 + 1} p(n,k_1,\ell_1) p(n,k_2,\ell_2) \left( \mathbbm{1}_{\ell_2 = \ell_1 + 1} - \mathbbm{1}_{\ell_2 = \ell_1 - 1} \right) \\
            &= 2 \sum_{\ell_1=0}^{k_1} \sum_{\ell_2=0}^{n} \sum_{k_2=\ell_2}^n  2^{(n-\ell_2 ) - (k_2 - \ell_2)} p(n,k_1,\ell_1) (-1)^{\ell_2} p(n - \ell_2,k_2 - \ell_2,0) \left( \mathbbm{1}_{\ell_2 = \ell_1 + 1} - \mathbbm{1}_{\ell_2 = \ell_1 - 1} \right),
        \end{align*}
        where we reordered the sums in the last step. Setting $j=k_2-\ell_2$, renaming $j$ as $k_2$, and applying \eqref{eq:Ayyer-Chhita-Thm4-1} with parameter $n-\ell_2$ instead of $n$, we obtain
        \begin{align*}
            &2 \sum_{\ell_1=0}^{k_1} \sum_{\ell_2=0}^{n} \sum_{k_2=0}^{n-\ell_2}  2^{(n-\ell_2 ) - k_2} p(n,k_1,\ell_1) (-1)^{\ell_2} p(n - \ell_2,k_2,0) \left( \mathbbm{1}_{\ell_2 = \ell_1 + 1} - \mathbbm{1}_{\ell_2 = \ell_1 - 1} \right)\\
            &= 2 \sum_{\ell_1=0}^{k_1} \sum_{\ell_2=0}^{n} p(n,k_1, \ell_1) (-1)^{\ell_2} \frac{1+(-1)^{n-\ell_2}}{2} \left( \mathbbm{1}_{\ell_2 = \ell_1 + 1} - \mathbbm{1}_{\ell_2 = \ell_1 - 1} \right)\\
            &=  \sum_{\ell_1=0}^{\min \{ k_1,n-1 \}} p(n,k_1, \ell_1) (-1)^{\ell_1 + 1} \left( 1+(-1)^{n-\ell_1-1} \right) - \sum_{\ell_1=1}^{\min \{ k_1, n \}} p(n,k_1, \ell_1) (-1)^{\ell_1 - 1} \left( 1+(-1)^{n-\ell_1+1} \right) \\
            &= - p(n,k_1,0) \left( 1+(-1)^{n-1} \right) = \left( (-1)^n-1 \right) p(n,k_1,0).
        \end{align*}
        Indeed, the two sums in the third line cancel term by term for $1\le \ell_1 \le \min \{ k_1, n-1 \}$ since $(-1)^{\ell_1+1}=(-1)^{\ell_1-1}$ and $(-1)^{n-\ell_1-1}=(-1)^{n-\ell_1+1}$. Thus only the boundary term $\ell_1=0$ from the first sum remains. The possible boundary term $\ell_1=n$ from the second sum contributes zero, because its prefactor $(1+(-1)^{n-n+1})$ vanishes.

        The second identity in part \cref{itm:simplified-2p} follows by interchanging $k_1$ with $k_2$ and $\ell_1$ with $\ell_2$; the sign changes because $
\mathbbm{1}_{\ell_2=\ell_1+1}-\mathbbm{1}_{\ell_2=\ell_1-1}$ is antisymmetric under this interchange.
    \item Using \cref{itm:simplified-2p} together with \eqref{eq:Ayyer-Chhita-Thm4-1} gives
    \begin{align*}
        &\sum_{k_1=0}^n \sum_{k_2 = 0}^n \sum_{\ell_1=0}^{k_1} \sum_{\ell_2=0}^{k_2} 2^{n - k_1 + 1} 2^{n - k_2 + 1} p(n,k_1,\ell_1) p(n,k_2,\ell_2) \left( \mathbbm{1}_{\ell_2 = \ell_1 + 1} - \mathbbm{1}_{\ell_2 = \ell_1 - 1} \right) \\
        &= \sum_{k_1 = 0}^n  2^{n - k_1 + 1} \left( (-1)^n -1 \right) p(n,k_1,0)=0.
    \end{align*}
\end{enumerate}
This finishes the proof of \cref{lem:p-simplifications}.
\end{proof}

We now continue with the proof of \cref{prop:final-invK}. We will use in the following that, under the empty-sum convention,
\begin{align}
    \sum_{m=0}^{\infty} \binom{n - k + q_2}{q_1 + q_2 - 2k - m - 1} = 2^{n - k + q_2} - \sum_{m=0}^{\infty} \binom{n - k + q_2}{q_1 + q_2 - 2k +m},
    \label{eq:rewirte-binom-m}
\end{align}
for any integers $q_1,q_2 \ge 0$ and $0\le k\le n$. The sums over $m$ are finite by the binomial coefficient convention, and clearly the upper entries of the binomial coefficients are never negative. To see \eqref{eq:rewirte-binom-m}, set $a=n-k+q_2$ and $b=q_1 + q_2 -2k-1$. The sum on the left hand side is $\sum_{i=0}^{b} \binom{a}{i}$, while the sum on the right hand side is $\sum_{i=b+1}^a \binom{a}{i}$. Hence, \eqref{eq:rewirte-binom-m} follows from $\sum_{m=0}^a \binom{a}{m}=2^a$.

\begin{itemize}[leftmargin=*]
    \item Suppose that $\epsilon_i,\epsilon_j=0$. Starting from \eqref{eq:K00Lemma} and applying \eqref{eq:rewirte-binom-m} to all binomial coefficients of that form gives
    \begin{align}
        K_n^{-1}(x,y)
        &= (-1)^{i_2+j_2} \Biggl(
        \sum_{k_1=0}^n \sum_{k_2=0}^n 2^{n-k_1+i_2} 2^{n-k_2+j_2} S_n(k_1,k_2) \notag \\
        &\quad - \sum_{k_1=0}^n \sum_{k_2=0}^n 2^{n-k_1+i_2} S_n(k_1,k_2) \sum_{m=0}^{\infty} \binom{n-k_2+j_2}{j_1+j_2-2k_2+m} \notag \\
        &\quad - \sum_{k_1=0}^n \sum_{k_2=0}^n 2^{n-k_2+j_2} S_n(k_1,k_2) \sum_{m=0}^{\infty} \binom{n-k_1+i_2}{i_1+i_2-2k_1+m} \notag \\
        &\quad + \sum_{k_1=0}^n \sum_{k_2=0}^n S_n(k_1,k_2) \sum_{m_1=0}^{\infty} \sum_{m_2=0}^{\infty} \binom{n-k_1+i_2}{i_1+i_2-2k_1+m_1} \binom{n-k_2+j_2}{j_1+j_2-2k_2+m_2} \notag \\
        &\quad + \left( [j_1+1]_2 \mathbbm{1}_{j_2=0} + 2^{j_2-1} \mathbbm{1}_{j_2\ge 1} \right) \left( \sum_{k=0}^n 2^{n-k+i_2} p(n,k,0) - \sum_{k=0}^n \sum_{m=0}^{\infty} p(n,k,0) \binom{n-k+i_2}{i_1+i_2-2k+m} \right) \notag \\
        &\quad - \left( [i_1+1]_2 \mathbbm{1}_{i_2=0} + 2^{i_2-1} \mathbbm{1}_{i_2 \ge 1} \right) \left( \sum_{k=0}^n 2^{n-k+j_2} p(n,k,0) - \sum_{k=0}^n \sum_{m=0}^{\infty} p(n,k,0) \binom{n-k+j_2}{j_1+j_2 -2k+m} \right) \notag \\
        &\quad - \sum_{\ell_1'=0}^{i_2-1}\sum_{\ell_2'=0}^{j_2-1} \binom{i_2-1}{\ell_1'}\binom{j_2-1}{\ell_2'} \operatorname{sgn}(i_1+\ell_1'-j_1-\ell_2') \notag \\
        &\quad +\mathbbm{1}_{j_2=0} \Xi^j_n(x,y) 
        +\mathbbm{1}_{i_2=0} \Xi^i_n(x,y)  
+\mathbbm{1}_{i_2=0}\mathbbm{1}_{j_2=0} \Xi^{i,j}_n(x,y)
        \Biggr).
        \label{eq:proof-second-rewriting-K00}
    \end{align}
    Using \cref{lem:p-simplifications} together with \eqref{eq:Ayyer-Chhita-Thm4-1}, we now simplify the terms obtained above. The first term involving $S_n(k_1,k_2)$ is zero by \cref{itm:double-simplified-2p}. The fourth term involving $S_n(k_1,k_2)$ is the double-binomial (``main'') term appearing in the statement. The second and third terms involving $S_n(k_1,k_2)$ yield, after using \cref{itm:simplified-2p}, two terms involving $p(n,k,0)$, which combine with the $p(n,k,0)$-terms that were already present in \eqref{eq:proof-second-rewriting-K00}. For the sum
    \begin{equation*}
        \sum_{k=0}^n p(n,k,0) \sum_{m=0}^{\infty} \binom{n-k+i_2}{i_1+i_2-2k+m},
    \end{equation*}
    the total coefficient is
    \begin{equation*}
        -2^{j_2-1} ((-1)^n-1) - \left( [j_1 + 1]_2 \, \mathbbm{1}_{j_2=0} + 2^{j_2-1} \mathbbm{1}_{j_2\ge 1} \right) = \mathbbm{1}_{j_2=0} \left( \frac{1-(-1)^n}{2} - [j_1+1]_2 \right) - \mathbbm{1}_{j_2\ge 1} (-1)^n 2^{j_2-1}.
    \end{equation*}
    Similarly, for the sum
    \begin{equation*}
        \sum_{k=0}^n p(n,k,0) \sum_{m=0}^{\infty} \binom{n-k+j_2}{j_1+j_2-2k+m},
    \end{equation*}
    the total coefficient is 
    \begin{equation*}
        -2^{i_2-1} (1-(-1)^n) + \left( [i_1+1]_2 \, \mathbbm{1}_{i_2=0} + 2^{i_2-1} \mathbbm{1}_{i_2\ge 1} \right) = \mathbbm{1}_{i_2=0} \left( \frac{(-1)^n-1}{2} + [i_1+1]_2 \right) + \mathbbm{1}_{i_2\ge 1} (-1)^n 2^{i_2-1}.
    \end{equation*}
    Substituting these coefficients into \eqref{eq:proof-second-rewriting-K00}, using \eqref{eq:Ayyer-Chhita-Thm4-1} for the remaining $p(n,k,0)$-terms, and carrying over the remaining finite sums from the last four lines unchanged gives the claimed formula.

    \item Suppose that $\epsilon_i=1,\epsilon_j=0$. Starting from \eqref{eq:K10Lemma} and applying \eqref{eq:rewirte-binom-m} gives for $K_n^{-1}(x,y)$
    \begin{align*}
        &(-1)^{i_2+j_2} \Biggl( \sum_{k_1 = 0}^n \sum_{k_2 = 0}^n 2^{n-k_2+j_2} S_n(k_1,k_2) \binom{n-k_1-i_2}{i_1-2k_1} \notag \\
        &\quad -  \sum_{k_1 = 0}^n \sum_{k_2 = 0}^n S_n(k_1,k_2) \binom{n-k_1-i_2}{i_1-2k_1} \sum_{m=0}^{\infty} \binom{n-k_2+j_2}{j_1+j_2-2k_2+m} \notag \\
        &\quad + \left( [j_1+1]_2 \mathbbm{1}_{j_2=0} + 2^{j_2-1} \mathbbm{1}_{j_2\ge1} \right) \sum_{k=0}^n p(n,k,0) \binom{n-k-i_2}{i_1-2k} - \mathbbm{1}_{i_2<j_2} \binom{j_2-i_2-1}{i_1-j_1} \Biggr).
    \end{align*}
By \cref{lem:p-simplifications}\cref{itm:simplified-2p}, we obtain
   \begin{align*}
        &(-1)^{i_2+j_2} \Biggl( \mathbbm{1}_{j_2=0} \sum_{k = 0}^n  \frac{(-1)^n -1}{2} p(n,k,0) \binom{n-k-i_2}{i_1-2k} + \mathbbm{1}_{j_2\ge 1} \, 2^{j_2-1} \sum_{k=0}^n \left((-1)^n-1\right) p(n,k,0) \binom{n-k-i_2}{i_1-2k} \notag \\
        &\quad -  \sum_{k_1 = 0}^n \sum_{k_2 = 0}^n S_n(k_1,k_2) \binom{n-k_1-i_2}{i_1-2k_1} \sum_{m=0}^{\infty} \binom{n-k_2+j_2}{j_1+j_2-2k_2+m} \notag \\
        &\quad + \left( [j_1+1]_2 \mathbbm{1}_{j_2=0} + 2^{j_2-1} \mathbbm{1}_{j_2\ge1} \right) \sum_{k=0}^n p(n,k,0) \binom{n-k-i_2}{i_1-2k} - \mathbbm{1}_{i_2<j_2} \binom{j_2-i_2-1}{i_1-j_1} \Biggr) \\
        &= (-1)^{i_2 + j_2} \Biggl( -\sum_{k_1 = 0}^n \sum_{k_2 = 0}^n S_n(k_1,k_2) \binom{n-k_1 - i_2}{i_1 - 2k_1} \sum_{m=0}^{\infty} \binom{n-k_2 + j_2}{j_1 + j_2 - 2k_2 +m} \notag \\
        &\quad + \left( \mathbbm{1}_{j_2=0} \left( [j_1 + 1]_2 + \frac{(-1)^n -1}{2} \right) + \mathbbm{1}_{j_2 \ge 1} (-1)^n \, 2^{j_2 - 1} \right) \sum_{k=0}^n p(n,k,0) \binom{n-k-i_2}{i_1-2k} \notag \\
        &\quad - \mathbbm{1}_{i_2<j_2} \binom{j_2 - i_2 - 1}{i_1 - j_1}\Biggr).
    \end{align*}

\end{itemize}
This finishes the proof of \cref{prop:final-invK}.
\end{proof}

With the entries of the inverse Kasteleyn matrix in this form, we now explain why only those values of $k$ are relevant which ensure that $n-k=\mathcal{O}(1)$. Write a generic vertex as
\begin{equation*}
    (q_1,q_1+2q_2+\epsilon_q),
\end{equation*}
for $(q_1,q_2)\in \{ (i_1,i_2),(j_1,j_2) \}$ with the coordinate restrictions given above. In the preceding proposition, every term involving $p(n,k,0)$ is accompanied by a binomial coefficient of one of the forms
\begin{align*}
    \binom{n-k-q_2}{q_1 - 2k} \quad \text{or} \quad \binom{n-k+q_2}{q_1+q_2-2k+m},
\end{align*}
where $m\ge 0$. The first binomial coefficient enforces the constraint
\begin{align}
     0 \le q_1-2k \le n-k-q_2,
\end{align}
since, as mentioned before, $\binom{a}{0}$ for $a <0$ never appears. Hence, in particular
\begin{align}
   k \ge q_1 + q_2 - n.
\end{align}
Since the original summation is over $0\le k\le n$, this lower bound is sufficient to control the relevant terms. Indeed, if $q_1=2n-\mathcal{O}(1)$ and $q_2 = \mathcal{O}(1)$, i.e., the vertex is within microscopic distance $\mathcal{O}(1)$ from the centre of the hexagon, which is the case for the events considered here as we will see in the next section, then
\begin{align}
    n-\mathcal{O}(1) \le k \le n.
\end{align}
Thus, only $\mathcal{O}(1)$ values of $k$ contribute, and in particular $n-k=\mathcal{O}(1)$.

Similarly, the second binomial coefficient enforces
\begin{align}
    0 \le q_1 + q_2 -2k + m \le n-k + q_2,
    \label{eq:bounds}
\end{align}
hence in particular
\begin{align}
    k \ge q_1 + m -n \ge q_1 -n,
\end{align}
since $m\ge 0$. Again, since $0 \le k \le n$, if $q_1=2n-\mathcal{O}(1)$, then $n-\mathcal{O}(1) \le k \le n$. Hence only $\mathcal{O}(1)$ values of $k$ contribute, and $n-k=\mathcal{O}(1)$. Moreover, \eqref{eq:bounds} implies 
\begin{equation}
     m\leq n+k-q_1=\mathcal{O}(1).
\end{equation}
Thus, the number of contributing terms is bounded independently of $n$. In particular, this means that all $m-$sums with infinite upper limits appearing in \cref{prop:simple-inverse-Kasteleyn} have only finitely many non-zero terms.

\section{Pfaffian representation of nested-arch probabilities}
\label{sec:dimer-events}

\subsection{Nested-arch events as prescribed-dimer events}
In \cref{sec:fonseca-zinn-justin-correspondence}, we expressed nested submatching probabilities in terms of counts of TSSCPPs whose first $k$ paths are horizontal. We now show that this constraint is equivalent to the presence of prescribed dimers in the corresponding dimer configurations.

\begin{figure}[ht]
    \centering
	\includegraphics[width=0.28\textwidth]{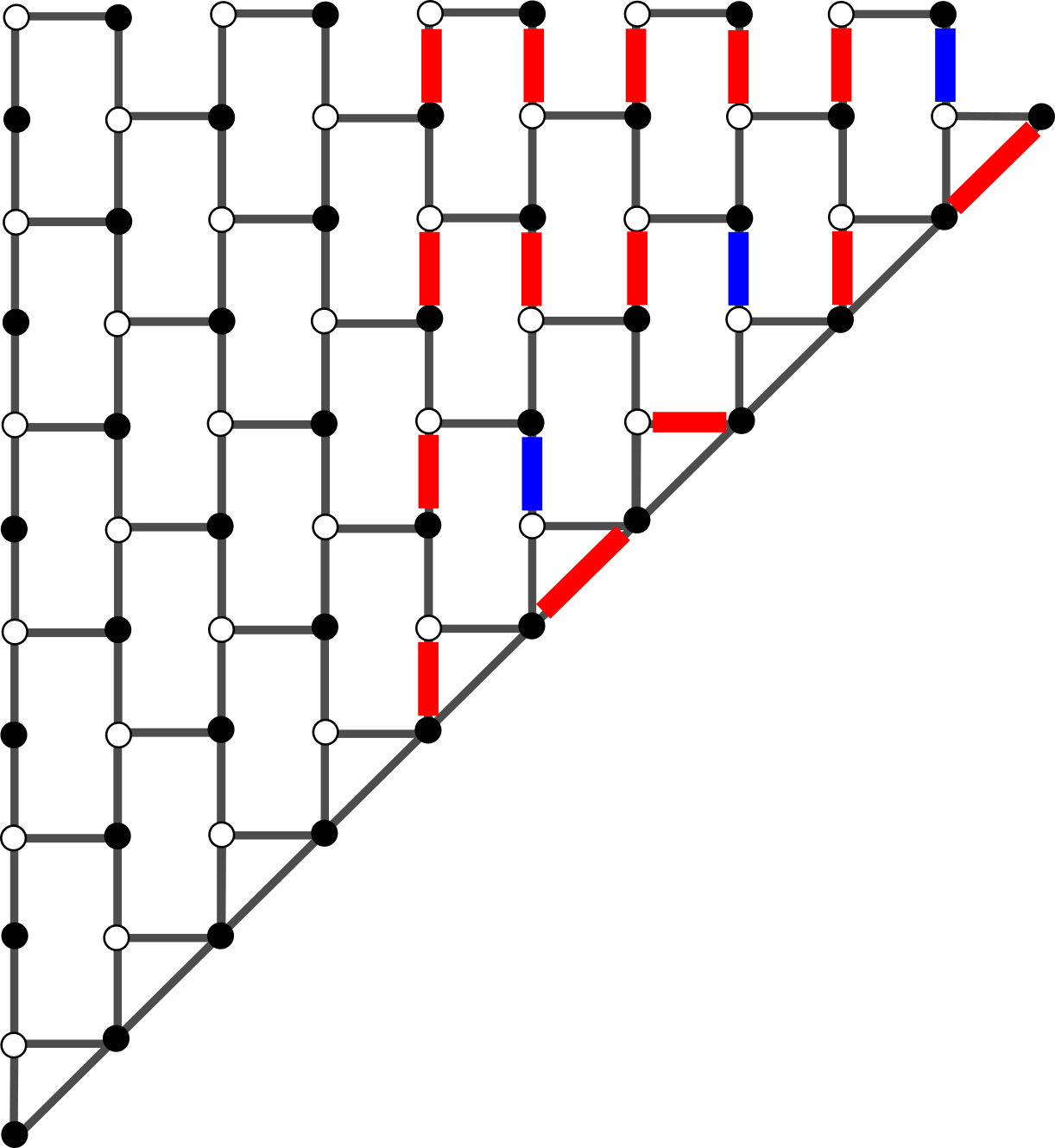}
         \caption{The prescribed dimers in $U_3^5$ (blue) and the additional dimers they force (red). The resulting pattern corresponds to the first three paths being horizontal.}
	\label{fig:enforced-dimers}
\end{figure}

Recall that $N=n+1$. For $1\leq k\leq n=N-1$ and $r=0,\ldots,k-1$, let
\begin{equation}
    e_r := \left( (2n-2r-1,\,2n-2r),(2n-2r-1,\,2n-2r+1) \right),
\end{equation}
and define
\begin{equation}
    U_k^n:=\{e_0,\ldots,e_{k-1}\}\subseteq E_n.
    \label{eq:prescribed-dimers}
\end{equation}
Thus, $U_k^n$ consists of $k$ prescribed edges of $\mathcal{G}_n$; see \cref{fig:enforced-dimers} for $U_3^5$. In the definition of $e_r$, the white endpoint is listed first and the black endpoint second.

\begin{lemma}
\label{lem:horizontal-paths-dimers}
    For $1\leq k\leq n$, under the bijection between TSSCPPs of order $N$ and dimer configurations on $\mathcal{G}_n$, the set $\mathcal{T}_{N,k}$ corresponds precisely to the dimer configurations containing $U_k^n$.
\end{lemma}

\begin{proof}
    Let $D_n \in \mathcal{M}(\mathcal{G}_n)$ be a dimer configuration containing $U_k^n$. For $r=0,\ldots,k-1$, define
    \begin{equation*}
        W_r^n:=\{ e_r+(0,2s) \mid 1\le s \le r \} \cup \{  e_r+(-1,2s) \mid 0\le s \le r \},
    \end{equation*}
    where $e+t$ denotes the edge obtained by adding the vector $t$ componentwise to each endpoint of $e$.

    We prove by induction on $r$ that imposing the dimers $e_0,\ldots,e_r$ forces every dimer in $W_r^n$, i.e., $W_r^n\subseteq D_n$.
    
    For $r=0$, the vertex $(2n-1,2n)$ is matched to $(2n-1,2n+1)$ by the prescribed dimer $e_0 \in U_k^n$. The only neighbours of $(2n-2,2n+1)$ are $(2n-1,2n+1)$ and $(2n-2,2n)$. It must therefore be matched to $(2n-2,2n)$, which proves that $W_0^n\subseteq D_n$.

    Now fix $1\leq r\leq k-1$ and assume that $W_{r-1}^n\subseteq D_{n}$. Define $a=2n-2r-1$, so that
    \begin{align*}
         e_r=((a,a+1),(a,a+2))\in U_k^n\subseteq D_n
    \end{align*}
    is the prescribed dimer. We first force the vertical dimers above $e_r$ in column $a$, proceeding upwards. For $s=1,\ldots,r$, the neighbours of $(a,a+2s+1)$ are $(a,a+2s), (a,a+2s+2)$ and $(a+1,a+2s+1)$. The lower neighbour is already matched by $e_r+(0,2s-2)$: this is the prescribed dimer $e_r$ when $s=1$, and the dimer forced at the preceding step when $s>1$. The right neighbour is already matched by $e_{r-1}+(-1,2s-2)\in W_{r-1}^n$. Thus $(a,a+2s+1)$ must be matched to $(a,a+2s+2)$. Consequently,
    \begin{align*}
        e_r+(0,2s)\in D_n, \qquad s=1,\ldots,r.
    \end{align*}
    We next force the required dimers in column $a-1$, proceeding downwards. For $s=r,r-1,\ldots,0$, consider the vertex $(a-1,a+2s+2)$. Its right neighbour $(a,a+2s+2)$ is already matched by $e_r+(0,2s)$. When $s=r$, the vertex lies on the top boundary and has no upper neighbour. When $s<r$, its upper neighbour is already matched by $e_r+(-1,2s+2)$, which was forced at the preceding step. Its only available neighbour is therefore $(a-1,a+2s+1)$. Hence
    \begin{align*}
        e_r+(-1,2s)\in D_n,     \qquad s=0,\ldots,r.
    \end{align*}
    Together, these two sweeps give $W_r^n\subseteq D_n$, completing the induction.

    Under the lozenge--path construction described in \cref{sec:fonseca-zinn-justin-correspondence}, the prescribed and forced dimers correspond to a tiling in which the first $k$ paths have no downward steps and are therefore horizontal. 
    
    Conversely, if the first $k$ paths are horizontal, the corresponding dimer configuration contains every edge of $U_k^n$.

    The prescribed dimers in $U_k^n$ also force additional dimers in the triangular part of $\mathcal{G}_n$, beyond those included in the sets $W_r^n$. These additional dimers are not needed for the preceding argument.
\end{proof}

The forcing mechanism is illustrated in \cref{fig:enforced-dimers} for $n=5$ and $k=3$.

\begin{corollary}
\label{cor:nested-submatching-pfaffian}
For $1\leq k\leq n$, let $D_n$ be uniformly distributed on $\mathcal{M}(\mathcal{G}_n)$. Write $e_r=(v_{2r+1}^n,v_{2r+2}^n)$ for $r=0,\ldots,k-1$, retaining the endpoint order in the definition of $e_r$. Then
    \begin{equation}
    \mathbb{P}\left( U_k^n \subseteq D_n \right) = \frac{|\mathcal{T}_{N,k}|}{|\mathcal{T}_N|} = \mathbb{P}_{\mathcal{C}_N}(F_k).
    \label{eq:three-equalities}
    \end{equation}
Moreover,
    \begin{align}
        \mathbb{P}_{\mathcal{C}_N}(F_k) &= \operatorname{Pf}\left( \left(K_n^{-1}(v_i^n,v_j^n)\right)_{1\leq i,j\leq 2k} \right).
    \label{eq:prob-nested-submatching-pfaffian}
    \end{align}
For $k=n$, this Pfaffian also gives $\mathbb{P}_{\mathcal{C}_{n+1}}(F_{n+1})$.

Furthermore, for each fixed $k\geq 1$,
    \begin{align}
        \mathbb{P}_{\mathbb{H}}(F_k) &=  \lim_{n\to\infty} \operatorname{Pf}\left( \left(K_n^{-1}(v_i^n,v_j^n)\right)_{1\leq i,j\leq 2k} \right).
   \label{eq:lim-prob-nested-submatching-pfaffian}
    \end{align}
\end{corollary}

\begin{proof}
By \cref{lem:horizontal-paths-dimers}, the TSSCPP--dimer bijection identifies $\mathcal{T}_{N,k}$ with the dimer configurations containing $U_k^n$. Since the bijection preserves the uniform measures, this gives the first equality in \eqref{eq:three-equalities}; the second follows from \cref{thm:fonseca-zinn-justin}.

For each prescribed edge, $v_{2r+1}^n$ is white and $v_{2r+2}^n$ is black, so $K_n(v_{2r+1}^n,v_{2r+2}^n)=-1$. Setting 
\begin{equation*}
    A=\left(K_n^{-1}(v_i^n,v_j^n)\right)_{1\leq i,j\leq 2k},
\end{equation*}
the Pfaffian correlation formula \eqref{eq:dimer-correlations} gives 
\begin{equation*}
    \mathbb{P}(U_k^n\subseteq D_n)=(-1)^k\operatorname{Pf}(A^{\mathsf T})=\operatorname{Pf}(A),
\end{equation*}
where the last equality follows from the skew-symmetry of $A$. This proves \eqref{eq:prob-nested-submatching-pfaffian}. For $k=n$, the additional equality follows from rotational invariance, as established in \eqref{eq:rotational-inv}. 

Finally, the convergence of cylindrical connectivity probabilities to their half-plane counterparts, as recalled in \cref{subsec:rationality-conj}, yields \eqref{eq:lim-prob-nested-submatching-pfaffian}.
\end{proof}

\subsection{Specialisation of the inverse Kasteleyn entries}
\label{subsec:relevant-inverse-K}

We now specialise the inverse Kasteleyn formulae of \cref{prop:final-invK} to the vertices appearing in \eqref{eq:prob-nested-submatching-pfaffian}. Recall the coordinate parametrisation
\begin{equation*}
    x=(i_1,i_1+2i_2+\epsilon_i), \qquad y=(j_1,j_1+2j_2+\epsilon_j),
\end{equation*}
where $\epsilon_i,\epsilon_j\in\{0,1\}$ indicate the vertex parities. 

We index the specialised entries by the cylinder parameter $N=n+1$. For $N>k\geq1$ and $r,s\in\{0,\ldots,k-1\}$, define
\begin{equation}
    f_N^{\epsilon_i,\epsilon_j}(r,s) :=K_{N-1}^{-1}\left( (2N-2r-3,\,2N-2r-1-\epsilon_i),(2N-2s-3,\,2N-2s-1-\epsilon_j) \right).
    \label{eq:relevant-kasteleyn-entries}
\end{equation}
These are precisely the inverse Kasteleyn entries corresponding to the endpoints of the prescribed edges in $U_k^{N-1}$; see \eqref{eq:prescribed-dimers}. In this parametrisation, the vertices satisfy
\begin{align*}
    &i_1=2N-2r-3,\qquad i_2=1-\epsilon_i,\\
    &j_1=2N-2s-3, \qquad j_2=1-\epsilon_j.
\end{align*}
Recall that $\epsilon_i$ and $\epsilon_j$ take the value $1$ at white vertices and $0$ at black vertices. Thus, each prescribed edge joins a white vertex nearest the diagonal to the black vertex immediately above it. For these vertices, several terms in \cref{prop:final-invK} simplify, and the error terms $\tilde{\Xi}^{i},\tilde{\Xi}^{j}$ and $\tilde{\Xi}^{i,j}$ vanish.

With this notation, the matrix in \eqref{eq:prob-nested-submatching-pfaffian} consists of $2\times 2$ blocks of the form
\begin{equation*}
    \begin{pmatrix}
        f_N^{1,1}(r,s) & f_N^{1,0}(r,s)\\
        f_N^{0,1}(r,s) & f_N^{0,0}(r,s)
    \end{pmatrix}, \qquad r,s\in\{0,\ldots,k-1\}.
\end{equation*}

We reindex the summation variables $k_1$ and $k_2$ in \cref{prop:final-invK} by setting $c_i=N-k_i$ for $i=1,2$. The resulting expressions involve the common inner sum
\begin{align}
    H_N(c_1,c_2)&:= \sum_{\ell_1=0}^{N-c_1} \sum_{\ell_2=0}^{N-c_2} p(N-1,N-c_1,\ell_1) \, p(N-1,N-c_2,\ell_2) \left( \mathbbm{1}_{\ell_2=\ell_1+1} - \mathbbm{1}_{\ell_2=\ell_1-1} \right) \label{eq:H-N-def} \\
    &= S_{N-1}(N-c_1,N-c_2), \notag
\end{align}
where $S_{n}$ is defined in \eqref{eq:S_n}. We extend $H_N(c_1,c_2)$ by zero whenever $c_1>N$ or $c_2>N$. Now we can state the final form of the relevant entries of the inverse Kasteleyn matrix.

\begin{lemma}
\label{lem:final-entries}
For $N>k\geq1$ and $r,s\in\{0,\ldots,k-1\}$, the specialised inverse Kasteleyn entries defined in \eqref{eq:relevant-kasteleyn-entries} satisfy
\begin{align}
    f_N^{1,1}(r,s) 
        &= \sum_{c_1=r+2}^{2r+2} \sum_{c_2=s+2}^{2s+2} H_N(c_1,c_2) \binom{c_1-1}{2r-c_1+2} \binom{c_2-1}{2s-c_2+2},
    \label{eq:finalf11}\\
        f_N^{0,0}(r,s)
        &=  \sum_{c_1 = 1}^{2r+2} \sum_{c_2 = 1}^{2s+2} H_N(c_1,c_2) \sum_{m_1=0}^{\infty} \binom{c_1}{2r-c_1-m_1+2} \sum_{m_2=0}^{\infty} \binom{c_2}{2s-c_2-m_2+2} \notag \\
        &\quad + (-1)^N \sum_{c_1=1}^{2r+2} p(N-1,N-c_1,0) \sum_{m = 0}^{\infty} \binom{c_1}{2r-c_1-m+2} \notag \\
        &\quad - (-1)^N \sum_{c_2=1}^{2s+2} p(N-1,N-c_2,0) \sum_{m = 0}^{\infty} \binom{c_2}{2s-c_2-m+2} - \operatorname{sgn} (s-r),
    \label{eq:finalf00}\\
    f_N^{1,0}(r,s) 
        &= \sum_{c_1 = r+2}^{2r+2} \sum_{c_2=1}^{2s+2} H_N(c_1,c_2) \binom{c_1-1 }{2r-c_1+2} \sum_{m=0}^{\infty} \binom{c_2}{2s-c_2-m+2} \notag \\
        &\quad + (-1)^N \sum_{c_1=r+2}^{2r+2} p(N-1,N-c_1,0) \binom{c_1-1}{2r-c_1+2} + \mathbbm{1}_{r = s},
        \label{eq:finalf10}\\
    f_N^{0,1}(r,s) &=-f_N^{1,0}(s,r).
    \label{eq:finalf01}
\end{align}
\end{lemma}

\begin{proof}
    Substituting the prescribed vertex coordinates into \cref{prop:final-invK}, using \eqref{eq:H-N-def} and using that by the symmetry of binomial coefficients $\binom{n}{k}=\binom{n}{n-k}$ gives this expression.
    
\end{proof}

For fixed $k$, the outer summation bounds in these expressions are independent of $N$, while the sums over $m$ are finite by the binomial coefficient convention. For $1\le c\le N$, using \eqref{eq:p-def}, one obtains that
\begin{align}
    (-1)^N p(N-1,N-c,0) 
    &= (-1)^c  (3c-1) C_{c-1} \frac{(2N-c)!\, (N+c-1)!}{(N-c)! \, (2N+c-1)!} \notag\\
    &= (-1)^c(3c-1)C_{c-1} \frac{\prod_{a=-c+1}^{c-1}(N+a)}{\prod_{a=-c+1}^{c-1}(2N+a)} \notag \\
    &= (-1)^c(3c-1)C_{c-1} \frac{1}{2} \prod_{a=1}^{c-1} \frac{N^2-a^2}{4N^2-a^2}.
    \label{eq:p-rational}
\end{align}
In the last step, the factors indexed by $a$ and $-a$ have been paired, while the central factors $N$ and $2N$ contribute the prefactor $1/2$. The equality between the first and last expressions also holds for $N<c\le 2N$, since both sides vanish by the zero-extension convention for $p$ and the factor at $a=N$ in the product, respectively. Hence, terms of this form are rational in $N$ for fixed $c$. Thus, it remains to establish the rational-function dependence of $H_N(c_1,c_2)$ on $N$. We address this using Gosper's and Zeilberger's algorithms in the following section.

\section{Hypergeometric evaluation and rationality}
\label{sec:hypergeomtric}
The aim of this section is to establish the following result.

\begin{proposition}
\label{prop:H_N-rational}
Let $N\ge 1$ and $1\le c_2 < c_1 \le 2N$ be integers. Then the function $H_N(c_1,c_2)$ defined in \eqref{eq:H-N-def} admits the representation    \begin{align}
    H_N(c_1,c_2)
        &= (-1)^{c_1-1} \frac{(3c_1-1)(3c_2-1)  (c_1-c_2)}{2(c_1+c_2)!}  \prod_{a=1}^{c_1-1} \frac{N^2-a^2}{4N^2-a^2} \Biggl((c_1+1) \prod_{k=c_1-2c_2+2}^{2c_1-c_2-2} k \notag \\
        &\quad + \sum_{j=1}^{c_2-1} (-1)^{j+1} \frac{C_{c_1-1} \, \left(j\, X_0(c_1,j) + N^2 X_1(c_1,j)\right)\,j\, (2c_1-c_2-2)! \, (c_1+j-2)! (2j)!}{4 (j+1)! \, (2c_1-j-2)! \, j! \, (2j-1)} \cdot \frac{1}{4N^2-j^2} \notag \\
        &\qquad \times \prod_{k=c_1-2c_2+2}^{c_1-2j-1} k \ \prod_{a=1}^{j-1} \frac{N^2-a^2}{4N^2-a^2} \Biggr),
        \label{eq:H_N}
\end{align}
for 
\begin{align}
    X_0(c_1,j)&:= 3c_1^3 (3j-1) + 3c_1^2 (3j^2-8j+3)+c_1(-17j^2+19j-6)-2j(2j-1)(j-1),\label{eq:X-0}\\
    X_1(c_1,j)&:=-12c_1^3 - 6c_1^2(4j-3) - 2c_1(12j^2-10j-1) - 2(3j+2)(2j-1)(j-1). \label{eq:X-1}
\end{align}
For $1\le c_1 < c_2 \le 2N$, the function is given by $H_N(c_1,c_2)=-H_N(c_2,c_1)$, with the right-hand side determined by the formula above. Furthermore, $H_N(c_1,c_2)=0$ when $c_1=c_2$. Thus, $H_N$ is antisymmetric.
\end{proposition}

To prove this, we will first put $H_N(c_1,c_2)$ into a suitable form in the next section, then recall the basics of Gosper's and Zeilberger's algorithms and finally apply these for evaluation and simplify the obtained expression.

\subsection{Reduction to a hypergeometric sum}
We begin by reducing $H_N(c_1,c_2)$ to a single sum. The indicator functions in \eqref{eq:H-N-def}, together with the convention that $p(n,k,\ell)=0$ outside $0\leq\ell\leq k\leq n$, give
\begin{align}
    H_N(c_1,c_2) &=\sum_{\ell=0}^{N-\mathrm{max}\{c_1,c_2\}} \Bigl( p(N-1,N-c_1,\ell) \, p(N-1,N-c_2,\ell+1) \notag \\
    &\hspace{2.5cm} - p(N-1,N-c_1,\ell+1) \, p(N-1,N-c_2,\ell)  \Bigl). \label{eq:H-ell}
\end{align}
In particular, $H_N(c_1,c_2)$ is antisymmetric in $c_1,c_2$, so $H_N(c_1,c_2)=0$ when $c_1=c_2$.

With the substitution $\omega=N-\ell$, \eqref{eq:p-ell} gives
\begin{equation*}
    p(N-1,N-c,N-\omega) = (-1)^{N-\omega}p(\omega-1,\omega-c,0).
\end{equation*}
Using \eqref{eq:p-def}, this becomes
\begin{align*}
p(N-1,N-c,N-\omega) &= (-1)^{N-c}(3c-1)C_{c-1} \frac{(2\omega-c)!(\omega+c-1)!}{(\omega-c)!(2\omega+c-1)!} \notag \\
                &=(-1)^{N-c}(3c-1)C_{c-1} \prod_{j=1}^{2c-1} \frac{\omega-c+j}{2\omega-c+j},
\end{align*}
for $1\leq c\leq\omega\leq N$. Substituting this expression into \eqref{eq:H-ell} and factoring the common product yields, for $c_1\neq c_2$,
\begin{align}
    H_N(c_1,c_2) = (-1)^{c_1+c_2}(3c_1-1)(3c_2-1) C_{c_1-1}C_{c_2-1} \sum_{\omega=\max\{c_1,c_2\}}^N W(\omega,c_1,c_2),
    \label{eq:H-W-relation}
\end{align}
where
\begin{align}
    W(\omega,c_1,c_2)&:=\prod_{a=1}^{2c_1-1} \frac{\omega-c_1+a}{2\omega-c_1+a} \prod_{a=1}^{2c_2-1} \frac{\omega-c_2+a}{2\omega-c_2+a} \notag \\
    &\qquad \times \Biggl( \frac{(\omega-c_2)(2\omega+c_2-2)(2\omega+c_2-1)}{(\omega+c_2-1)(2\omega-c_2-1)(2\omega-c_2)} - \frac{(\omega-c_1)(2\omega+c_1-2)(2\omega+c_1-1)}{(\omega+c_1-1)(2\omega-c_1-1)(2\omega-c_1)} \Biggr).
    \label{eq:W}
\end{align}
Here $W$ is defined for distinct positive integers $c_1,c_2$ and $\omega\geq\max\{c_1,c_2\}$. Its product factor is symmetric in $c_1,c_2$, while the difference in parentheses is antisymmetric. Hence
\begin{equation*}
    W(\omega,c_1,c_2)=-W(\omega,c_2,c_1).
\end{equation*}
It therefore suffices to consider $c_1>c_2$.

For integers $c_1>c_2\geq1$, define the shifted summand
\begin{equation}
    f(u,c_1,c_2):=W(u+c_1,c_1,c_2),
    \qquad u\geq0.
    \label{eq:f}
\end{equation}
We extend $f$ to the diagonal by setting $f(u,c,c)=0$ for $c\geq1$ and $u\geq0$.

With the change of summation index $u=\omega-c_1$, \eqref{eq:H-W-relation} becomes
\begin{equation}
    H_N(c_1,c_2) =(-1)^{c_1+c_2}(3c_1-1)(3c_2-1) C_{c_1-1}C_{c_2-1} \sum_{u=0}^{N-c_1}f(u,c_1,c_2), \qquad c_1>c_2\geq1,
    \label{eq:H-f-relation}
\end{equation}
where the remaining cases follow from the antisymmetry of $H_N$. Thus, the evaluation of $H_N$ reduces to the hypergeometric sum $\sum_{u=0}^{N-c_1}f(u,c_1,c_2)$, which we analyse in the following subsections.

\subsection{Hypergeometric summation methods}
\label{subsec:algorithms}

We briefly recall Gosper's and Zeilberger's algorithms, which we use to evaluate the hypergeometric sum arising in our analysis. A detailed account of these algorithms can be found, for instance, in \cite{petkovsek1996b}, which we use as our main reference.

A function $F(u)$ is called a \emph{hypergeometric term in $u$} if $F(u+1)/F(u)$ is a rational function of $u$, wherever the quotient is defined. Additional variables are treated as parameters.

Gosper's algorithm determines whether a hypergeometric term $F(u)$ admits a hypergeometric antidifference $G(u)$ satisfying
\begin{equation*}
    F(u) = G(u+1) - G(u). 
\end{equation*}
If such a term exists, the algorithm constructs it in the form $G(u)=R(u)F(u)$, where $R(u)$ is rational. Consequently,
\begin{equation*}
    \sum_{u=0}^n F(u) = G(n+1) - G(0).
\end{equation*}
The sum is thus reduced to boundary terms, which facilitates both its evaluation and the analysis of its asymptotic behaviour as $n\to\infty$. This will be useful in our subsequent analysis.

Zeilberger's algorithm, also known as \emph{creative telescoping}, extends this approach by seeking recurrences in an additional parameter. Suppose $F(u,k)$ is hypergeometric in both variables, meaning that $F(u+1,k)/F(u,k)$ and $F(u,k+1)/F(u,k)$ are rational functions of $u$ and $k$. The algorithm seeks polynomials $a_0(k),\ldots,a_J(k)$, not all zero, and a hypergeometric term $G(u,k)$ such that
\begin{align}
    \sum_{j=0}^J a_j(k) F(u,k+j) = G(u+1,k) - G(u,k).
    \label{eq:Zeilberger}
\end{align}
Thus, summing this identity over a finite range in $u$ yields a recurrence for the corresponding partial sums, with boundary terms arising from the telescoping expression. In our application, the resulting recurrence is of first order ($J=1$); solving it and evaluating the boundary terms will yield the desired closed form.

A sufficient condition for the existence of a relation of the form \eqref{eq:Zeilberger} is that $F$ is a \emph{proper hypergeometric term}; see \cite[Theorem 6.2.1]{petkovsek1996b}. We recall the definition.

\begin{definition}
\label{def:proper-hyp}
    A function $F(u,k)$ is called a \emph{proper hypergeometric term} if it can be written as
    \begin{align}
        F(u,k) = P(u,k) \frac{\prod_{i=1}^v (a_i u + b_i k + c_i)!}{\prod_{i=1}^w (f_i u + g_i k + h_i)!} x^u,
    \end{align}
    where $P(u,k)$ is a polynomial, $x$ is a parameter, the coefficients $a_i,b_i,c_i,f_i,g_i,h_i$ are fixed integers, and $v$ and $w$ are fixed, nonnegative integers.
\end{definition}

The following lemma verifies this condition for the summand arising in our problem.

\begin{lemma}
    For every fixed integer $c_1\ge 2$, the function $f(u,c_1,c_2)$ defined in \eqref{eq:f}, for integers $u\geq 0$ and $1\leq c_2<c_1$, is a proper hypergeometric term in $u$ and $c_2$.
    \label{lem:proper-hyp}
\end{lemma}

\begin{proof}
Recall that $f(u,c_1,c_2) = W(u+c_1,c_1,c_2)$, where $W(\omega,c_1,c_2)$ is defined in \eqref{eq:W}. For integers $1\le c\le \omega$,
\begin{equation*}
    \prod_{i=1}^{2c-1} \frac{\omega-c+i}{2\omega-c+i} = \frac{(\omega+c-1)!\, (2\omega-c)!}{(\omega-c)! \, (2\omega+c-1)!}.
\end{equation*}
Thus, after substituting $\omega=u+c_1$, the finite products in \eqref{eq:W} can be written as ratios of factorials whose arguments are linear expressions in $u$ and $c_2$, with constant terms that may depend on the fixed parameter $c_1$. 

The factor in parentheses in \eqref{eq:W}, under the same substitution, is a rational function of $u$ and $c_2$. After putting its two terms over a common denominator, its numerator is a polynomial, while its denominator is a product of linear expressions in $u$ and $c_2$. Each factor $L(u,c_2)$ in the denominator can be written in factorial form using
\begin{equation*}
    \frac{1}{L(u,c_2)} = \frac{(L(u,c_2)-1)!}{L(u,c_2)!},
\end{equation*}
since all these factors are positive integers for $u\geq 0$ and $1\leq c_2<c_1$. Hence $f(u,c_1,c_2)$ has the form required in \cref{def:proper-hyp}.
\end{proof}

\subsection{Evaluation of the hypergeometric sum}
\begin{proof}[Proof of \cref{prop:H_N-rational}]
We first consider $1\le c_2<c_1\le N$. Applying Zeilberger's algorithm to the summand $f(u,c_1,c_2)$ defined in \eqref{eq:f} yields the telescoping relation
\begin{align}
    A(c_1,c_2) f(u,c_1,c_2) + B(c_1,c_2) f(u,c_1,c_2+1) = G(u+1,c_1,c_2) - G(u,c_1,c_2), \label{eq:telescoping}
\end{align}
for
\begin{align}
    A(c_1,c_2) &:= -(c_1-2c_2)(c_1-2c_2+1)(c_1-c_2-1)(c_2+1)(c_1+c_2-1)(c_1+c_2), \\
    B(c_1,c_2) &:= -2(c_1-c_2)(2c_1-c_2-2)(c_1+c_2-1)(c_1+c_2)(c_1+c_2+1)(2c_2-1), \\
    G(u,c_1,c_2)&:= - \frac{D_1(u,c_1,c_2)}{D_2(u,c_1,c_2)},\\
    D_1(u,c_1,c_2)&:= (c_1-c_2-1)(c_1-c_2)\, c_2 \, u \, (c_1-c_2+u)(3c_1+2u-2)(3c_1+2u-1) (2c_1+c_2+2u-1) \notag \\
                &\quad \times \Bigl(12c_1^5+2(-3c_2+2c_2^2+1)\left(c_2^2+2(u-1)^2+3c_2(u-1)^2\right)  +6c_1^4(4c_2+4u-7) \notag \\
                &\qquad+c_1^3\left(15c_2^2-60u+12u^2+c_2(48u-65)+46\right) +c_1^2\Bigl((3c_2^3  +c_2^2(48u-34)\notag \\
                &\qquad-2(-16u+9u^2+5)+c_2(-88u+24u^2+49)\Bigr)+c_1\Bigl(c_2^3(24u-7) +c_2(28u-20u^2-2) \notag \\
                &\qquad-2(-6u+u^2+5)+c_2^2(-68u+24u^2+25)\Bigr)\Bigr)(u+1)^{\overline{2c_1-1}}(c_1-c_2+u+1)^{\overline{2c_2-1}}, \\
    D_2(u,c_1,c_2)&:=(2c_1+u-1)(c_1+c_2+u-1)(c_1+2u-1)(c_1+2u)(2c_1-c_2+2u-2) \notag \\
    &\quad \times (2c_1-c_2+2u-1)(2c_1-c_2+2u)(c_1+2u+1)^{\overline{2c_1-1}}(2c_1-c_2+2u+1)^{\overline{2c_2-1}},
\end{align}
for $u\ge 0$, where $x^{\overline{n}}$ denotes the rising factorial,
\begin{equation*}
    x^{\overline{n}}:=x(x+1)\cdots(x+n-1),
\end{equation*}
with $x^{\overline{0}}:=1$. Since $c_1>c_2\ge 1$ and $u\ge 0$, every factor in $D_2(u,c_1,c_2)$ is strictly positive, so $G(u,c_1,c_2)$ is well-defined. We will subsequently evaluate $G$ at $u=N-c_1+1\ge 1$.

The function $G$ is the telescoping term obtained from the rational certificate produced by Zeilberger's algorithm. The telescoping relation \eqref{eq:telescoping} is verified by direct substitution and simplification. The boundary case $c_2=c_1-1$ is included, since $A(c_1,c_1-1)=G(u,c_1,c_1-1)=0$ and $f(u,c_1,c_1)=0$. Note that $G(0,c_1,c_2)=0$. Thus, summing over $u=0,\dots,N-c_1$ gives
\begin{align}
    A(c_ 1,c_ 2) \sum_{u=0}^{N-c_1} f(u,c_1,c_2) + B(c_ 1,c_ 2) \sum_{u=0}^{N-c_1} f(u,c_1,c_2+1) = G(N-c_ 1+1, c_ 1,c_2),
    \label{eq:first-order-recurrence}
\end{align}
which is a first-order recurrence relation. Moreover, since $c_1> c_2 \ge1$, we have
\begin{equation*}
    B(c_1,j) \neq 0 \ \text{for } 1\le j \le c_2-1.
\end{equation*}
Hence, \eqref{eq:first-order-recurrence} may be divided by $B(c_1,j)$ and iterated from $j=1$ to $j=c_2-1$. This yields
\begin{align}
\label{eq:prae-hypergeometric-sum}
    \sum_{u=0}^{N-c_1} f(u,c_1,c_2) = \sum_{j=1}^{c_2 - 1} \frac{G(N-c_1 + 1, c_ 1, j)}{B(c_1,j)} \prod_{k=j+1}^{c_2 - 1} \left( -\frac{A(c_1,k)}{B(c_ 1,k)} \right) +\sum_{u=0}^{N-c_1} f(u,c_1,1) \prod_{k=1}^{c_2 - 1} \left( -\frac{A(c_1,k)}{B(c_1,k)} \right).
\end{align}
It remains to evaluate the initial value $\sum_{u=0}^{N-c_1}f(u,c_1,1)$. By \cref{lem:proper-hyp}, the specialisation $f(u,c_1,1)$ is proper hypergeometric in $u$, and hence, in particular, it is hypergeometric in $u$. Applying Gosper's algorithm with respect to $u$, we obtain a function $G_0$ satisfying
\begin{equation*}
    f(u,c_1,1) = G_0(u+1,c_1) - G_0(u,c_1)
\end{equation*}
for 
\begin{align}
 G_0(u,c_1):= \frac{u\,(3c_1+2u-2)(3c_1+2u-1)(u+1)^{\overline{2c_1-1}}}{2(2c_1+u-1)(c_1+2u-1)(c_1+2u)(c_1+2u+1)^{\overline{2c_1-1}}},
 \label{eq:Gosper-relation}
\end{align}
where the denominator is strictly positive since $u \ge 0$ and $c_1>1$. Direct substitution and simplification verify \eqref{eq:Gosper-relation}. Since $G_0(0,c_1)=0$, one obtains
\begin{align}
    \sum_{u=0}^{N-c_1} f(u,c_1,1) &= \sum_{u=0}^{N-c_1} \left( G_0(u+1,c_1) - G_0(u,c_1) \right) \notag \\
    &= G_0(N-c_1+1,c_1).
\end{align}
Putting this back into \eqref{eq:prae-hypergeometric-sum} yields
\begin{equation}
     \sum_{u=0}^{N-c_1} f(u,c_1,c_2) = \sum_{j=1}^{c_2 - 1} \frac{G(N-c_1 + 1, c_ 1, j)}{B(c_1,j)} \prod_{k=j+1}^{c_2 - 1} \left( -\frac{A(c_1,k)}{B(c_ 1,k)} \right) +G_0(N-c_1+1,c_1) \prod_{k=1}^{c_2 - 1} \left( -\frac{A(c_1,k)}{B(c_1,k)} \right).
     \label{eq:explicit-sum}
\end{equation}
Recall that empty sums are $0$ and empty products are $1$ by our conventions.

Although the expression in \eqref{eq:explicit-sum} is explicit, its dependence on $N$ is not yet transparent. We therefore simplify it, beginning with the product in the second term. Cancelling
common factors in $A(c_1,k)$ and $B(c_1,k)$ yields
\begin{equation}
    -\frac{A(c_1,k)}{B(c_1,k)} = - \frac{(c_1-2k)(c_1-2k+1)(c_1-k-1)(k+1)}{2(c_1-k)(2c_1-k-2)(c_1+k+1)(2k-1)}.
\end{equation}
The numerator factors satisfy
\begin{equation*}
    \prod_{k=1}^{c_2-1}(c_1-2k)(c_1-2k+1) = \prod_{k=c_1-2c_2+2}^{c_1-1} k.
\end{equation*}
This product vanishes when $c_2<c_1\leq 2c_2-2$. We retain the product form, since the corresponding factorial ratio would involve $(c_1-2c_2+1)!$, which is undefined in this range.

The remaining factors simplify as follows:
\begin{align*}
        \prod_{k=1}^{c_2-1} (2c_1-k-2) &=\frac{(2c_1-3)!}{(2c_1-c_2-2)!},\\
        \prod_{k=1}^{c_2-1} \frac{c_1-k-1}{c_1-k} &= \frac{c_1-c_2}{c_1-1},\\
        \prod_{k=1}^{c_2-1}\frac{1}{2k-1} &= \frac{2^{c_2-1}(c_2-1)!}{(2c_2-2)!},\\
        \prod_{k=1}^{c_2-1} \frac{1}{(c_1+k+1)} &= \frac{(c_1+1)!}{(c_1+c_2)!}, \\
        \prod_{k=1}^{c_2-1}(k+1)&=c_2!.
\end{align*}
The factor $2^{c_2-1}$ in the third equation cancels the factor $2^{c_2-1}$ contributed by the explicit $2$ in each denominator $B(c_1,k)$. Cancelling common factors, we thereby obtain
\begin{align*}
    \prod_{k=1}^{c_2-1} \left( -\frac{A(c_1,k)}{B(c_1,k)}\right) &= (-1)^{c_2-1} \frac{ (c_1-c_2) \, c_2! \, (2c_1-c_2-2)! \, (c_1+1)! \, (c_2-1)! \, \prod_{k=c_1-2c_2+2}^{c_1-1} k}{(c_1-1) (2c_1-3)! \, (c_1+c_2)! \, (2c_2-2)!} \\
    &= (-1)^{c_2-1} \frac{ (c_1-c_2) \, (2c_1-c_2-2)! \, (c_1+1)! \, \prod_{k=c_1-2c_2+2}^{c_1-1} k}{(c_1-1) (2c_1-3)! \, (c_1+c_2)! \, C_{c_2-1}},
\end{align*}
where $C_n = \frac{(2n)!}{n!(n+1)!}$ denotes the $n$-th Catalan number as before. To introduce $C_{c_1-1}$, we use
\begin{equation*}
    \frac{(c_1+1)!}{(c_1-1)(2c_1-3)!} = \frac{2(c_1+1)}{C_{c_1-1}(c_1-1)!}.
\end{equation*}
Since
\begin{equation*}
    \frac{(2c_1-c_2-2)!}{(c_1-1)!} = \prod_{k=c_1}^{2c_1-c_2-2}k,
\end{equation*}
whose range is consecutive with that of the preceding product, the two products combine to give
\begin{equation}
    \prod_{k=1}^{c_2-1} \left( -\frac{A(c_1,k)}{B(c_1,k)}\right) = (-1)^{c_2-1} \frac{2 (c_1-c_2) \, (c_1+1) }{(c_1+c_2)! \, C_{c_1-1} \, C_{c_2-1}} \prod_{k=c_1-2c_2+2}^{2c_1-c_2-2} k.
\end{equation}
For $\prod_{k=j+1}^{c_2-1} \left(-\tfrac{A(c_1,k)}{B(c_1,k)} \right)$, we repeat the steps from before, which gives
\begin{align*}
\prod_{k=j+1}^{c_2-1}(c_1-2k)(c_1-2k+1) &= \prod_{c_1-2c_2+2}^{c_1-2j-1} k,
\\
\prod_{k=j+1}^{c_2-1}(2c_1-k-2) &= \frac{(2c_1-j-3)!}{(2c_1-c_2-2)!},
\\
\prod_{k=j+1}^{c_2-1} \frac{c_1-k-1}{c_1-k} &= \frac{c_1-c_2}{c_1-j-1},
\\
\prod_{k=j+1}^{c_2-1}\frac{1}{2k-1}&=\frac{2^{c_2-j-1}(c_2-1)!(2j)!}
     {j!(2c_2-2)!},
\\
\prod_{k=j+1}^{c_2-1}\frac{1}{c_1+k+1}&=\frac{(c_1+j+1)!}{(c_1+c_2)!},
\\
\prod_{k=j+1}^{c_2-1}(k+1)&=\frac{c_2!}{(j+1)!}.
\end{align*}
The factor $2^{c_2-j-1}$ in the fourth equation cancels, just as before, the factor $2^{c_2-j-1}$ contributed by the explicit $2$ in each denominator $B(c_1,k)$. Hence,
\begin{align*}
    \prod_{k=j+1}^{c_2-1} \left( -\frac{A(c_1,k)}{B(c_1,k)}\right) &= (-1)^{c_2-j-1} \frac{ (c_1-c_2) \, c_2! \, (2c_1-c_2-2)! \, (c_1+j+1)! \, (c_2-1)! \ (2j)! \, \prod_{k=c_1-2c_2+2}^{c_1-2j-1} k}{(c_1-j-1) \, (j+1)! \, (2c_1-j-3)! \, (c_1+c_2)! \, j! \, (2c_2-2)!} \\
    &= (-1)^{c_2-j-1} \frac{ (c_1-c_2)  \, (2c_1-c_2-2)! \, (c_1+j+1)! \ (2j)! \, \prod_{k=c_1-2c_2+2}^{c_1-2j-1} k}{ (c_1-j-1) \, (j+1)! \, (2c_1-j-3)! \, (c_1+c_2)! \, j! \, C_{c_2-1}}
\end{align*}
Note that the last expression is well-defined since we always have $1 \le j\le c_2-1 < c_1-1$.

Before we continue, we also rewrite the $N$-dependent terms appearing in $\sum_{u=0}^{N-c_1}f(u,c_1,c_2)$. For $G_0(N-c_1+1,c_1)$, one obtains
\begin{align*}
    G_0(N-c_1+1,c_1)&=\frac{(N-c_1+1)(2N+c_1)(2N+c_1+1)(N-c_1+2)^{\overline{2c_1-1}}}{2(N+c_1)(2N-c_1+1)(2N-c_1+2)(2N-c_1+3)^{\overline{2c_1-1}}} \\
    &= \frac{\prod_{a=-c_1+1}^{c_1-1}(N+a)}{2\prod_{a=-c_1+1}^{c_1-1}(2N+a)} \\
    &= \frac{1}{4} \prod_{a=1}^{c_1-1} \frac{N^2-a^2}{4N^2-a^2},
\end{align*}
where the factors indexed by $a$ and $-a$ have been paired in the last step, while the central factors $N$ and $2N$ contribute the prefactor $1/2$.

Next, we consider the ratio $\frac{G(N-c_1+1,c_1,j)}{B(c_1,j)}$, using the polynomials $X_0$ and $X_1$ defined in \eqref{eq:X-0} and \eqref{eq:X-1}.

After substituting $u=N-c_1+1$ and cancelling common factors, we obtain
\begin{align*}
    \frac{G(N-c_1+1,c_1,j)}{B(c_1,j)} &= \frac{(c_1-j-1)\, j\, (c_1-N-1) (j-N-1) (c_1+2N)(c_1+2N+1)}{2(2c_1-j-2)(c_1+j-1)(c_1+j)(c_1+j+1)(2j-1)(c_1-2N-1)(j-2N-1)} \notag \\
    &\quad \times \frac{(j+2N+1) \, \left( j\, X_0(c_1,j) + N^2 \, X_1(c_1,j) \right)}{(j-2N)(c_1+N)(j+N)(c_1-2(N+1))(j-2(N+1))} \notag \\
    &\quad\times
    \frac{(-c_1+N+2)^{\overline{2c_1-1}}(-j+N+2)^{\overline{2j-1}}}{(-c_1+2N+3)^{\overline{2c_1-1}}(-j+2N+3)^{\overline{2j-1}}}.
\end{align*}
Rewriting all linear factors so that $N$ has a positive coefficient introduces an overall minus sign. Combining the linear factors with the corresponding rising factorials, we obtain
\begin{align}
    \frac{G(N-c_1+1,c_1,j)}{B(c_1,j)} &= -\frac{\left( j\, X_0(c_1,j) + N^2 \, X_1(c_1,j) \right)}{2}\frac{(c_1-j-1)\,j}{(2c_1-j-2)(c_1+j-1)(c_1+j)(c_1+j+1)(2j-1)} \notag \\
    &\quad \times \frac{\prod_{a=-c_1+1}^{c_1-1}(N+a)\prod_{a=-j+1}^{j-1}(N+a)}{\prod_{a=-c_1+1}^{c_1-1}(2N+a)\prod_{a=-j}^{j}(2N+a)} \notag \\
     &= - \frac{\left( j\, X_0(c_1,j) + N^2 \, X_1(c_1,j) \right)(c_1-j-1)\,j}{2^3(2c_1-j-2)(c_1+j-1)(c_1+j)(c_1+j+1)(2j-1)} \notag \\
    &\quad \times \frac{1}{4N^2-j^2} \prod_{a=1}^{c_1-1} \frac{N^2-a^2}{4N^2-a^2} \prod_{a=1}^{j-1} \frac{N^2-a^2}{4N^2-a^2},
\end{align}
where we paired the factors indexed by $a$ and $-a$ in the last step. 

Combining the preceding identities and cancelling common factors yields \eqref{eq:H_N} for $1\leq c_2<c_1\leq N$. If $N<c_1\leq 2N$, then $H_N(c_1,c_2)=0$ by \eqref{eq:H-N-def}, while the right-hand side of \eqref{eq:H_N} also vanishes because the product
$$
\prod_{a=1}^{c_1-1}\frac{N^2-a^2}{4N^2-a^2}
$$
contains a zero factor at $a=N$. All denominators remain nonzero in this range.
\end{proof}

\subsection{Proofs of the main theorems}
\label{subsec:proofs}
\begin{proof}[Proof of \cref{thm:cylinder-rationality}]
Fix $k \ge 1$ and let $N>k$. By \cref{cor:nested-submatching-pfaffian}, 
\begin{equation}
    \mathbb{P}_{\mathcal{C}_N} (F_k)
    = \operatorname{Pf}\left( \left( K_{N-1}^{-1}(v_i^{N-1},v_j^{N-1})\right)_{1\leq i,j\leq 2k} \right).
\end{equation}

By \cref{lem:final-entries}, the matrix entries are finite sums involving $H_N(c_1,c_2)$ and $(-1)^N p(N-1,N-c,0)$, with summation bounds independent of $N$ for fixed $k$. Since all indices satisfy $1 \le c,c_1,c_2\leq 2k<2N$, \cref{prop:H_N-rational} and \eqref{eq:p-rational} imply that each entry
is a rational function of $N$.

As the Pfaffian is a polynomial in its matrix entries with integer coefficients, there exists
$R_k\in\mathbb{Q}(N)$ such that
\begin{equation*}
    \mathbb{P}_{\mathcal{C}_N}(F_k)=R_k(N)
\end{equation*}
for every integer $N>k$.
\end{proof}
The explicit expressions in \eqref{eq:p-rational} and \cref{prop:H_N-rational} also show that $R_k$ is regular at $N=k$, since none of the denominator factors vanishes there. Thus, the boundary value $R_k(k)$ appearing in \cref{conj:cylinder-boundary} is well-defined.

\begin{proof}[Proof of \cref{thm:half-plane}]
Fix $k\ge 1$. 
By \cref{cor:nested-submatching-pfaffian} and \cref{thm:cylinder-rationality},
\begin{align}
    \mathbb{P}_{\mathbb{H}}(F_k) = \lim_{N\to \infty} \mathbb{P}_{\mathcal{C}_N} (F_k) = \lim_{N\to \infty} R_k (N).
\end{align}
Since $R_k \in \mathbb{Q}(N)$ is a rational function of $N$ and the limit exists, it is rational.
\end{proof}

The limiting entries of the inverse Kasteleyn matrix can also be determined explicitly. For fixed $c\ge 1$, $c_1>c_2\ge 1$, we have
\begin{align}
    \lim_{N\to \infty} (-1)^N p(N-1,N-c,0) &=  (-1)^c(3c-1)C_{c-1} \frac{1}{2^{2c-1}}, \label{eq:lim-p} \\
    H_{\infty}(c_1,c_2):=\lim_{N\to \infty} H_N(c_1,c_2) &=  (-1)^{c_1-1} \frac{(3c_1-1)(3c_2-1)  (c_1-c_2)}{2^{2c_1-1}(c_1+c_2)!} \Biggl((c_1+1) \prod_{k=c_1-2c_2+2}^{2c_1-c_2-2} k \notag \\
        &\ + \sum_{j=1}^{c_2-1} (-1)^{j+1} \frac{C_{c_1-1} \, X_1(c_1,j)\,j\, (2c_1-c_2-2)! \, (c_1+j-2)! (2j)!}{2^{2j+2} (j+1)! \, (2c_1-j-2)! \, j! \, (2j-1)} \prod_{k=c_1-2c_2+2}^{c_1-2j-1} k \Biggr). \label{eq:lim-H}
\end{align}
Substituting these expressions into the inverse Kasteleyn entries specified in \cref{lem:final-entries} yields the limits of the corresponding inverse Kasteleyn entries.

\section{Further observations on dyadic rationality}
\label{sec:observations}
\subsection{A related conjecture}

The denominator in \cref{conj:half-plane} admits a natural interpretation in terms of a rescaling of the limiting inverse Kasteleyn entries. For fixed $k\ge 1$, write
\begin{equation*}
    K^{-1}_{\infty} := \left( \lim_{N\to\infty}K_{N-1}^{-1}(v_i^{N-1},v_j^{N-1}) \right)_{1\le i,j\le 2k}.
\end{equation*}
Define
    \begin{align*}
        \lim_{N\to \infty} f^{\epsilon_i,\epsilon_j}_N (r,s) &:= f^{\epsilon_i,\epsilon_j}_{\infty}(r,s),\\
        \tilde{f}^{\epsilon_i,\epsilon_j}_{\infty}(r,s) &:= 2^{2r+2s+3}f^{\epsilon_i,\epsilon_j}_{\infty}(r,s),
    \end{align*}
and let $\tilde{K}^{-1}_{\infty}$ denote the corresponding rescaled $2k\times 2k$ matrix. Setting
\begin{equation*}
    B=\operatorname{diag}\left(
    2^{-2r-\frac{3}{2}}I_2
    \right)_{0\le r \le k-1},
\end{equation*}
where $I_2$ denotes the $2\times 2$ identity matrix, we have 
\begin{equation*}
    K^{-1}_{\infty} = B\tilde{K}^{-1}_{\infty}B^{\mathsf{T}}.
\end{equation*}
Moreover,
\begin{equation*}
    \det(B)
    =\prod_{r=0}^{k-1}
    \left(2^{-2r-\frac{3}{2}}\right)^2
    =2^{-2k^2-k}.
\end{equation*}
Hence, by the transformation formula for Pfaffians,\begin{align}
    \mathbb{P}_{\mathbb{H}}(F_k)
    &= \operatorname{Pf}(K^{-1}_{\infty})\notag\\
    &= \det(B)\operatorname{Pf}(\tilde{K}^{-1}_{\infty})\notag\\
    &= \frac{\operatorname{Pf}(\tilde{K}^{-1}_{\infty})}{2^{2k^2+k}}.
\end{align}

Consequently, \cref{conj:half-plane} is equivalent to the assertion that $\operatorname{Pf}(\tilde{K}^{-1}_{\infty})$ is an odd integer.

Exact symbolic computations in Mathematica for $1\le k\le 13$ show that the rescaled Pfaffian is an odd integer, although the only nonzero integer entries of the corresponding matrices are $\tilde{f}^{1,0}_{\infty}(0,0)=3$ and $\tilde{f}^{0,1}_{\infty}(0,0)=-3$. All remaining nonzero entries in these computations are non-integer dyadic rationals.

A possible intermediate step towards \cref{conj:dyadic-rat} is to establish dyadicity of the individual inverse Kasteleyn entries. Symbolic computations in Mathematica suggest the following conjecture.

\begin{conjecture}
\label{conj:weak-dyadicity}
    For integers $1\le c_2 < c_1$,
    \begin{equation}
        H_{\infty}(c_1,c_2) \in 2^{-2c_1-2c_2+2} \, \mathbb{Z}.
    \end{equation}
\end{conjecture}

By antisymmetry, the conjectured bound also holds when $c_1<c_2$, and $H_{\infty}(c_1,c_2)=0$ when $c_1=c_2$. Assuming \cref{conj:weak-dyadicity}, substituting these bounds together with \eqref{eq:lim-p} into the inverse Kasteleyn entries specified in \cref{lem:final-entries} yields
\begin{equation*}
    f^{\epsilon_i,\epsilon_j}_{\infty}(r,s) \in 2^{-4r-4s-6} \, \mathbb{Z}, \qquad \epsilon_i,\epsilon_j \in \{ 0,1\}.
\end{equation*}
Consequently,
\begin{equation*}
    \tilde{f}^{\epsilon_i,\epsilon_j}_{\infty}(r,s) \in 2^{-2r-2s-3} \, \mathbb{Z}, \qquad \epsilon_i,\epsilon_j \in \{ 0,1\}.
\end{equation*}

Each term in the Pfaffian expansion is a product of $k$ matrix entries in which every matrix index occurs exactly once, either as a row or as a column index. Since each block contains two matrix indices, the total denominator exponent of each term is at most
\begin{equation*}     
    \sum_{t=0}^{k-1}2\left(2t+\frac{3}{2}\right) =2k^2+k
\end{equation*}

Thus, 
\begin{equation*}
    \mathbb{P}_{\mathbb{H}} (F_k) \in 2^{-4k^2-2k}\mathbb{Z}.
\end{equation*}
This would establish the dyadic rationality predicted by Romik, see \cref{conj:dyadic-rat}, in the case of nested-arch events, although not the exact denominator conjectured in \cref{conj:half-plane}.

\subsection{Partial results and limitations of an entrywise approach}
The known value
\begin{equation*}
    \mathbb{P}_{\mathbb{H}}(F_1)=\frac{3}{2^3}
\end{equation*}
provides a simple first case in which to examine whether dyadicity can be established directly from the inverse Kasteleyn entries, without determining the denominator in lowest terms.

For $k=1$, we have $r,s=0$, so the relevant entries of the inverse Kasteleyn matrix as specified in \cref{lem:final-entries} involve only $H_{\infty}(2,1)$, up to antisymmetry and vanishing diagonal terms. It therefore suffices to establish \cref{conj:weak-dyadicity} for $c_1=2$ and $c_2=1$, which would imply
\begin{equation*}
    \mathbb{P}_{\mathbb{H}} (F_1) \in 2^{-6}\mathbb{Z}.
\end{equation*}

For $c_2=1$, the second sum in \cref{prop:H_N-rational} vanishes. Simplifying the remaining expression yields, for $1< c_1 \le 2N$,
\begin{align*}
    H_N(c_1,1) &= (-1)^{c_1-1} \frac{(3c_1-1)  (c_1-1) \, (c_1+1) }{(c_1+1)!}   \prod_{a=1}^{c_1-1} \, \frac{N^2-a^2}{4N^2-a^2} \prod_{k=c_1}^{2c_1-3} k \notag \\
    &= (-1)^{c_1-1} \frac{(3c_1-1)(2c_1-2)!}{2 c_1! (c_1-1)!} \prod_{a=1}^{c_1-1} \, \frac{N^2-a^2}{4N^2-a^2} \notag \\
    &= (-1)^{c_1-1} \frac{(3c_1-1)C_{c_1-1}}{2} \prod_{a=1}^{c_1-1} \,  \frac{N^2-a^2}{4N^2-a^2}.
\end{align*}
Taking $N\to\infty$, we obtain
\begin{equation*}
    H_{\infty}(c_1,1)
    =(-1)^{c_1-1}
    \frac{(3c_1-1)C_{c_1-1}}{2^{2c_1-1}}
    \in 2^{-2c_1}\mathbb{Z}.
\end{equation*}
Thus, \cref{conj:weak-dyadicity} holds for $c_2=1$ and arbitrary $c_1\geq 2$. In particular,
\begin{equation*}
    H_{\infty}(2,1)=-\frac{5}{2^3}\in 2^{-4}\mathbb{Z},
\end{equation*}
which establishes the required bound for the case $k=1$.

For $k=2$, the relevant parameters satisfy $1\leq c_2<c_1\leq 4$. While these individual cases can be evaluated directly, proving \cref{conj:weak-dyadicity} in general requires bounds valid for arbitrary parameters. We therefore sought analogous expressions for fixed $c_2$ and arbitrary $c_1>c_2$.

The case $c_2=2$ can also be treated by a direct, although substantially more involved, calculation. This requires expressing the arising polynomials in the basis of falling factorials, identifying coefficient sequences with the aid of the \textit{On-Line Encyclopedia of Integer Sequences}, and verifying the necessary divisibility properties. We omit the details, as the argument is specific to this case and does not suggest an effective extension to arbitrary $c_2$. These results provide further evidence for \cref{conj:weak-dyadicity}, but the rapidly increasing complexity suggests that a different approach may be needed to establish it in general.

\appendix
\section{Explicit inverse Kasteleyn entries for nested-arch events}
\label{appendix-A}
We provide explicit values of the limiting inverse Kasteleyn entries used to compute the probabilities $\mathbb{P}_{\mathbb{H}}(F_k)$, $1\leq k\leq 8$, listed in \cref{tab:half-plane-probabilities}.

For $\epsilon_i,\epsilon_j\in\{0,1\}$, define 
\begin{equation*}    
    f_\infty^{\epsilon_i,\epsilon_j}(r,s):= \lim_{N\to\infty}f_N^{\epsilon_i,\epsilon_j}(r,s), 
\end{equation*} 
where the entries $f_N^{\epsilon_i,\epsilon_j}(r,s)$ are given in \cref{lem:final-entries}. With the ordering used in \eqref{eq:lim-prob-nested-submatching-pfaffian}, the limiting inverse Kasteleyn entries appearing in the Pfaffian are arranged in $2\times 2$ blocks of the form
\begin{equation*}
    \begin{pmatrix}
        f_\infty^{1,1}(r,s) & f_\infty^{1,0}(r,s)\\
        f_\infty^{0,1}(r,s) & f_\infty^{0,0}(r,s)
    \end{pmatrix}, \qquad r,s\in\{0,\ldots,k-1\}. 
\end{equation*}
Since $f^{0,1}_{\infty}(r,s)=-f^{1,0}_{\infty}(s,r)$, it suffices to evaluate the remaining three functions.

\newpage 
The values of $f_\infty^{1,1}(r,s)$ for $0\le r,s \le 7$ are:

\[
\scalebox{1.2}{$
\begin{array}{c|cccc}
r \backslash s & 0 & 1 & 2 & 3 \\ \hline
0 & 0
  & \frac{259}{2^{10}}
  & -\frac{589}{2^{14}}
  & -\frac{53789}{2^{18}} \\[7pt]
1 & -\frac{259}{2^{10}}
  & 0
  & \frac{8035}{2^{17}}
  & \frac{778881}{2^{22}} \\[7pt]
2 & \frac{589}{2^{14}}
  & -\frac{8035}{2^{17}}
  & 0
  & \frac{115375}{2^{25}} \\[7pt]
3 & \frac{53789}{2^{18}}
  & -\frac{778881}{2^{22}}
  & -\frac{115375}{2^{25}}
  & 0 \\[7pt]
4 & -\frac{244243}{2^{20}}
  & \frac{7486053}{2^{25}}
  & -\frac{104771957}{2^{30}}
  & -\frac{478161617}{2^{33}} \\[7pt]
5 & \frac{1657289}{2^{25}}
  & -\frac{57171995}{2^{30}}
  & \frac{61871849}{2^{30}}
  & -\frac{7869930147}{2^{37}} \\[7pt]
6 & \frac{164744183}{2^{30}}
  & -\frac{296362367}{2^{31}}
  & \frac{222268449}{2^{34}}
  & \frac{117906079935}{2^{40}} \\[7pt]
7 & -\frac{3249113203}{2^{34}}
  & \frac{47944916563}{2^{38}}
  & -\frac{92501363095}{2^{41}}
  & -\frac{6117878764101}{2^{46}}
\end{array}
$}
\]
\\
\[
\scalebox{1.2}{$
\begin{array}{c|cccc}
r \backslash s & 4 & 5 & 6 & 7 \\ \hline
0 & \frac{244243}{2^{20}}
  & -\frac{1657289}{2^{25}}
  & -\frac{164744183}{2^{30}}
  & \frac{3249113203}{2^{34}} \\[7pt]
1 & -\frac{7486053}{2^{25}}
  & \frac{57171995}{2^{30}}
  & \frac{296362367}{2^{31}}
  & -\frac{47944916563}{2^{38}} \\[7pt]
2 & \frac{104771957}{2^{30}}
  & -\frac{61871849}{2^{30}}
  & -\frac{222268449}{2^{34}}
  & \frac{92501363095}{2^{41}} \\[7pt]
3 & \frac{478161617}{2^{33}}
  & \frac{7869930147}{2^{37}}
  & -\frac{117906079935}{2^{40}}
  & \frac{6117878764101}{2^{46}} \\[7pt]
4 & 0
  & \frac{34788948289}{2^{40}}
  & \frac{2079692100391}{2^{44}}
  & -\frac{75847897150279}{2^{49}} \\[7pt]
5 & -\frac{34788948289}{2^{40}}
  & 0
  & \frac{2312952191759}{2^{47}}
  & \frac{843616745718929}{2^{53}} \\[7pt]
6 & -\frac{2079692100391}{2^{44}}
  & -\frac{2312952191759}{2^{47}}
  & 0
  & \frac{2548451627477271}{2^{56}} \\[7pt]
7 & \frac{75847897150279}{2^{49}}
  & -\frac{843616745718929}{2^{53}}
  & -\frac{2548451627477271}{2^{56}}
  & 0
\end{array}
$}
\]

\newpage
The corresponding values of $f_\infty^{1,0}(r,s)$ are:

\[
\scalebox{1.2}{$
\begin{array}{c|cccc}
r \backslash s & 0 & 1 & 2 & 3 \\ \hline
0 & \frac{3}{2^3}
  & \frac{131}{2^{10}}
  & -\frac{2925}{2^{14}}
  & \frac{28771}{2^{18}} \\[7pt]
1 & \frac{325}{2^{10}}
  & \frac{77}{2^8}
  & \frac{25187}{2^{17}}
  & -\frac{491647}{2^{22}} \\[7pt]
2 & -\frac{1155}{2^{14}}
  & \frac{14973}{2^{17}}
  & \frac{22881}{2^{16}}
  & \frac{4187823}{2^{25}} \\[7pt]
3 & -\frac{30779}{2^{18}}
  & -\frac{162113}{2^{22}}
  & \frac{4720769}{2^{25}}
  & \frac{1411707}{2^{22}} \\[7pt]
4 & \frac{153217}{2^{20}}
  & \frac{592997}{2^{25}}
  & -\frac{125416757}{2^{30}}
  & \frac{1370502767}{2^{33}} \\[7pt]
5 & -\frac{1176783}{2^{25}}
  & \frac{2765925}{2^{30}}
  & \frac{29900107}{2^{30}}
  & -\frac{10119019683}{2^{37}} \\[7pt]
6 & -\frac{97312169}{2^{30}}
  & -\frac{41977719}{2^{31}}
  & \frac{1158286737}{2^{34}}
  & -\frac{32271513921}{2^{40}} \\[7pt]
7 & \frac{1968860149}{2^{34}}
  & \frac{5425440659}{2^{38}}
  & -\frac{189050123815}{2^{41}}
  & \frac{4691813654715}{2^{46}}
\end{array}
$}
\]
\\
\[
\scalebox{1.2}{$
\begin{array}{c|cccc}
r \backslash s & 4 & 5 & 6 & 7 \\ \hline
0 & \frac{39323}{2^{20}}
  & -\frac{4245321}{2^{25}}
  & \frac{93006793}{2^{30}}
  & \frac{494268275}{2^{34}} \\[7pt]
1 & -\frac{874165}{2^{25}}
  & \frac{123380123}{2^{30}}
  & -\frac{174926325}{2^{31}}
  & -\frac{6711279187}{2^{38}} \\[7pt]
2 & -\frac{50800971}{2^{30}}
  & -\frac{16392613}{2^{30}}
  & \frac{477529311}{2^{34}}
  & -\frac{12527070825}{2^{41}} \\[7pt]
3 & \frac{1187147473}{2^{33}}
  & -\frac{12621578781}{2^{37}}
  & \frac{33220324593}{2^{40}}
  & \frac{2145792671301}{2^{46}} \\[7pt]
4 & \frac{1407118661}{2^{32}}
  & \frac{166337705025}{2^{40}}
  & -\frac{1439205381561}{2^{44}}
  & -\frac{5834980110663}{2^{49}} \\[7pt]
5 & \frac{144358513071}{2^{40}}
  & \frac{92974300487}{2^{38}}
  & \frac{18708034100751}{2^{47}}
  & -\frac{543871644055407}{2^{53}} \\[7pt]
6 & -\frac{1059659119239}{2^{44}}
  & \frac{19563529502353}{2^{47}}
  & \frac{23525170472145}{2^{46}}
  & \frac{9973125866502423}{2^{56}} \\[7pt]
7 & \frac{3166836640983}{2^{49}}
  & -\frac{774586365338321}{2^{53}}
  & \frac{10507070361809241}{2^{56}}
  & \frac{1491032412033171}{2^{52}}
\end{array}
$}
\]

\newpage
Finally, the values of $f^{0,0}_{\infty}(r,s)$ are:

\[
\scalebox{1.2}{$
\begin{array}{c|cccc}
r \backslash s & 0 & 1 & 2 & 3 \\ \hline
0 & 0
  & \frac{59}{2^{10}}
  & \frac{1635}{2^{14}}
  & -\frac{19781}{2^{18}} \\[7pt]
1 & -\frac{59}{2^{10}}
  & 0
  & \frac{9347}{2^{17}}
  & -\frac{2495}{2^{22}} \\[7pt]
2 & -\frac{1635}{2^{14}}
  & -\frac{9347}{2^{17}}
  & 0
  & \frac{3156863}{2^{25}} \\[7pt]
3 & \frac{19781}{2^{18}}
  & \frac{2495}{2^{22}}
  & -\frac{3156863}{2^{25}}
  & 0 \\[3pt]
4 & \frac{17593}{2^{20}}
  & \frac{443573}{2^{25}}
  & -\frac{12113525}{2^{30}}
  & -\frac{644382097}{2^{33}} \\[7pt]
5 & -\frac{2526543}{2^{25}}
  & -\frac{16682267}{2^{30}}
  & \frac{61485319}{2^{30}}
  & -\frac{3828502499}{2^{37}} \\[7pt]
6 & \frac{57398423}{2^{30}}
  & \frac{16549373}{2^{31}}
  & -\frac{702165295}{2^{34}}
  & \frac{43988482319}{2^{40}} \\[7pt]
7 & \frac{276178165}{2^{34}}
  & \frac{1363019667}{2^{38}}
  & -\frac{24982339111}{2^{41}}
  & -\frac{248639212357}{2^{46}}
\end{array}
$}
\]
\\
\[
\scalebox{1.2}{$
\begin{array}{c|cccc}
r \backslash s & 4 & 5 & 6 & 7 \\ \hline
0 & -\frac{17593}{2^{20}}
  & \frac{2526543}{2^{25}}
  & -\frac{57398423}{2^{30}}
  & -\frac{276178165}{2^{34}} \\[7pt]
1 & -\frac{443573}{2^{25}}
  & \frac{16682267}{2^{30}}
  & -\frac{16549373}{2^{31}}
  & -\frac{1363019667}{2^{38}} \\[7pt]
2 & \frac{12113525}{2^{30}}
  & -\frac{61485319}{2^{30}}
  & \frac{702165295}{2^{34}}
  & \frac{24982339111}{2^{41}} \\[7pt]
3 & \frac{644382097}{2^{33}}
  & \frac{3828502499}{2^{37}}
  & -\frac{43988482319}{2^{40}}
  & \frac{248639212357}{2^{46}} \\[7pt]
4 & 0
  & \frac{86975044945}{2^{40}}
  & \frac{114016256935}{2^{44}}
  & -\frac{13201346218711}{2^{49}} \\[7pt]
5 & -\frac{86975044945}{2^{40}}
  & 0
  & \frac{11913341257839}{2^{47}}
  & \frac{104593352009425}{2^{53}} \\[7pt]
6 & -\frac{114016256935}{2^{44}}
  & -\frac{11913341257839}{2^{47}}
  & 0
  & \frac{5615451538169511}{2^{56}} \\[7pt]
7 & \frac{13201346218711}{2^{49}}
  & -\frac{104593352009425}{2^{53}}
  & -\frac{5615451538169511}{2^{56}}
  & 0
\end{array}
$}
\]

\clearpage

\bibliographystyle{alpha}
\bibliography{References}

\bigskip

\begin{singlespace}
\noindent
Melda Akyazi-Zavadski\\
Faculty of Mathematics, University of Vienna, Vienna, Austria\\
\textit{Email:} melda.akyazi-zavadski@univie.ac.at

\medskip

\noindent
Nathana\"el Berestycki\\
Faculty of Mathematics, University of Vienna, Vienna, Austria\\
\textit{Email:} nathanael.berestycki@univie.ac.at

\medskip

\noindent
Sunil Chhita\\
Department of Mathematical Sciences, Durham University, Durham, UK\\
\textit{Email:} sunil.chhita@durham.ac.uk
\end{singlespace}

\end{document}